\documentclass{article}
\usepackage{amsmath}
\usepackage{amssymb}

\begin{document}

\renewcommand{\theequation}{\arabic{section}.\arabic{equation}}%
\begin{equation*}
\end{equation*}

\begin{center}
{\Large Minimal time sliding mode control for evolution equations in Hilbert
spaces}

\bigskip

Gabriela Marinoschi

\medskip

\textquotedblleft Gheorghe Mihoc-Caius Iacob\textquotedblright Institute of
Mathematical Statistics and

Applied Mathematics of the Romanian Academy,

Calea 13 Septembrie 13, Bucharest, Romania

gabriela.marinoschi@acad.ro
\end{center}

\bigskip

{\small \noindent Abstract. This work is concerned with the time optimal
control problem for evolution equations in Hilbert spaces. The attention is
focused on the maximum principle for the time optimal controllers having the
dimension smaller that of the state system, in particular for minimal time
sliding mode controllers, which is one of the novelties of this paper. We
provide the characterization of the controllers by the optimality conditions
determined for some general cases. The proofs rely on a set of hypotheses
meant to cover a large class of applications. Examples of control problems
governed by parabolic equations with potential and drift terms, porous media
equation or reaction-diffusion systems with linear and nonlinear
perturbations, describing real world processes, are presented at the end.}

\medskip


\medskip


\section{Problem presentation}

\setcounter{equation}{0}

The purpose of this paper is to study the time optimal control for a family
of evolution equations in Hilbert spaces. In time optimal control the
optimality criterion is the elapsed time. Here, by the time optimal control
problem we mean to search for a constrained internal controller able to
drive the trajectory of the solution from an initial state to a given target
set in the shortest time, while controlling over the complete timespan.

Minimum time control problems have been initiated by Fattorini in the paper 
\cite{Fatt-64} and developed later in the monograph \cite{Fattorini}. A list
of only few titles dealing with this subject, in special for problems
governed by parabolic type equations includes \cite{Kunish-Wang-2012}, \cite%
{Kunish-Wang-2013}, \cite{Micu-Ro-Tu}, \cite{Wang-Wang-2003}, \cite%
{Wang-Wang-2007}, \cite{Wang-Zuazua}. In what concerns problems governed by
abstract evolution equations, we cite \cite{VB-SIAM-91} and the monographs 
\cite{VB-93}, \cite{VB-optim}, \cite{VB-control}. In \cite{VB-SIAM-91} the
existence and uniqueness of a viscosity solution was provided for the
Bellman equation associated with the time-optimal control problem for a
semilinear evolution equation in a Hilbert space, while in \cite{VB-97} the
time optimal control was studied for the Navier-Stokes equations. The
existence of the optimal time control for a phase-field system was proved in 
\cite{Wang-Wang-2003} for a regular double-well potential, by using the
Carleman inequality and the maximum principle was established by using two
controls acting in subsets of the space domain. The asymptotic behavior of
the solutions of a class of abstract parabolic time optimal control problems
when the generators converge, in an appropriate sense, to a given strictly
negative operator was studied in \cite{Tu-Wa-Wu-2016}. For a large class of
problems and aspects related to this subject we refer the reader to the
recent monograph \cite{Wang-book}.

From the perspective of applications, many processes in engineering,
physics, biology, medicine, environmental sciences, ecology require
solutions relying on time optimal control problems. The theoretical results
in this paper aim to cover models governed by parabolic equations with
potential and drift terms and various reaction-diffusion systems with linear
or nonlinear perturbations, or nonlocal control problems, presented in the
last section.

Especially of interest in applications is to control a system using a
controller whose dimension is smaller than that of the state system. In this
case the initial datum is steered not into a point, but into a linear
manifold of the state space, situation which is relevant for the sliding
mode control (for some references see e.g., \cite{BCGMR-sli}, \cite%
{sli-Cecilia}, \cite{CGMR-tu}). The solution to such a problem, which is
more challenging from the mathematical point of view, is a central point in
our theoretical approach.

We prove the existence of the time optimal control and the first order
necessary conditions of optimality in relation with the evolution equation
on a Hilbert space $H,$ 
\begin{equation}
y^{\prime }(t)+Ay(t)=Bu(t),\mbox{ }t>0,  \label{g4}
\end{equation}%
\begin{equation}
y(0)=y_{0}.  \label{g4-1}
\end{equation}%
Here, $A$ is a nonlinear and unbounded operator over a Hilbert space $H,$ $B$
is a linear operator from a Banach space $U$ to $H,$ $u$ is a controller
constrained to belong to $U$ and $y$ is the solution to (\ref{g4})
corresponding to the initial datum $y_{0}$ and controlled by $Bu.$ The
following minimization problem is studied: 
\begin{equation}
\mbox{Minimize }\left\{ J(T,u):=T;\mbox{ }(T,u)\in \mathcal{U}_{ad},\mbox{ }%
Py(T)=Py^{tar}\right\}  \tag{$\mathcal{P}$}
\end{equation}%
where $y$ is the solution to (\ref{g4})-(\ref{g4-1}) and 
\begin{equation}
\mathcal{U}_{ad}:=\{(T,u);\mbox{ }T\in \mathbb{R},\mbox{ }T>0,\mbox{ }u\in
L^{\infty }(0,T;U),\mbox{ }\left\Vert u(t)\right\Vert _{U}\leq \rho 
\mbox{
a.e. }t>0\}.  \label{Uad}
\end{equation}%
As it will be further explained, $P$ is an algebraic projection of the
solution from $H$ to $H$ or to a subspace of it, $y^{tar}$ is a prescribed
target for the state and $\rho $ is a positive constant at our choice.

Relying on certain hypotheses ensuring the well-posedness of (\ref{g4})-(\ref%
{g4-1}) and on the hypothesis of a not empty admissible set $\mathcal{U}%
_{ad},$ the existence of optimal controls is proved. A maximum principle is
first provided for an intermediate approximating problem. This generates a
sequence of approximating optimal solutions which converges to a precise
optimal pair to $(\mathcal{P}),$ whose characterization is a central point
of the paper.

The results are more relevant in the case of states with many components.
That is why, for the sake of a clearer explanation and for a simpler
notation, let us first assume that the state $y$ in (\ref{g4})-(\ref{g4-1})
has two components, $y=(y_{1},y_{2})\in H=H_{1}\times H_{2}$ and $%
u=(u_{1},u_{2})\in U=U_{1}\times U_{2},$ where $H_{i}$ are Hilbert spaces
and $U_{i}$ are Banach spaces, $i=1,2.$

We focus on two problems. The more challenging case is to steer, by the
action of only one control $Bu=(B_{1}u_{1},0)$, only the first component $%
y_{1}(t)$ of the state $y(t)$ from its initial value into a manifold $%
\mathcal{S},$ within a minimal time $T^{\ast }$. In this case, the target
manifold is $y_{1}=y_{1}^{\mbox{target}}.$ This action may be realized using
effectively one controller acting in the first equation. Thus, the state is
forced to reach the manifold $\mathcal{S}=\{y;$ $y_{1}=y_{1}^{\mbox{target}%
}\}$ on which it may continue to slide, for $t\geq T^{\ast },$ possibly
under supplementary conditions and by performing a controller slight
modification after the time $T^{\ast }$. This turns out to be in fact the
sliding mode control and it will be detailed for a reaction-diffusion model
in Section 6, Example 3.

Another possibility is to control both state components, forcing them to
reach a prescribed point target $y^{tar}:=(y_{1}^{\mbox{target}},y_{2}^{%
\mbox{target}}),$ by employing two controllers, with $%
Bu=(B_{1}u_{1},B_{2}u_{2}).$

Because the intention is to simultaneously prove the objectives stated
before, these are formalized by means of the minimization problem $(\mathcal{%
P})$ involving a mapping $P\in L(H,H)$ covering each of the following
situations:

$(i)$ $P(y_{1},y_{2})=(y_{1},y_{2})\in H_{1}\times H_{2},$ $%
B(u_{1},u_{2})=(B_{1}u_{1},B_{2}u_{2}),$ $y^{tar}=(y_{1}^{\mbox{target}%
},y_{2}^{\mbox{target}})=Py^{tar},$ in the case when both state components
are controlled by two controllers;

$(ii)$ $P(y_{1},y_{2})=(y_{1},0),$ $y_{1}\in H_{1},$ $y_{2}\in H_{2},$ $%
B(u_{1},u_{2})=(B_{1}u_{1},0),$ $Py^{tar}:=(y_{1}^{\mbox{target}},0),$ in
the situation when only the first component is controlled by one controller.

\noindent With this notation, (\ref{g4}) can be rewritten 
\begin{equation*}
y_{1}^{\prime }(t)+(Ay(t))_{1}=B_{1}u_{1}(t),\mbox{ }y_{2}^{\prime
}(t)+(Ay(t))_{2}=B_{2}u_{2}(t),\mbox{ a.e. }t>0
\end{equation*}%
in the case $(i)$ and 
\begin{equation*}
y_{1}^{\prime }(t)+(Ay(t))_{1}=B_{1}u_{1}(t),\mbox{ }y_{2}^{\prime
}(t)+(Ay(t))_{2}=0,\mbox{ a.e. }t>0
\end{equation*}%
in the case $(ii),$ where $Ay(t)=((Ay(t))_{1},(Ay(t))_{2}).$

Actually, in case $(ii),$ $P$ is an algebraic projection of $H$ into $%
H_{1}\times \{0\}\subset H$, mapping $y=(y_{1},y_{2})$ into its first
component $y_{1}$ which is the only one controlled in this case. For the
second component of $Py$ corresponding to the second equation which is not
controlled we set the value zero, in order to ensure a compact notation
consistent in the calculations (e.g., of the type $\left\Vert Pv\right\Vert
_{H}\leq C\left\Vert Bv\right\Vert _{H})$ with the value $Bu_{2}=0.$

In both cases we agree to use the same notation $y^{tar}$ in order to allow
a compact writing. In the first case, $y^{tar}$ contains the targets for
each state component. In the second case, the essential role is played by
the first component of $Py^{tar}$ while the second component of $y^{tar}$
plays no role. We can set the latter zero, even if this is not a target,
because its action, as well as that of the second component of $Py$ will be
cancelled in the calculations by the second zero component of $Bu$.

This explanation can be extended to the case with the state $y$ having $n$
components, either when all $n$ components are controlled by $n$
controllers, or when only $k$ trajectories $(y_{1},...,y_{k})$ are led into $%
(y_{1}^{\mbox{target}},...,y_{k}^{\mbox{target}}),$ by using $k$
controllers, via $Bu=(B_{1}u_{1},...B_{k}u_{k},0,...,0).$ We note that we
can change the notation, indicating the vector $(y_{1},...,y_{k})$ of the
first $k$ components still by $y_{1}$ and the vector $(y_{k+1},...,y_{n})$
by $y_{2}$ and we can use a similar notation for $u$ and $B.$ So, the
general case can be reduced to that with two state components. To conclude,
for the writing simplicity, we shall refer in the sequel to the case with
two state components.

The paper is organized as follows. The theoretical results rely on a set of
hypotheses, $(a_{1})-(a_{6}),$ $(b_{1}),$ $(c_{1}),$ listed in Section 2.
For the passing to the limit proof in Theorem 5.5, Section 5.3, there are
necessary some technical assumptions $(d_{1})-(d_{5}),$ including the
hypothesis (\ref{g5-000}). This is essential for the characterization of the
controller if only one state is controlled ($P\neq I,$ $Bu=(B_{1}u_{1},0)),$
which may be the more relevant in applications. This is one of the novelty
of this paper, besides the results characterizing the controller for
evolution equations with some general nonlinear operators. Section 3
includes some results of existence, beginning with the well-posedness of the
state system (\ref{g4})-(\ref{g4-1}), in Theorem 3.2. The existence of the
minimum time is provided in Theorem 3.3. In Section 4, we employ an
approximating problem $(\mathcal{P}_{\varepsilon })$ indexed along a small
parameter $\varepsilon $ occurring in some penalization terms of the
functional $J.$ After giving a basic result in Theorem 4.1 for the existence
of a solution to $(\mathcal{P}_{\varepsilon }),$ the convergence of $(%
\mathcal{P}_{\varepsilon })$ to $(\mathcal{P})$ is proved in Theorem 4.2.
This result is strong by asserting that if one fix an optimal pair $(T^{\ast
},u^{\ast })$ in $(\mathcal{P}),$ the sequence of optimal pairs in $(%
\mathcal{P}_{\varepsilon })$ tends exactly to $(T^{\ast },u^{\ast })$. The
necessary conditions of optimality for $(\mathcal{P}_{\varepsilon })$ are
determined in Proposition 5.4 at the end of an extremely technical
procedure, while in Theorem 5.5 the necessary conditions of optimality for $(%
\mathcal{P})$ are obtained as the limit of the previous ones, as $%
\varepsilon \rightarrow 0,$ after sharp estimates for the approximating
solution. A particular case for $U$, usually encountered, is treated in
Corollary 5.6. Applications of these results, including a detailed example
of minimum time sliding mode control, are presented in the last section. In
the Appendix we provide some definitions and general results necessary in
the paper.

\section{Functional framework and basic hypotheses}

\setcounter{equation}{0}

\paragraph{\textbf{Functional framework.}}

Let $V_{i}$ and $H_{i},$ $i=1,2,$ be Hilbert spaces and consider the
standard triplet $V_{i}\subset H_{i}\equiv H_{i}^{\ast }\subset V_{i}^{\ast
},$ with compact embeddings, where $V_{i}^{\ast }$ is the dual of $V_{i}.$
Let $U_{i},$ $i=1,2,$ be Banach spaces with the duals $U_{i}^{\ast }$
uniformly convex, implying that $U_{i}^{\ast }$ and $U_{i}$ are reflexive
(see e.g., \cite{vb-springer-2010}, p. 2). Let us denote%
\begin{equation*}
V=V_{1}\times V_{2},\mbox{ }H=H_{1}\times H_{2},\mbox{ }V^{\ast
}=V_{1}^{\ast }\times V_{2}^{\ast },\mbox{ }U=U_{1}\times U_{2},\mbox{ }%
U^{\ast }=U_{1}^{\ast }\times U_{2}^{\ast }.
\end{equation*}

\noindent We recall that the operator $P$ was defined in the introduction
(see $(i)-(ii)),$ but as a matter of fact we can define it on any space $%
X=X_{1}\times X_{2},$ where $X_{i}$ can be $V_{i},$ $H_{i},$ $V_{i}^{\ast }.$
Also, we use the same symbol $P$ for $X_{i}=U_{i}.$ Thus, 
\begin{equation*}
P:X=X_{1}\times X_{2}\rightarrow X,\mbox{ }P\in L(X,X),\mbox{ }
\end{equation*}%
and it is defined as%
\begin{equation*}
Py=(y_{1},y_{2})\mbox{ or }Py=(y_{1},0),\mbox{ for }y=(y_{1},y_{2})\in X.
\end{equation*}%
It can be easily seen $P^{2}=P$ and $\left\Vert Py\right\Vert _{X}\leq
\left\Vert y\right\Vert _{X}.$

Let $A:V\rightarrow V^{\ast }.$ We denote by $A_{H}:D(A_{H})\subset
H\rightarrow H$ the \textit{restriction} of $A$ on $H$ defined by $A_{H}y=Ay$
for $y\in D(A_{H})=\{y\in V;$ $Ay\in H\}.$ In the sequel, $\Gamma
:V\rightarrow V^{\ast }$ is the duality mapping from $V$ to $V^{\ast }$ and $%
\Gamma _{H}$ is the restriction of $\Gamma $ to $H$ (see \ref{Gamma1}).

Let $A\in C^{1}(V,V^{\ast }).$ The G\^{a}teaux derivative of $A$ is the
linear operator $A^{\prime }(y):V\rightarrow V^{\ast }$ defined by 
\begin{equation*}
A^{\prime }(y)z=\lim_{\lambda \rightarrow 0}\frac{A(y+\lambda z)-Ay}{\lambda 
}\mbox{ strongly in }V^{\ast },\mbox{ for all }y,z\in V
\end{equation*}%
and $A_{H}^{\prime }(y)$ is the restriction of $A^{\prime }(y)$ to $H$ (see (%
\ref{g0})-(\ref{g-00})). Properties of these operators are given in the
Appendix.

\medskip

\noindent \textbf{Notation. }Let $X$ and $Y$ be Banach spaces. By\textbf{\ }$%
L^{p}(0,T;X)$ we denote the space of $p$-summable functions from $(0,T)$ to $%
X$, for $1\leq p\leq \infty .$ $W^{1,p}(0,T;X)=\{f;$ $f:[0,T]\rightarrow X,$
absolutely continuous, $f,$ $df/dt\in L^{p}(0,T;X)\}.$ $C(X,Y)$ and $%
C^{1}(X,Y)$ are the spaces of continuous and differentiable G\^{a}teaux,
respectively, operators from $X$ to $Y$. $L(X,Y)$ is the space of linear
continuous operators from $X$ to $Y.$ We denote the scalar product and norm
in the space $X$ by $(\cdot ,\cdot )_{X}$ and $\left\Vert \cdot \right\Vert
_{X},$ respectively.

\noindent We shall denote by $C,$ $C_{i},$ $\alpha _{i},$ $\gamma _{i},$ $%
i=0,1,2,...$ positive constants that may change from line to line.

Some other notation and definitions related to the hypotheses below can be
found in the Appendix.

\medskip

\textbf{Hypotheses }$(a_{1}),(a_{2}),(b_{1}),(c_{1})$

$(a_{1})$ $A:V\rightarrow V^{\ast }$ is demicontinuous, $A0=0,$ 
\begin{equation}
\left\langle Ay-A\overline{y},y-\overline{y}\right\rangle _{V^{\ast },V}\geq
\alpha _{1}\left\Vert y-\overline{y}\right\Vert _{V}^{2}-\alpha
_{2}\left\Vert y-\overline{y}\right\Vert _{H}^{2},\mbox{ for all }y,%
\overline{y}\in V,\mbox{ }\alpha _{1}>0,  \label{g5}
\end{equation}%
\begin{equation}
A\mbox{ is bounded on bounded subsets of }V,  \label{g3}
\end{equation}%
\begin{equation}
D(A_{H})=D(\Gamma _{H}):=D_{H}.  \label{D0}
\end{equation}

$(a_{2})$%
\begin{equation}
(A_{H}y,\Gamma _{H}y)_{H}\geq \alpha _{3}\left\Vert \Gamma _{H}y\right\Vert
_{H}^{2}-\alpha _{4}\left\Vert y\right\Vert _{V}^{2},\mbox{ for all }y\in
D_{H}.  \label{A0H}
\end{equation}

$(b_{1})$%
\begin{equation*}
B\in L(U,H).
\end{equation*}

$(c_{1})$ For each $y^{tar}\in H,$ and $\rho $ positive large enough, there
exists\textbf{\ }$T_{\ast }>0$ and $u\in L^{\infty }(0,T;U)$ with $%
\sup\nolimits_{t\geq 0}\left\Vert u(t)\right\Vert _{U}\leq \rho ,$ such that 
$Py(T_{\ast })=Py^{tar},$ where $y$ is the solution to (\ref{g4})-(\ref{g4-1}%
).

\medskip

Hypotheses $(a_{1}),$ $(a_{2}),$ $(b_{1})$ are necessary to prove the state
system well-posedness and the existence of the solution to $(\mathcal{P})$,
in Section 3. The minimization problem $(\mathcal{P})$ is relevant if the
set $\mathcal{U}_{ad}$ is not empty. Hypothesis $(c_{1})$ ensures that $%
\mathcal{U}_{ad}\neq \varnothing $ and it is used in the proof of the
control existence. We specify that the proof of the controllability of (\ref%
{g4})-(\ref{g4-1}) is beyond the objective of this paper. However, for the
reader convenience, the existence of a least a pair $(T,u)$ in the
admissible set, or equivalently an example of proving the controllability of
(\ref{g4})-(\ref{g4-1}) in some cases, is given in Appendix, Proposition
7.1. Next, in the Examples, the reliability of $(c_{1})$ is commented in
each case.

\textbf{Hypotheses }$(a_{3})-(a_{6})$

$(a_{3})$ $A\in C^{1}(V;V^{\ast })$ and $A_{H}\in C^{1}(D_{H};H).$

$(a_{4})$ $A^{\prime }(y)$ and $A_{H}^{\prime }(y)$ defined by (\ref{g0})
and (\ref{g-00}) respectively, satisfy%
\begin{equation}
\left\Vert A^{\prime }(y)z\right\Vert _{V^{\ast }}\leq C\left\Vert
z\right\Vert _{V}(1+C_{1}\left\Vert y\right\Vert _{V}^{\kappa }),%
\mbox{ for
all }y,z\in V,\mbox{ some }\kappa \in \mathbb{R},\mbox{ }\kappa \geq 0,
\label{g7}
\end{equation}%
\begin{equation}
\left\Vert A_{H}^{\prime }(y)z\right\Vert _{H}\leq C\left\Vert z\right\Vert
_{D_{H}}(1+C_{1}\left\Vert y\right\Vert _{V}^{\kappa }),\mbox{ for all }%
y,z\in D_{H},\mbox{ }\kappa \geq 0.  \label{g7-1}
\end{equation}

$(a_{5})$ $A^{\prime }$ is strongly continuous from $V$ to $L_{s}(V,V^{\ast
}),$ (see (\ref{g-000})) and $A_{H}^{\prime }(y)$ is strongly continuous
from $V$ to $L(D_{H},H),$ namely%
\begin{eqnarray}
\left\Vert A_{H}^{\prime }(y_{n})\psi -A_{H}^{\prime }(y)\psi \right\Vert
_{H} &\rightarrow &0\mbox{ for all }\psi \in D_{H},\mbox{ }y_{n},y\in D_{H},
\label{g6} \\
\mbox{as }y_{n} &\rightarrow &y\mbox{ strongly in }V.  \notag
\end{eqnarray}

$(a_{6})$ The adjoint operator $(A^{\prime }(y))^{\ast }:V\rightarrow
V^{\ast }$ satisfies the condition%
\begin{equation}
\left\langle (A^{\prime }(y))^{\ast }z,\Gamma _{\nu }z\right\rangle
_{V^{\ast },V}\geq \gamma _{1}\left\Vert \Gamma _{\nu }z\right\Vert
_{H}^{2}-\gamma _{2}\left\Vert z\right\Vert _{V}^{2}(1+\gamma _{3}\left\Vert
y\right\Vert _{V}^{l})-\gamma _{4},\mbox{ }  \label{p-regular}
\end{equation}%
for all $y,z\in V,$ some $l\geq 0,$ $\gamma _{1},\gamma _{2}>0,$ where $%
\Gamma _{\nu },$ $\nu >0,$ is the Yosida approximation of $\Gamma ,$ see (%
\ref{Yosida}).

\medskip

Hypotheses $(a_{3})-(a_{6})$ are necessary for the proof of the existence of
the system in variations, the adjoint system and the determination of the
approximating optimality conditions.

\textbf{Hypotheses }$(d_{1})-(d_{5})$

Assume that $U\neq H$ and that there exists $\alpha \in (0,1)$ such that

$(d_{1})$ 
\begin{equation}
\left\Vert P\Gamma _{H}^{-\alpha /2}v\right\Vert _{H}\leq C\left\Vert
B^{\ast }v\right\Vert _{U^{\ast }},\mbox{ }v\in H\mbox{.}  \label{g74}
\end{equation}

$(d_{2})$ $(A_{H}^{\prime }(y))^{\ast }$ satisfies the relations 
\begin{equation}
\left( (A_{H}^{\prime }(y\right) )^{\ast }v,\Gamma _{H}^{-\alpha }v)_{H}\geq
C_{1}\left\Vert \Gamma _{H}^{(1-\alpha )/2}v\right\Vert
_{H}^{2}-C_{2}\left\Vert v\right\Vert _{H}^{2}(1+C_{3}\left\Vert
y\right\Vert _{V}^{l}),  \label{g74-1}
\end{equation}%
for all $y,v\in D_{H},$ some $l\geq 0,$ and%
\begin{equation}
\left( (A_{H}^{\prime }(y\right) )^{\ast }v,\Gamma _{H}^{-1}v)_{H}\geq
C_{1}\left\Vert v\right\Vert _{H}^{2}-C_{2}\left\Vert v\right\Vert _{V^{\ast
}}^{2}(1+C_{3}\left\Vert y\right\Vert _{V}^{l}),\mbox{ }  \label{g74-10}
\end{equation}

$(d_{3})$%
\begin{equation}
\left\Vert Pv\right\Vert _{V^{\ast }}\leq C^{\ast }\left\Vert B^{\ast
}v\right\Vert _{U^{\ast }},\mbox{ for }v\in H.  \label{g74-2}
\end{equation}

$(d_{4})$ Let $Py=(y_{1},0)$ and let $\rho $ be sufficiently large. For each 
$u\in L^{\infty }(0,\infty ;U),$ with $\left\Vert u(t)\right\Vert _{U}\leq
\rho $ there exists $\widehat{z},$ possibly depending on $u,$ such that 
\begin{equation}
\widehat{y}(t)=(y_{1}^{\mbox{target}},\widehat{z})  \label{y^}
\end{equation}%
satisfies%
\begin{equation}
\left\langle Ay(t)-A\widehat{y},P(y(t)-\widehat{y})\right\rangle _{V^{\ast
},V}\geq -C_{3}\left\Vert P(y(t)-\widehat{y})\right\Vert _{H}^{2},\mbox{ }
\label{g5-000}
\end{equation}%
for all $y\in V,$ and $t\in (0,T_{\ast }+\delta ),$ with the choice 
\begin{equation}
\rho >\rho _{1},\mbox{ }\rho _{1}:=C^{\ast }\left\Vert A_{H}\widehat{y}%
\right\Vert _{V}.  \label{Ay^}
\end{equation}

$(d_{5})$ Let $P=I,$ assume 
\begin{equation}
\rho >\rho _{1},\mbox{ }\rho _{1}:=C^{\ast }\left\Vert
A_{H}y^{tar}\right\Vert _{V},  \label{Aytar}
\end{equation}%
and relation (\ref{g5-000}), where $\widehat{y}$ is replaced by $%
y^{tar}=(y_{1}^{\mbox{target}},y_{2}^{\mbox{target}}).$

\medskip

We specify that $C^{\ast }$ in (\ref{Ay^}) and (\ref{Aytar}) is exactly the
constant $C^{\ast }$ occurring in (\ref{g74-2}), depending on the domain $%
\Omega $ and $B,$ $T_{\ast }$ is the time specified in the controllability
hypothesis $(c_{1}),$ $\delta $ is arbitrary and $y(t)$ is the solution to (%
\ref{g4})-(\ref{g4-1}) corresponding to $u.$

Assumption (\ref{g5-000}) is a basic statement in the proof of the
characterization of the controller in the case when only one state is
controlled by one controller. This is the case when the state is allowed to
reach a sliding manifold.

We also note that if $B_{1}=B_{2}=I$ or $B_{1}=I$ and $B_{2}=0$ and the
spaces are such that $V\subset U$ or $H\subset U,$ then (\ref{g74-2}) is
automatically satisfied. The case $U=H$ will be treated in Corollary 5.6.

\medskip

Immediate consequences of the previous hypotheses are:

The operator $\lambda I+A,$ for $\lambda $ positive large enough, is
coercive and $A_{H}$ is quasi $m$-accretive on $H\times H,$ implied by $%
(a_{1}).$

The operator $A^{\prime }(y)$ satisfies the estimate%
\begin{equation}
\left\langle A^{\prime }(y)z,z\right\rangle _{V^{\ast },V}\geq \alpha
_{1}\left\Vert z\right\Vert _{V}^{2}-\alpha _{2}\left\Vert y\right\Vert
_{H}^{2},\mbox{ for all }z\in V.  \label{g11}
\end{equation}

The operator $\left. A^{\prime }(y)\right\vert _{H}=A_{H}^{\prime }(y)$ is
quasi $m$-accretive for each $y\in V$ and%
\begin{equation}
(A_{H}^{\prime }(y)z,z)_{H}\geq \alpha _{1}\left\Vert z\right\Vert
_{V}^{2}-\alpha _{2}\left\Vert y\right\Vert _{H}^{2},%
\mbox{ for all
\thinspace }y,\mbox{ }z\in D_{H}.  \label{g11-1}
\end{equation}

By (\ref{g7}) we have 
\begin{equation}
\left\Vert A^{\prime }(y)\right\Vert _{L(V,V^{\ast })}\leq
C(1+C_{1}\left\Vert y\right\Vert _{V}^{\kappa }).  \label{g11-2}
\end{equation}

\medskip

\section{Existence results}

\setcounter{equation}{0}

In this section we provide the proofs of the existence of the solution to
the state system and of a solution to the minimization problem $(\mathcal{P}%
).$ All over in this section, we assume $(a_{1}),$ $(a_{2}),$ $(b_{1}),$ $%
(c_{1}).$ Let 
\begin{equation}
X_{T}=C([0,T];H)\cap W^{1,2}(0,T;H)\cap L^{\infty }(0,T;V)\cap
L^{2}(0,T;D_{H}).  \label{g14}
\end{equation}

\noindent \textbf{Definition 3.1. }A strong solution to the Cauchy problem (%
\ref{g4})-(\ref{g4-1}) is a continuous function $y:[0,T]\rightarrow H,$
which is a.e. differentiable and satisfies (\ref{g4}) a.e. $t\in (0,T)$ and (%
\ref{g4-1})$.$

\medskip

\noindent \textbf{Theorem 3.2.} \textit{Let }$T>0,$ $u\in L^{2}(0,T;U),$ $%
y_{0}\in V.$\textit{\ Then, }(\ref{g4})-(\ref{g4-1}) \textit{has a unique
strong solution }$y\in X_{T},$ \textit{satisfying} 
\begin{eqnarray}
&&\left\Vert y(t)\right\Vert _{V}^{2}+\int_{0}^{t}\left\Vert A_{H}y(\tau
)\right\Vert _{H}^{2}d\tau +\int_{0}^{t}\left\Vert y(\tau )\right\Vert
_{D_{H}}^{2}d\tau +\int_{0}^{t}\left\Vert y^{\prime }(\tau )\right\Vert
_{H}^{2}d\tau  \label{g15} \\
&\leq &C\left( \left\Vert y_{0}\right\Vert _{V}^{2}+\rho ^{2}T\right)
e^{CT}:=C_{T},\mbox{ \textit{for all} }t\in \lbrack 0,T],  \notag
\end{eqnarray}%
\textit{with }$C$\textit{\ a positive constant. Moreover, for two solutions }%
$y$\textit{\ and }$\overline{y}$\textit{\ corresponding to }$u$\textit{\ and 
}$\overline{u}$\textit{\ we have} 
\begin{equation}
\left\Vert (y-\overline{y})(t)\right\Vert _{H}^{2}+\int_{0}^{t}\left\Vert (y-%
\overline{y})(\tau )\right\Vert _{V}^{2}d\tau \leq C\left( \left\Vert y_{0}-%
\overline{y_{0}}\right\Vert _{H}^{2}+\left\Vert u-\overline{u}\right\Vert
_{L^{2}(0,T;U)}^{2}\right) ,\mbox{ \textit{for all} }t\in \lbrack 0,T].
\label{g15-0}
\end{equation}%
\textit{Finally, if }$u_{n}\in L^{\infty }(0,T;U),$\textit{\ }$%
u_{n}\rightarrow u$\textit{\ weak-star in }$L^{\infty }(0,T;U),$\textit{\
then the solution }$y_{n}$\textit{\ corresponding to }$u_{n}$\textit{\ tends
to }$y,$\textit{\ the solution corresponding to }$u,$ \textit{namely}%
\begin{eqnarray}
y_{n} &\rightarrow &y\mbox{ \textit{weakly in} }W^{1,2}(0,T;H)\cap
L^{2}(0,T;D_{H}),  \notag \\
&&\mbox{\textit{weak-star in} }L^{\infty }(0,T;V),%
\mbox{ \textit{strongly in}
}L^{2}(0,T;H),  \label{g15-1} \\
y_{n}(t) &\rightarrow &y(t)\mbox{ \textit{strongly in }}H,%
\mbox{
\textit{uniformly in} }\lbrack 0,T].  \notag
\end{eqnarray}

\medskip

\noindent \textbf{Proof. }We recall that\textbf{\ }$A_{H}$ is quasi $m$%
-accretive on $H\times H.$ Assume first that the right-hand side of (\ref{g4}%
), $f=Bu$ is in $W^{1,1}(0,T;H)$ and $y_{0}\in D_{H}=D(A_{H})$. In this case
we obtain a unique solution $y\in W^{1,\infty }(0,T;H)\cap L^{\infty
}(0,T;D_{H})$ (see e.g., \cite{vb-springer-2010}, p. 151, Theorem 4.9),
implying that $A_{H}y\in L^{\infty }(0,T;H).$ A\ first estimate is obtained
by testing the equation (\ref{g4}) by $y(t)$ and integrating it over $(0,t),$
\begin{eqnarray*}
&&\frac{1}{2}\left\Vert y(t)\right\Vert _{H}^{2}+\alpha
_{1}\int_{0}^{t}\left\Vert y(\tau )\right\Vert _{V}^{2}d\tau \\
&\leq &\frac{1}{2}\left\Vert y_{0}\right\Vert _{H}^{2}+\alpha
_{2}\int_{0}^{t}\left\Vert y(\tau )\right\Vert _{H}^{2}d\tau
+\int_{0}^{t}\left\Vert Bu(\tau )\right\Vert _{H}\left\Vert y(\tau
)\right\Vert _{H}d\tau ,
\end{eqnarray*}%
which yields 
\begin{equation}
\left\Vert y(t)\right\Vert _{H}^{2}+\int_{0}^{t}\left\Vert y(\tau
)\right\Vert _{V}^{2}d\tau \leq C\left( \left\Vert y_{0}\right\Vert
_{H}^{2}+\int_{0}^{T}\left\Vert u(t)\right\Vert _{U}^{2}dt\right) e^{Ct}.
\label{g16}
\end{equation}%
Then, we multiply (\ref{g4}) in $H$ by $\Gamma _{H}y(t),$ use (\ref{A0H})
and integrate over $(0,t),$ obtaining%
\begin{eqnarray*}
&&\left\Vert \Gamma _{H}^{1/2}y(t)\right\Vert _{H}^{2}+\alpha
_{3}\int_{0}^{t}\left\Vert \Gamma _{H}y(\tau )\right\Vert _{H}^{2}d\tau \\
&\leq &\left\Vert \Gamma _{H}^{1/2}y_{0}\right\Vert
_{H}^{2}+\int_{0}^{t}\left\Vert Bu(\tau )\right\Vert _{H}\left\Vert \Gamma
_{H}y(\tau )\right\Vert _{H}d\tau +\alpha _{4}\int_{0}^{t}\left\Vert y(\tau
)\right\Vert _{V}^{2}d\tau \\
&\leq &\left\Vert y_{0}\right\Vert _{V}^{2}+\frac{1}{2}\int_{0}^{t}\left%
\Vert \Gamma _{H}y(\tau )\right\Vert _{H}^{2}d\tau
+C_{1}\int_{0}^{T}\left\Vert u(t)\right\Vert
_{U}^{2}dt+C_{2}\int_{0}^{t}\left\Vert y(\tau )\right\Vert _{V}^{2}d\tau .
\end{eqnarray*}%
Using (\ref{g16}) we get 
\begin{equation}
\left\Vert y(t)\right\Vert _{V}^{2}+\int_{0}^{t}\left\Vert \Gamma _{H}y(\tau
)\right\Vert _{H}^{2}d\tau \leq C\left( \left\Vert y_{0}\right\Vert
_{V}^{2}+\int_{0}^{T}\left\Vert u(t)\right\Vert _{U}^{2}dt\right) e^{Ct}\leq
C_{T},\mbox{ }  \label{g15-10}
\end{equation}%
for all $t\in \lbrack 0,T],$ $T>0.$ We note that $C_{T}$ is continuous and
increasing with respect to $T,$ but it can vary from line to line via the
constant $C.$ This implies that 
\begin{equation}
\int_{0}^{T}\left\Vert y(\tau )\right\Vert _{D_{H}}^{2}d\tau \leq C_{T}%
\mbox{
and }\int_{0}^{T}\left\Vert A_{H}y(\tau )\right\Vert _{H}^{2}d\tau \leq
C_{T}.  \label{g16-0}
\end{equation}%
By comparison with (\ref{g4}) we deduce that $\int_{0}^{T}\left\Vert
y^{\prime }(\tau )\right\Vert _{H}^{2}d\tau \leq C_{T}.$ By gathering all
estimates we obtain (\ref{g15}).

To prove (\ref{g15-0}) we consider two solutions corresponding to $(y_{0},u)$
and $(\overline{y_{0}},\overline{u}),$ write the difference of the equations
for these solutions, test it by $(y-\overline{y})(t),$ integrate over $(0,t)$
and apply the Gronwall lemma.

We proceed further by a density argument. We take $y_{0}^{n}\in D_{H}$ and $%
u_{n}\in W^{1,1}(0,T;U)$ such that $y_{0}^{n}\rightarrow y_{0}$ strongly in $%
V$ and $u_{n}\rightarrow u$ strongly in $L^{2}(0,T;U),$ the latter implying $%
Bu_{n}\rightarrow Bu$ strongly in $L^{2}(0,T;H).$ It follows that the
solution to (\ref{g4}) with $Bu_{n}$ instead of $Bu$ and with the initial
datum $y_{0}^{n}$ has a unique strong solution $y_{n}\in W^{1,\infty
}(0,T;H)\cap L^{\infty }(0,T;D_{H}),$ satisfying (\ref{g15}) and (\ref{g15-0}%
). From here, it follows that $y_{n}\rightarrow y$ strongly in $%
C([0,T];H)\cap L^{2}(0,T;V),$ as $n\rightarrow \infty ,$ and the estimate (%
\ref{g15-0}) for $y_{n}$ is preserved at limit. The right-hand side of (\ref%
{g15}) is bounded and so $A_{H}y_{n}\rightarrow A_{H}y$ weakly in $%
L^{2}(0,T;H),$ since $A_{H}$ is strongly-weakly closed, and $y_{n}^{\prime
}\rightarrow y^{\prime }$ weakly in $L^{2}(0,T;H).$ The estimate (\ref{g15})
is preserved at limit by the lower weakly continuity of the norms.

Let $u_{n}\in L^{\infty }(0,T;U),$\textit{\ }$u_{n}\rightarrow u$\textit{\ }%
weak-star in $L^{\infty }(0,T;U).$\textit{\ }Then, (\ref{g4})-(\ref{g4-1})
has a unique solution $y_{n}$ satisfying (\ref{g15}). Since the estimates
are uniform, on a subsequence we get the convergences in the first line of (%
\ref{g15-1}). The strongly convergence follows by the Aubin-Lions lemma and
the last one by the Arzel\`{a}-Ascoli theorem. Passing to the limit in (\ref%
{g4})-(\ref{g4-1}) written for $y_{n}$ we get (\ref{g4})-(\ref{g4-1})
corresponding to $y.$ \hfill $\square $

\medskip

\noindent We observe that by (\ref{g15-1}) we deduce that\textbf{\ }$y\in
C_{w}([0,T];V),$ that is, except for a subset of zero measure, $y$ is a weak
continuous function from $[0,T]$ in $V.$ We recall that $y\in
C_{w}([0,T];V), $ if when $t_{n}\rightarrow t,$ as $n\rightarrow \infty ,$
it follows that $y(t_{n})\rightarrow y(t)$ weakly in $V.$ Indeed, in the
proof of Theorem 3.2 we have for $t_{n}\rightarrow t,$ that $%
y(t_{n})\rightarrow y(t)$ strongly in $H.$ On the other hand, $\left\Vert
y(t_{n})\right\Vert _{V}\leq C_{T}$ and so $y(t_{n})\rightarrow \xi $ weakly
in $V.$ But the limit is unique and so $\xi =y(t)\in V,$ for all $t\geq 0.$

Similarly, we deduce that $Ay\in C_{w}([0,T];V^{\ast }).$ By (\ref{g3}), $%
\left\Vert Ay(t_{n})\right\Vert _{V^{\ast }}$ is bounded, since $\left\Vert
y(t_{n})\right\Vert _{V}\leq C_{T}$ and so $Ay(t_{n})\rightarrow \zeta $
weakly in $V^{\ast }.$ On the other hand, $Ay(t_{n})\rightarrow Ay(t)$
strongly in \thinspace $D_{H}^{\ast },$ the dual of $D_{H}$ because $%
y(t_{n})\rightarrow y(t)$ in $H.$ Thus, $\zeta =Ay(t)\in V^{\ast },$ for all 
$t\geq 0.$

\medskip

Now we prove the existence of the minimum in $(\mathcal{P}).$ Recall that $%
y^{tar}:=(y_{1}^{\mbox{target}},y_{2}^{\mbox{target}})$ if $P=I$ and $%
y^{tar}:=(y_{1}^{\mbox{target}},z),$ $z\in H$ if $P\neq I.$

\medskip

\noindent \textbf{Theorem 3.3.} \textit{Let }%
\begin{equation*}
y_{0}\in V,\mbox{ }Py^{tar}\in P(H),\mbox{ }Py_{0}\neq Py^{tar}.
\end{equation*}%
\textit{Then, problem }$(\mathcal{P})$\textit{\ has at least one solution }$%
(T^{\ast },u^{\ast })$ \textit{with the corresponding state} $y^{T^{\ast
},u^{\ast }}:=y^{\ast }.$

\medskip

\noindent \textbf{Proof. }We recall that we have assumed\textbf{\ }$(c_{1})$
asserting that the admissible set $\mathcal{U}_{ad}\neq \varnothing .$ The
functional $J(T,u)=T$ is nonnegative, hence it has an infimum. We denote $%
\inf J(T,u)=T^{\ast }\geq 0.$ Let us consider a minimizing sequence $%
(T_{n},u_{n}),$ 
\begin{equation*}
T_{n}>0,\mbox{ }u_{n}\in L^{\infty }(0,\infty ;U),\mbox{ }\sup\limits_{t\in
(0,\infty )}\left\Vert u_{n}(t)\right\Vert _{U}\leq \rho ,\mbox{ }%
Py^{T_{n},u_{n}}(T_{n})=Py^{tar},
\end{equation*}%
where $y^{T_{n},u_{n}}$ is the solution to the state system corresponding to 
$(T_{n},u_{n}),$ such that%
\begin{equation}
T^{\ast }\leq J(T_{n},u_{n})=T_{n}\leq T^{\ast }+\frac{1}{n},\mbox{ }n\geq 1.
\label{g18}
\end{equation}%
On a subsequence it follows that 
\begin{equation}
u_{n}\rightarrow u^{\ast }\mbox{ weak-star in }L^{\infty }(0,T;U),\mbox{ }%
Bu_{n}\rightarrow Bu^{\ast }\mbox{ weak-star in }L^{\infty }(0,T;H),\mbox{
for all }T>0.  \label{g19}
\end{equation}%
We see that $T_{n}\rightarrow T^{\ast }$ and passing to the limit in (\ref%
{g18}) we get that $J(T^{\ast },u^{\ast })=T^{\ast }.$ The state system
corresponding to each $T>0$ and $u_{n}$ has a unique solution $y^{T,u_{n}}$
satisfying (\ref{g15}). In particular, this is true for $T=T_{n}$ and $%
T=T^{\ast }.$ We note that the restriction of the solution $y^{T_{n},u_{n}}$
to $(0,T^{\ast })$ is in fact the solution $y^{T^{\ast },u_{n}}.$ We have by
(\ref{g15}) 
\begin{equation}
\left\Vert y^{T^{\ast },u_{n}}\right\Vert _{X_{T^{\ast }}}^{2}+\left\Vert
A_{H}y^{T^{\ast },u_{n}}\right\Vert _{L^{2}(0,T^{\ast };H)}^{2}\leq
\left\Vert y^{T_{n},u_{n}}\right\Vert _{X_{T_{n}}}^{2}+\left\Vert
A_{H}y^{T_{n},u_{n}}\right\Vert _{L^{2}(0,T_{n};H)}^{2}\leq C_{T_{n}}\leq
C_{T^{\ast }+1},  \notag
\end{equation}%
where $C_{T}$ depends continuously and increasingly on $T$ (see (\ref{g15}%
)). Therefore, by selecting a subsequence and recalling (\ref{g15-1}) we
have 
\begin{eqnarray}
y^{T^{\ast },u_{n}} &\rightarrow &y^{\ast }:=y^{T^{\ast },u^{\ast }}%
\mbox{
weakly in }W^{1,2}(0,T^{\ast };H)\cap L^{2}(0,T^{\ast };D_{H}),  \notag \\
&&\mbox{weak-star in }L^{\infty }(0,T^{\ast };V)\mbox{ and strongly in }%
L^{2}(0,T^{\ast };H),  \label{g21-0}
\end{eqnarray}%
\begin{equation*}
A_{H}y^{T^{\ast },u_{n}}\rightarrow A_{H}y^{\ast }\mbox{ weakly in }%
L^{2}(0,T^{\ast };H),
\end{equation*}%
since\ $A_{H}$ is strongly-weakly closed. By Ascoli-Arzel\`{a} theorem we
still get%
\begin{equation}
y^{T^{\ast },u_{n}}(t)\rightarrow y^{\ast }(t)\mbox{ strongly in }H,%
\mbox{
uniformly in }t\in \lbrack 0,T^{\ast }].  \label{g21-1}
\end{equation}%
Also, by the last assertion in Theorem 3.2 we infer that $y^{\ast }$ is the
solution to the state system corresponding to $(T^{\ast },u^{\ast }).$ We
show next the convergence of $y^{T_{n},u_{n}}$ to $y^{\ast }$ as $%
n\rightarrow \infty .$ For any $v\in H$ we have 
\begin{eqnarray*}
&&\left\vert \int_{0}^{T_{n}}(y^{T_{n},u_{n}}(t),v)_{H}dt-\int_{0}^{T^{\ast
}}(y^{\ast }(t),v)_{H}dt\right\vert \leq \left\vert
\int_{0}^{T_{n}}(y^{T_{n},u_{n}}(t),v)_{H}dt-\int_{0}^{T^{\ast
}}(y^{T_{n},u_{n}}(t),v)_{H}dt\right\vert _{H} \\
&&+\left\vert \int_{0}^{T^{\ast
}}(y^{T_{n},u_{n}}(t),v)_{H}dt-\int_{0}^{T^{\ast }}(y^{\ast
}(t),v)_{H}dt\right\vert \\
&=&\left\vert \int_{T^{\ast
}}^{T_{n}}(y^{T_{n},u_{n}}(t),v)_{H}dt\right\vert +\left\vert
\int_{0}^{T^{\ast }}(y^{T_{n},u_{n}}(t)-y^{\ast }(t),v)_{H}dt\right\vert \\
&\leq &C\left\Vert y^{T_{n},u_{n}}\right\Vert _{L^{\infty
}(0,T_{n};H)}\left\vert T_{n}-T^{\ast }\right\vert +\left\vert
\int_{0}^{T^{\ast }}\left( y^{T^{\ast },u_{n}}(t)-y^{\ast }(t),v\right)
_{H}dt\right\vert \rightarrow 0,
\end{eqnarray*}%
because $\left\Vert y^{T_{n},u_{n}}(t)\right\Vert _{L^{\infty
}(0,T_{n};H)}\leq C_{T_{n}}\leq C_{T^{\ast }+1}$ and by (\ref{g15-1})$.$ We
took into account that $y^{T_{n},u_{n}}$ restricted to $(0,T^{\ast })$
coincides with $y^{T^{\ast },u_{n}.}$ In a similar way, we can prove the
weak convergences of the other sequences, that is $(y^{T_{n},u_{n}})^{\prime
}\rightarrow (y^{\ast })^{\prime },$ $A_{H}y^{T_{n},u_{n}}\rightarrow
A_{H}y^{\ast }$ weakly in $L^{2}(0,T^{\ast };H)$ and that $%
y^{T_{n},u_{n}}\rightarrow y^{\ast }$ strongly in $L^{2}(0,T^{\ast };H).$

It remains to prove that $Py^{\ast }(T^{\ast })=Py^{tar}.$ We have 
\begin{eqnarray*}
&&\left\Vert P(y^{T_{n},u_{n}}(T_{n})-y^{\ast }(T^{\ast }))\right\Vert
_{H}\leq \left\Vert P(y^{T_{n},u_{n}}(T_{n})-y^{T_{n},u_{n}}(T^{\ast
}))\right\Vert _{H}+\left\Vert P(y^{T_{n},u_{n}}(T^{\ast })-y^{\ast
}(T^{\ast }))\right\Vert _{H} \\
&=&\left\Vert \int_{T^{\ast }}^{T_{n}}(Py^{T_{n},u_{n}})^{\prime
}dt\right\Vert +\left\Vert P(y^{T^{\ast },u_{n}}(T^{\ast })-y^{\ast
}(T^{\ast }))\right\Vert _{H} \\
&\leq &\sqrt{T_{n}-T^{\ast }}\left\Vert P(y^{T_{n},u_{n}})^{\prime
}\right\Vert _{L^{2}(0,T_{n};H)}+\left\Vert P(y^{T^{\ast },u_{n}}(T^{\ast
})-y^{\ast }(T^{\ast }))\right\Vert _{H}
\end{eqnarray*}%
which tend to zero since $\left\Vert P(y^{T_{n},u_{n}})^{\prime }\right\Vert
_{L^{2}(0,T_{n};H)}\leq C_{T^{\ast }+1}$ and by (\ref{g15-1}). Hence 
\begin{equation*}
\lim\limits_{n\rightarrow \infty }Py^{T_{n},u_{n}}(T_{n})=Py^{tar}=Py^{\ast
}(T^{\ast }).
\end{equation*}%
From here, we also deduce that $T^{\ast }>0.$ Otherwise, we would have $%
Py^{\ast }(T^{\ast })=Py^{\ast }(0)=Py_{0},$ that is $Py_{0}=Py^{tar}$ which
contradicts the hypothesis that $Py_{0}\neq Py^{tar}.$ Thus, we have
obtained $T^{\ast }>0,$ $u^{\ast }$ with the restriction $\left\Vert u^{\ast
}\right\Vert _{L^{\infty }(0,T^{\ast };U)}\leq \rho $, and $J(T^{\ast
},u^{\ast })=T^{\ast }.$ We have got that $T^{\ast }$ is the unique infimum
time at which $Py^{\ast }(T^{\ast })=Py^{tar}.$ This ends the proof. \hfill $%
\square $

\section{The approximating problem}

\setcounter{equation}{0}

Let $\varepsilon $ be positive and consider the problem 
\begin{equation}
\mbox{Minimize }\left\{ J_{\varepsilon }(T,u)\right. \mbox{ }T>0,\mbox{ }%
u\in L^{\infty }(0,\infty ;U),\mbox{ }ess\sup_{t\in (0,\infty )}\left\Vert
u(t)\right\Vert _{U}\leq \rho \}, \tag{$\mathcal{P}_{\varepsilon }$}
\end{equation}%
subject to (\ref{g4})-(\ref{g4-1}), where 
\begin{eqnarray}
J_{\varepsilon }(T,u) &=&T+\frac{1}{2\varepsilon }\left\Vert
Py(T)-Py^{tar}\right\Vert _{H}^{2}  \label{Jeps} \\
&&\mbox{ \ \ }+\frac{\varepsilon }{2}\int_{0}^{T}\left\Vert Pu(t)\right\Vert
_{U}^{2}dt+\frac{1}{2}\int_{0}^{T}\left\Vert \int_{0}^{t}P(u(\tau )-u^{\ast
}(\tau ))d\tau \right\Vert _{U}^{2}dt.  \notag
\end{eqnarray}

\medskip

\noindent \textbf{Theorem 4.1. }\textit{Let }$y_{0}\in V,$ $Py^{tar}\in
P(H), $ $Py_{0}\neq Py^{tar}.$ \textit{Then, problem }$(\mathcal{P}%
_{\varepsilon }) $\textit{\ has at least a solution }$(T_{\varepsilon
}^{\ast },u_{\varepsilon }^{\ast })$\textit{, with the corresponding state }$%
y^{T_{\varepsilon }^{\ast },u_{\varepsilon }^{\ast }}:=y_{\varepsilon
}^{\ast }.$

\medskip

\noindent \textbf{Proof. }Since $J_{\varepsilon }(T,u)$ is nonnegative,
there exists $d_{\varepsilon }=\inf J_{\varepsilon }(T,u)$ and it is
positive. Indeed, we note that if $J_{\varepsilon }(T,u)=0,$ each term
should be equal with 0. This implies that in the second term of $%
J_{\varepsilon },$ $P(y(T=0))-Py^{tar}=0$ which is a contradiction with the
fact that $Py_{0}\neq Py^{tar}.$ We conclude that the optimal $%
T_{\varepsilon }^{\ast }$ must be positive.

\noindent We consider a minimizing sequence $(T_{\varepsilon
}^{n},u_{\varepsilon }^{n})$ with $T_{\varepsilon }^{n}>0$ and $\left\Vert
u_{\varepsilon }^{n}(t)\right\Vert _{U}\leq \rho ,$ satisfying 
\begin{equation}
d_{\varepsilon }\leq J_{\varepsilon }(T_{\varepsilon }^{n},u_{\varepsilon
}^{n})\leq d_{\varepsilon }+\frac{1}{n},\mbox{ }n\geq 1.  \label{g22}
\end{equation}%
Hence, $T_{\varepsilon }^{n}\rightarrow T_{\varepsilon }^{\ast },$ as $%
n\rightarrow \infty .$ Then, for any $\delta >0$ there exists $n_{\delta }$
such that $T_{\varepsilon }^{n}\geq T_{\varepsilon }^{\ast }-\delta ,$ with $%
\delta $ arbitrarily small, for $n\geq n_{\delta }.$ On a subsequence%
\begin{equation*}
u_{\varepsilon }^{n}\rightarrow u_{\varepsilon }^{\ast },\mbox{ }%
Pu_{\varepsilon }^{n}\rightarrow Pu_{\varepsilon }^{\ast }%
\mbox{ weak-star
in }L^{\infty }(0,T;U),\mbox{ }Bu_{\varepsilon }^{n}\rightarrow
Bu_{\varepsilon }^{\ast }\mbox{ weak-star in }L^{\infty }(0,T;H),
\end{equation*}%
for all $T>0.$ Then, the state system corresponding to any $T>d_{\varepsilon
}+1$ and $u_{\varepsilon }^{n}$ has a unique continuous solution satisfying (%
\ref{g15}) on $(0,T),$ and it tends, as $n\rightarrow \infty ,$ to the
solution corresponding to $(T,u_{\varepsilon }^{\ast }).$ In particular,
this happens for $T=T_{\varepsilon }^{\ast }-\delta $, with $\delta $
arbitrary small. Then, on a subsequence denoted still by $n,$ we have 
\begin{eqnarray}
&&y^{T_{\varepsilon }^{\ast }-\delta ,u_{\varepsilon }^{n}}\rightarrow
y^{T_{\varepsilon }^{\ast }-\delta ,u_{\varepsilon }^{\ast }}%
\mbox{ strongly
in }L^{2}(0,T_{\varepsilon }^{\ast }-\delta ;H),\mbox{ }  \label{g22-0} \\
&&\mbox{weakly in }W^{1,2}(0,T_{\varepsilon }^{\ast }-\delta ;H)\cap
L^{2}(0,T_{\varepsilon }^{\ast }-\delta ;D_{H}),\mbox{ and weak-star in }%
L^{\infty }(0,T_{\varepsilon }^{\ast }-\delta ;V),  \notag
\end{eqnarray}%
\begin{equation}
A_{H}y^{T_{\varepsilon }^{n}-\delta ,u_{\varepsilon }^{n}}\rightarrow
A_{H}y_{\varepsilon }^{\ast }\mbox{ weakly in }L^{2}(0,T_{\varepsilon
}^{\ast }-\delta ;H),  \label{g22-1}
\end{equation}%
\begin{equation}
y^{T_{\varepsilon }^{\ast }-\delta ,u_{\varepsilon }^{n}}(t)\rightarrow
y_{\varepsilon }^{\ast }(t)\mbox{ strongly in }H,\mbox{ uniformly in }t\in
\lbrack 0,T_{\varepsilon }^{\ast }-\delta ].  \label{g22-2}
\end{equation}%
Next, we proceed in a similar way as in Theorem 3.3 to show that $%
y^{T_{\varepsilon }^{n},u_{\varepsilon }^{n}}\rightarrow y^{T_{\varepsilon
}^{\ast },u_{\varepsilon }^{\ast }}:=y_{\varepsilon }^{\ast }$ in the
corresponding spaces and that%
\begin{equation}
Py^{T_{\varepsilon }^{n},u_{n}}(T_{\varepsilon }^{n})\rightarrow
Py^{T_{\varepsilon }^{\ast },u_{\varepsilon }^{\ast }}(T_{\varepsilon
}^{\ast })\mbox{ strongly in }H.  \label{g22-2-0}
\end{equation}%
These imply that $y_{\varepsilon }^{\ast }$ is the solution to the state
system corresponding to $(T_{\varepsilon }^{\ast },u_{\varepsilon }^{\ast
}). $

Let us denote 
\begin{equation*}
h_{\varepsilon }^{n}(t)=\int_{0}^{t}P(u_{\varepsilon }^{n}-u^{\ast })(\tau
)d\tau ,\mbox{ }t\geq 0.
\end{equation*}%
Taking $\psi \in U^{\ast }$ we have 
\begin{eqnarray*}
\left\langle \int_{0}^{t}P(u_{\varepsilon }^{n}-u^{\ast })(\tau )d\tau ,\psi
\right\rangle _{U,U^{\ast }} &=&\int_{0}^{t}\left\langle P(u_{\varepsilon
}^{n}-u^{\ast })(\tau ),\psi \right\rangle _{U,U^{\ast }}d\tau \rightarrow \\
\int_{0}^{t}\left\langle P(u_{\varepsilon }^{\ast }-u^{\ast })(\tau ),\psi
\right\rangle _{U,U^{\ast }}d\tau &=&\left\langle
\int_{0}^{t}P(u_{\varepsilon }^{\ast }-u^{\ast })(\tau )d\tau ,\psi
\right\rangle _{U,U^{\ast }},\mbox{ for all }t\geq 0.
\end{eqnarray*}%
So, we get that $h_{\varepsilon }^{n}(t)\rightarrow h_{\varepsilon }^{\ast
}(t)$ weakly in $U,$ for all $t\geq 0$ and 
\begin{equation}
h_{\varepsilon }^{\ast }(t)=\int_{0}^{t}P(u_{\varepsilon }^{\ast }-u^{\ast
})(\tau )d\tau .  \label{g22-4}
\end{equation}%
Passing to the limit in (\ref{g22}) we get on the basis of the previous
convergences and of the weakly lower semicontinuity of the norms, that $%
J_{\varepsilon }(T_{\varepsilon }^{\ast },u_{\varepsilon }^{\ast
})=d_{\varepsilon },$ that is $(T_{\varepsilon }^{\ast },u_{\varepsilon
}^{\ast })$ is an optimal controller in $(P_{\varepsilon }).$ \hfill $%
\square $

\medskip

\noindent \textbf{Theorem 4.2. }\textit{Assume }$y_{0}\in V,$ $Py^{tar}\in
P(H),$ $Py_{0}\neq Py^{tar}$\textit{. Let }$(T^{\ast },u^{\ast },y^{\ast })$ 
\textit{be optimal in }$(\mathcal{P})$\textit{\ and }$(T_{\varepsilon
}^{\ast },u_{\varepsilon }^{\ast },y_{\varepsilon }^{\ast })$\textit{\ be
optimal in }$(\mathcal{P}_{\varepsilon }).$ \textit{Then, }%
\begin{equation}
T_{\varepsilon }^{\ast }\rightarrow T^{\ast },\mbox{ }u_{\varepsilon }^{\ast
}\rightarrow u^{\ast }\mbox{ \textit{weak-star in} }L^{\infty }(0,T^{\ast
};U),\mbox{ }Bu_{\varepsilon }^{\ast }\rightarrow Bu^{\ast }%
\mbox{
\textit{weak-star in} }L^{\infty }(0,T^{\ast };H),  \label{g25}
\end{equation}%
\begin{eqnarray}
y_{\varepsilon }^{\ast } &\rightarrow &y^{\ast }\mbox{ \textit{strongly in} }%
L^{2}(0,T^{\ast };H),\mbox{ }  \label{g26} \\
&&\mbox{\textit{weakly in} }W^{1,2}(0,T^{\ast };H)%
\mbox{ \textit{and
weak-star in} }L^{\infty }(0,T^{\ast };V),  \notag
\end{eqnarray}%
\begin{equation}
Ay_{\varepsilon }^{\ast }\rightarrow Ay^{\ast }\mbox{ \textit{weakly in} }%
L^{2}(0,T^{\ast };H),  \label{g27}
\end{equation}%
\begin{equation}
y_{\varepsilon }^{\ast }(T^{\ast })\rightarrow y^{\ast }(T^{\ast })%
\mbox{ 
\textit{strongly in} }H,\mbox{ }Py^{\ast }(T^{\ast })=P^{tar}.  \label{g29}
\end{equation}

\medskip

\noindent \textbf{Proof. }Let $(T_{\varepsilon }^{\ast },u_{\varepsilon
}^{\ast },y_{\varepsilon }^{\ast })$\textit{\ }be optimal in $%
(P_{\varepsilon }).$ Then, 
\begin{eqnarray}
&&J_{\varepsilon }(T_{\varepsilon }^{\ast },u_{\varepsilon }^{\ast
})=T_{\varepsilon }^{\ast }+\frac{1}{2\varepsilon }\left\Vert
Py_{\varepsilon }^{\ast }(T_{\varepsilon }^{\ast })-Py^{tar}\right\Vert
_{H}^{2}  \label{g30} \\
&&+\frac{\varepsilon }{2}\int_{0}^{T_{\varepsilon }^{\ast }}\left\Vert
Pu_{\varepsilon }^{\ast }(t)\right\Vert _{U}^{2}dt+\frac{1}{2}%
\int_{0}^{T_{\varepsilon }^{\ast }}\left\Vert \int_{0}^{t}P(u_{\varepsilon
}^{\ast }(\tau )-u^{\ast }(\tau ))d\tau \right\Vert _{U}^{2}dt  \notag \\
&\leq &J_{\varepsilon }(T,u)=T+\frac{1}{2\varepsilon }\left\Vert
Py(T)-Py^{tar}\right\Vert _{H}^{2}  \notag \\
&&+\frac{\varepsilon }{2}\int_{0}^{T}\left\Vert Pu(t)\right\Vert _{U}^{2}dt+%
\frac{1}{2}\int_{0}^{T}\left\Vert \int_{0}^{t}P(u(\tau )-u^{\ast }(\tau
))d\tau \right\Vert _{U}^{2}dt,  \notag
\end{eqnarray}%
for any $T>0$ and $u\in L^{\infty }(0,\infty ;U),$ $\left\Vert
u(t)\right\Vert _{U}\leq \rho $ a.e. $t>0,$ where $y_{\varepsilon }^{\ast }$
is the solution to the state system corresponding to $(T_{\varepsilon
}^{\ast },u_{\varepsilon }^{\ast })$ and $y$ is the solution to the state
system corresponding to $(T,u).$ Let us set in (\ref{g30}), $T=T^{\ast }$
and $u=u^{\ast },$ an optimal controller in $(P).$ Thus, the second and the
last term on the right-hand side of (\ref{g30}) are zero and 
\begin{eqnarray}
&&J_{\varepsilon }(T_{\varepsilon }^{\ast },u_{\varepsilon }^{\ast
})=T_{\varepsilon }^{\ast }+\frac{1}{2\varepsilon }\left\Vert
Py_{\varepsilon }^{\ast }(T_{\varepsilon }^{\ast })-Py^{tar}\right\Vert
_{H}^{2}+\frac{\varepsilon }{2}\int_{0}^{T_{\varepsilon }^{\ast }}\left\Vert
Pu_{\varepsilon }^{\ast }(t)\right\Vert _{U}^{2}dt  \label{g31} \\
&&+\frac{1}{2}\int_{0}^{T_{\varepsilon }^{\ast }}\left\Vert
\int_{0}^{t}P(u_{\varepsilon }^{\ast }(\tau )-u^{\ast }(\tau ))d\tau
\right\Vert _{U}^{2}dt\leq T^{\ast }+\frac{\varepsilon }{2}\int_{0}^{T^{\ast
}}\left\Vert Pu^{\ast }(t)\right\Vert _{U}^{2}dt.  \notag
\end{eqnarray}%
Then, $T_{\varepsilon }^{\ast }\rightarrow T^{\ast \ast }$ and $%
T_{\varepsilon }^{\ast }\geq T^{\ast \ast }-\delta ,$ with $\delta $
arbitrarily small, and selecting a subsequence indicated still by $%
\varepsilon ,$ we have%
\begin{eqnarray*}
u_{\varepsilon }^{\ast } &\rightarrow &u^{\ast \ast }\mbox{ weak-star in }%
L^{\infty }(0,T;U),\mbox{ }\left\Vert u^{\ast \ast }(t)\right\Vert _{U}\leq
\rho , \\
Bu_{\varepsilon }^{\ast } &\rightarrow &Bu^{\ast \ast }\mbox{ weak-star in }%
L^{2}(0,T;H),\mbox{ for all }T>0.
\end{eqnarray*}%
The solution $y_{\varepsilon }^{\ast }$ satisfies the estimates 
\begin{equation}
\left\Vert y_{\varepsilon }^{\ast }\right\Vert _{X_{T_{\varepsilon }^{\ast
}}}^{2}+\left\Vert A_{H}y_{\varepsilon }^{\ast }\right\Vert
_{L^{2}(0,T_{\varepsilon }^{\ast };H)}^{2}\leq C_{T_{\varepsilon }^{\ast }},
\label{g31-1}
\end{equation}%
where $X_{T}:=C([0,T],H)\cap L^{\infty }(0,T;V)\cap W^{1,2}(0,T;H)$, and $%
y_{\varepsilon }^{\ast }\rightarrow y^{\ast \ast }$ in the spaces $%
X_{T^{\ast \ast }-\delta }$ defined on $(0,T^{\ast \ast }-\delta ),$ for $%
\delta $ arbitrary$.$ As in the previous proof we show that all the
convergences (\ref{g22-0})-(\ref{g22-2}) take place also in the spaces
defined on $(0,T^{\ast \ast }).$ Also, we have%
\begin{equation*}
\left\Vert y_{\varepsilon }^{\ast }(T_{\varepsilon }^{\ast })-y^{\ast \ast
}(T^{\ast \ast }-\delta )\right\Vert _{H}\leq \sqrt{T_{\varepsilon }^{\ast
}-T^{\ast \ast }+\delta }\left\Vert (y_{\varepsilon }^{\ast })^{\prime
}\right\Vert _{L^{2}(0,T_{\varepsilon }^{\ast };H)}+\left\Vert
y_{\varepsilon }^{\ast }(T_{\varepsilon }^{\ast }-\delta )-y^{\ast \ast
}(T^{\ast \ast }-\delta )\right\Vert _{H},\mbox{ }
\end{equation*}%
implying that $\left\Vert y_{\varepsilon }^{\ast }(T_{\varepsilon }^{\ast
})-y^{\ast \ast }(T^{\ast \ast })\right\Vert _{H}\rightarrow 0,$ strongly in 
$H,$ as $\varepsilon \rightarrow 0.$ By (\ref{g31}) we have 
\begin{equation*}
\left\Vert Py_{\varepsilon }^{\ast }(T_{\varepsilon }^{\ast
})-Py^{tar}\right\Vert _{H}^{2}\leq 2\varepsilon T^{\ast }+\varepsilon
^{2}\int_{0}^{T^{\ast }}\left\Vert Pu^{\ast }(t)\right\Vert _{U}^{2}dt
\end{equation*}%
and so $Py_{\varepsilon }^{\ast }(T_{\varepsilon }^{\ast })\rightarrow
Py^{tar}$ strongly in $H,$ implying the relation $Py^{\ast \ast }(T^{\ast
\ast })=Py^{tar}.$ Again by (\ref{g31}), 
\begin{equation*}
T_{\varepsilon }^{\ast }\leq J_{\varepsilon }(T_{\varepsilon }^{\ast
},u_{\varepsilon }^{\ast })\leq T^{\ast }+\frac{\varepsilon }{2}%
\int_{0}^{T^{\ast }}\left\Vert Pu^{\ast }(t)\right\Vert _{U}^{2}dt,
\end{equation*}%
whence we get at limit that $T^{\ast \ast }\leq T^{\ast }.$ Now $T^{\ast
\ast }$ and $u^{\ast \ast }$ satisfy the restrictions required in problem $%
(P),$ that is $T^{\ast \ast }>0,$ $\left\Vert u^{\ast \ast }(t)\right\Vert
_{U}\leq \rho ,$ and $Py^{\ast \ast }(T^{\ast \ast })=Py^{tar},$ and since $%
T^{\ast }$ is the infimum in $(P)$ it follows that $T^{\ast \ast }=T^{\ast }.
$ Recalling (\ref{g22-4}) we define%
\begin{equation*}
E_{\varepsilon }(T):=\int_{0}^{T}\left\Vert \int_{0}^{t}P(u_{\varepsilon
}^{\ast }-u^{\ast })(\tau )d\tau \right\Vert
_{U}^{2}dt=\int_{0}^{T}\left\Vert h_{\varepsilon }^{\ast }(t)\right\Vert
_{U}^{2}dt,\mbox{ for all }T>0.
\end{equation*}%
We have by (\ref{g31}) that 
\begin{equation*}
T_{\varepsilon }^{\ast }+E_{\varepsilon }(T_{\varepsilon }^{\ast })\leq
J_{\varepsilon }(T_{\varepsilon }^{\ast },u_{\varepsilon }^{\ast })\leq
T^{\ast }+\frac{\varepsilon }{2}\int_{0}^{T^{\ast }}\left\Vert Pu^{\ast
}(t)\right\Vert _{U}^{2}dt
\end{equation*}%
and so $T^{\ast \ast }+\limsup\limits_{\varepsilon \rightarrow
0}E_{\varepsilon }(T_{\varepsilon }^{\ast })\leq T^{\ast },$ hence 
\begin{equation}
\limsup\limits_{\varepsilon \rightarrow 0}\int_{0}^{T_{\varepsilon }^{\ast
}}\left\Vert h_{\varepsilon }^{\ast }(t)\right\Vert
_{U}^{2}dt=\limsup\limits_{\varepsilon \rightarrow
0}\int_{0}^{T_{\varepsilon }^{\ast }}\left\Vert \int_{0}^{t}P(u_{\varepsilon
}^{\ast }-u^{\ast })(\tau )d\tau \right\Vert _{U}^{2}dt=0.  \label{g33}
\end{equation}%
Therefore, 
\begin{equation}
h_{\varepsilon }^{\ast }\rightarrow 0\mbox{ strongly in }L^{2}(0,T^{\ast
};U),\mbox{ as }\varepsilon \rightarrow 0.  \label{g34}
\end{equation}%
On the other hand, we know that $u_{\varepsilon }^{\ast }\rightarrow u^{\ast
\ast }$ weakly in $L^{2}(0,T^{\ast };U),$ so that 
\begin{equation*}
\int_{0}^{t}\left\langle P(u^{\ast \ast }-u^{\ast })(\tau ),\psi
\right\rangle _{U,U^{\ast }}d\tau =0,\mbox{ for all }t\in (0,T^{\ast }),%
\mbox{ }\psi \in U^{\ast },
\end{equation*}%
implying that $u^{\ast \ast }=u^{\ast }$ on $(0,T^{\ast }).$

For a later use we prove that 
\begin{equation}
h_{\varepsilon }^{\ast }(T_{\varepsilon }^{\ast })\rightarrow 0%
\mbox{
strongly in }U.  \label{g34-0}
\end{equation}%
We write 
\begin{equation*}
h_{\varepsilon }^{\ast }(T_{\varepsilon }^{\ast })-h_{\varepsilon }^{\ast
}(t)=\int_{t}^{T_{\varepsilon }^{\ast }}P(u_{\varepsilon }^{\ast }-u^{\ast
})(s)ds,\mbox{ for all }t\in \lbrack 0,T_{\varepsilon }^{\ast }].
\end{equation*}%
Then,%
\begin{equation*}
\left\Vert h_{\varepsilon }^{\ast }(T_{\varepsilon }^{\ast })\right\Vert
_{U}\leq \left\Vert h_{\varepsilon }^{\ast }(t)\right\Vert
_{U}+\int_{t}^{T_{\varepsilon }^{\ast }}\left\Vert P(u_{\varepsilon }^{\ast
}-u^{\ast })(s)\right\Vert _{U}ds\leq \left\Vert h_{\varepsilon }^{\ast
}(t)\right\Vert _{U}+2\rho (T_{\varepsilon }^{\ast }-t).
\end{equation*}%
Let us take $t\in \lbrack T_{\varepsilon }^{\ast }-\tau ,T_{\varepsilon
}^{\ast }]$ with $\tau >\varepsilon ,$ and integrate the previous inequality
along with $t$ in this interval. We have%
\begin{equation*}
\tau \left\Vert h_{\varepsilon }^{\ast }(T_{\varepsilon }^{\ast
})\right\Vert _{U}\leq \int_{T_{\varepsilon }^{\ast }-\tau }^{T_{\varepsilon
}^{\ast }}\left\Vert h_{\varepsilon }^{\ast }(t)\right\Vert _{U}dt+2\rho
\tau ^{2}\leq \int_{0}^{T_{\varepsilon }^{\ast }}\left\Vert h_{\varepsilon
}^{\ast }(t)\right\Vert _{U}dt+2\rho \tau ^{2}\leq \sqrt{T_{\varepsilon
}^{\ast }}\left( \int_{0}^{T_{\varepsilon }^{\ast }}\left\Vert
h_{\varepsilon }^{\ast }(t)\right\Vert _{U}^{2}dt\right) ^{1/2}+2\rho \tau
^{2}.
\end{equation*}%
Let us make $\varepsilon $ goes to zero and get by (\ref{g33}) that%
\begin{equation*}
\limsup\limits_{\varepsilon \rightarrow 0}\left\Vert h_{\varepsilon }^{\ast
}(T_{\varepsilon }^{\ast })\right\Vert _{U}\leq 2\rho \tau ,
\end{equation*}%
which yields (\ref{g34-0}), since $\tau $ is arbitrary. On the basis of (\ref%
{g31}) we write that 
\begin{equation*}
T_{\varepsilon }^{\ast }+\frac{1}{2\varepsilon }\left\Vert Py_{\varepsilon
}^{\ast }(T_{\varepsilon }^{\ast })-Py^{tar}\right\Vert _{H}^{2}\leq
J_{\varepsilon }(T_{\varepsilon }^{\ast },u_{\varepsilon }^{\ast })\leq
T^{\ast }+\frac{\varepsilon }{2}\int_{0}^{T^{\ast }}\left\Vert Pu^{\ast
}(t)\right\Vert _{U}^{2}dt
\end{equation*}%
whence 
\begin{equation}
\limsup_{\varepsilon \rightarrow 0}\frac{1}{2\varepsilon }\left\Vert
Py_{\varepsilon }^{\ast }(T_{\varepsilon }^{\ast })-Py^{tar}\right\Vert
_{H}^{2}=0.  \label{g36-0}
\end{equation}%
We conclude that $\lim\limits_{\varepsilon \rightarrow 0}J_{\varepsilon
}(T_{\varepsilon }^{\ast \ast },u_{\varepsilon }^{\ast \ast })=T^{\ast \ast
}=J(T^{\ast \ast },u^{\ast \ast })$ and so $(T^{\ast \ast },u^{\ast \ast })$
is optimal in $(P)$. But, $T^{\ast }$ is also optimal and unique and it
follows $T^{\ast \ast }=T^{\ast }$ and $u^{\ast \ast }=u^{\ast }$ a.e. on $%
(0,T^{\ast }).$ Eventually, we also have obtained (\ref{g26})-(\ref{g29}),
as claimed.\hfill $\square $

\section{The maximum principle}

\setcounter{equation}{0}

In this section, besides $(a_{1}),$ $(a_{2}),$ $(b_{1}),$ $(c_{1}),$ we
assume $(a_{3})-(a_{6}).$

\subsection{The system of first order variations and the dual system}

Let us introduce the Cauchy problem 
\begin{eqnarray}
Y^{\prime }(t)+A^{\prime }(y_{\varepsilon }^{\ast }(t))Y(t) &=&Bv(t),%
\mbox{
a.e. }t>0,  \label{g47} \\
Y(0) &=&0,  \notag
\end{eqnarray}%
where $v\in L^{\infty }(0,\infty ;U),$ $\left\Vert v(t)\right\Vert _{U}\leq
C,$ a.e. $t\geq 0.$

\medskip

\noindent \textbf{Proposition 5.1. }\textit{Problem} (\ref{g47}) \textit{has
a unique solution} 
\begin{equation}
Y\in C([0,T];H)\cap W^{1,2}(0,T;V^{\ast })\cap L^{2}(0,T;V),%
\mbox{\textit{\
for all }}T>0.  \label{g47-1}
\end{equation}

\medskip

\noindent \textbf{Proof}. We recall that $A^{\prime }(y_{\varepsilon }^{\ast
}(t))$ is continuous from $V$ to $V^{\ast },$ for all $t\geq 0,$ and has the
properties (\ref{g11}) and (\ref{g7}), $\left\Vert A^{\prime
}(y_{\varepsilon }^{\ast }(t))z\right\Vert _{V^{\ast }}\leq C\left\Vert
z\right\Vert _{V},$ due to $\left\Vert y_{\varepsilon }^{\ast
}(t)\right\Vert _{V}\leq C_{T},$ by (\ref{g15}). Then, the result claimed in
the statement\ is ensured by the Lions theorem. \hfill $\square $

\medskip

\noindent Let $(T_{\varepsilon }^{\ast },u_{\varepsilon }^{\ast })$ be an
optimal controller in $(\mathcal{P}_{\varepsilon }).$ For $\lambda >0,$ we
set%
\begin{equation}
u_{\varepsilon }^{\lambda }=Pu_{\varepsilon }^{\ast }+\lambda v,%
\mbox{ where 
}v=P(\overline{u}-u_{\varepsilon }^{\ast }),\mbox{ }\left\Vert \overline{u}%
\right\Vert _{L^{\infty }(0,\infty ;U)}\leq \rho .  \label{u-lambda}
\end{equation}%
In this way we can give variations to both controllers, if $P=I,$ or to the
first component in the case $P\neq I.$ We define $Y_{\lambda }=\frac{%
y_{\varepsilon }^{\lambda }-y_{\varepsilon }^{\ast }}{\lambda },$ where $%
y_{\varepsilon }^{\lambda }$ is the solution to the state system (\ref{g4})-(%
\ref{g4-1}) corresponding to $u_{\varepsilon }^{\lambda }$ and $%
T_{\varepsilon }^{\ast }.$

\medskip

\noindent \textbf{Proposition 5.2. }\textit{Let }$Y$\textit{\ be the
solution to} (\ref{g47}) \textit{and let} $T>0$. \textit{We have }%
\begin{equation}
\lim_{\lambda \rightarrow 0}\frac{y_{\varepsilon }^{\lambda }-y_{\varepsilon
}^{\ast }}{\lambda }=Y\mbox{ \textit{strongly in} }C([0,T];H)\cap
L^{2}(0,T;V),\mbox{ \textit{as} }\lambda \rightarrow 0,  \label{g47-4}
\end{equation}%
\textit{which ensures that} (\ref{g47}) \textit{is just the system of first
order variations related to} (\ref{g4})-(\ref{g4-1}).

\medskip

\noindent \textbf{Proof. }Let us define 
\begin{equation*}
\zeta _{\lambda }=\frac{y_{\varepsilon }^{\lambda }-y_{\varepsilon }^{\ast }%
}{\lambda }-Y.
\end{equation*}%
We write the equation for $y_{\varepsilon }^{\lambda }$ subtract the
equation for $y_{\varepsilon }^{\ast },$ divide by $\lambda $ and subtract
the equation (\ref{g47}). The equation verified by $\zeta _{\lambda }$ reads%
\begin{eqnarray}
\zeta _{\lambda }^{\prime }+A^{\prime }(y_{\varepsilon }^{\ast })\zeta
_{\lambda }+\frac{Ay_{\varepsilon }^{\lambda }-Ay_{\varepsilon }^{\ast }}{%
\lambda }-A^{\prime }(y_{\varepsilon }^{\ast })\frac{y_{\varepsilon
}^{\lambda }-y_{\varepsilon }^{\ast }}{\lambda } &=&0,  \label{g47-2} \\
\zeta _{\lambda }(0) &=&0.  \notag
\end{eqnarray}%
Now, we can represent the third term as 
\begin{equation}
Ay_{\varepsilon }^{\lambda }-Ay_{\varepsilon }^{\ast
}=\int_{0}^{1}(A^{\prime }(\nu y_{\varepsilon }^{\lambda }+(1-\nu
)y_{\varepsilon }^{\ast })(y_{\varepsilon }^{\lambda }-y_{\varepsilon
}^{\ast })d\nu   \label{g47-2-0}
\end{equation}%
and so, the equation becomes 
\begin{equation}
\zeta _{\lambda }^{\prime }+A^{\prime }(y_{\varepsilon }^{\ast })\zeta
_{\lambda }+\int_{0}^{1}(A^{\prime }(\nu y_{\varepsilon }^{\lambda }+(1-\nu
)y_{\varepsilon }^{\ast })-A^{\prime }(y_{\varepsilon }^{\ast }))\frac{%
y_{\varepsilon }^{\lambda }-y_{\varepsilon }^{\ast }}{\lambda }d\nu =0.
\label{g47-2-1}
\end{equation}%
We test (\ref{g47-2-1}) by $\zeta _{\lambda }(t),$ integrate with respect to 
$t,$ and get, by (\ref{g11}) and (\ref{g7}) that 
\begin{eqnarray*}
&&\frac{1}{2}\left\Vert \zeta _{\lambda }(t)\right\Vert _{H}^{2}+\alpha
_{1}\int_{0}^{t}\left\Vert \zeta _{\lambda }(\tau )\right\Vert _{V}^{2}d\tau
\leq \alpha _{2}\int_{0}^{t}\left\Vert \zeta _{\lambda }(\tau )\right\Vert
_{H}^{2}d\tau  \\
&&+\int_{0}^{1}\int_{0}^{t}\left\vert \left\langle (A^{\prime }((\nu
y_{\varepsilon }^{\lambda }+(1-\nu )y_{\varepsilon }^{\ast })(\tau
))-A^{\prime }(y_{\varepsilon }^{\ast }(\tau )))\frac{y_{\varepsilon
}^{\lambda }-y_{\varepsilon }^{\ast }}{\lambda }(\tau ),\zeta _{\lambda
}(\tau )\right\rangle _{V^{\ast },V}\right\vert d\tau d\nu  \\
&\leq &\alpha _{2}\int_{0}^{t}\left\Vert \zeta _{\lambda }(\tau )\right\Vert
_{H}^{2}d\tau +\int_{0}^{1}\int_{0}^{T}\left\Vert A^{\prime }((\nu
y_{\varepsilon }^{\lambda }+(1-\nu )y_{\varepsilon }^{\ast })(\tau
))-A^{\prime }(y_{\varepsilon }^{\ast }(\tau ))\right\Vert _{L(V,V^{\ast
})}\left\Vert \frac{y_{\varepsilon }^{\lambda }-y_{\varepsilon }^{\ast }}{%
\lambda }(\tau )\right\Vert _{V}d\tau .
\end{eqnarray*}%
By Gronwall's lemma we obtain 
\begin{eqnarray}
&&\left\Vert \zeta _{\lambda }(t)\right\Vert _{H}^{2}+\int_{0}^{t}\left\Vert
\zeta _{\lambda }(\tau )\right\Vert _{V}^{2}d\tau   \label{500} \\
&\leq &C_{T_{\varepsilon }^{\ast }}\int_{0}^{1}\left\{ \left(
\int_{0}^{T}\left\Vert A^{\prime }((\nu y_{\varepsilon }^{\lambda }+(1-\nu
)y_{\varepsilon }^{\ast })(\tau ))-A^{\prime }(y_{\varepsilon }^{\ast }(\tau
))\right\Vert _{L(V,V^{\ast })}^{2}d\tau \right) ^{1/2}\right.   \notag \\
&&\mbox{ \ \ \ \ \ \ \ }\left. \times \left( \int_{0}^{T}\left\Vert \frac{%
y_{\varepsilon }^{\lambda }-y_{\varepsilon }^{\ast }}{\lambda }(\tau
)\right\Vert _{V}^{2}d\tau \right) ^{1/2}d\tau \right\} d\nu .  \notag
\end{eqnarray}%
We recall that by (\ref{g15}), we have $\left\Vert y_{\varepsilon }^{\ast
}(t)\right\Vert _{V}\leq e^{CT}\left( \left\Vert y_{0}\right\Vert
_{V}^{2}+\rho ^{2}T\right) \leq C_{T},$ for all $t\in \lbrack 0,T],$ and 
\begin{equation*}
\left\Vert y_{\varepsilon }^{\lambda }(t)\right\Vert _{V}\leq C\left(
\left\Vert y_{0}\right\Vert _{V}^{2}+\int_{0}^{T}\left\Vert u^{\lambda
}(t)\right\Vert _{U}^{2}dt\right) e^{CT}\leq e^{CT}\left( \left\Vert
y_{0}\right\Vert _{V}^{2}+\rho ^{2}T\right) \leq C_{T}.
\end{equation*}%
Here, $C_{T}$ may change from line to line. Moreover, by (\ref{g11-2})%
\begin{equation*}
\left\Vert A^{\prime }(y_{\varepsilon }^{\ast }(\tau ))\right\Vert
_{L(V,V^{\ast })}\leq (1+C\left\Vert y_{\varepsilon }^{\ast }(\tau
)\right\Vert _{V}^{\kappa })\leq C_{T},\mbox{ and}
\end{equation*}%
\begin{equation*}
\left\Vert A^{\prime }((\nu y_{\varepsilon }^{\lambda }+(1-\nu
)y_{\varepsilon }^{\ast })(\tau ))\right\Vert _{L(V,V^{\ast })}\leq
(1+C\left\Vert (\nu y_{\varepsilon }^{\lambda }+(1-\nu )y_{\varepsilon
}^{\ast })(\tau )\right\Vert _{V}^{\kappa })\leq C_{T}.
\end{equation*}%
We also recall (\ref{g15-0}) which yields, for all $t\in \lbrack
0,T_{\varepsilon }^{\ast }],$ that%
\begin{equation}
\left\Vert (y_{\varepsilon }^{\lambda }-y_{\varepsilon }^{\ast
})(t)\right\Vert _{H}^{2}+\int_{0}^{t}\left\Vert (y_{\varepsilon }^{\lambda
}-y_{\varepsilon }^{\ast })(\tau )\right\Vert _{V}^{2}d\tau \leq C\left\Vert
u_{\varepsilon }^{\lambda }-u_{\varepsilon }^{\ast }\right\Vert
_{L^{2}(0,T;U)}^{2}\leq CT\lambda ^{2}\rho ^{2}\mbox{ \ \ \ }  \label{g47-6}
\end{equation}%
and so 
\begin{equation*}
\int_{0}^{t}\left\Vert \frac{y_{\varepsilon }^{\lambda }-y_{\varepsilon
}^{\ast }}{\lambda }(\tau )\right\Vert _{V}^{2}d\tau \leq CT\rho ^{2},%
\mbox{
for all }t\in \lbrack 0,T].
\end{equation*}%
\bigskip Therefore, 
\begin{equation}
y_{\varepsilon }^{\lambda }\rightarrow y_{\varepsilon }^{\ast }%
\mbox{
strongly in }C([0,T];H)\cap L^{2}(0,T;V),\mbox{ as }\lambda \rightarrow 0
\label{g47-3-0}
\end{equation}%
and%
\begin{equation}
\nu y_{\varepsilon }^{\lambda }+(1-\nu )y_{\varepsilon }^{\ast }\rightarrow
y_{\varepsilon }^{\ast }\mbox{ strongly in }C([0,T];H)\cap L^{2}(0,T;V),%
\mbox{ as }\lambda \rightarrow 0,  \label{g47-3}
\end{equation}%
for $\nu $ fixed, implying that 
\begin{equation}
(\nu y_{\varepsilon }^{\lambda }+(1-\nu )y_{\varepsilon }^{\ast
})(t)\rightarrow y_{\varepsilon }^{\ast }(t)\mbox{ strongly in }V%
\mbox{,
a.e. }t\in (0,T).  \label{g47-3-1}
\end{equation}%
This yields that 
\begin{equation*}
A^{\prime }((\nu y_{\varepsilon }^{\lambda }+(1-\nu )y_{\varepsilon }^{\ast
})(t))\rightarrow A^{\prime }(y_{\varepsilon }^{\ast }(t))%
\mbox{ strongly in 
}L(V,V^{\ast }),\mbox{ a.e. }t\in (0,T)
\end{equation*}%
by $(a_{4})$ and (\ref{g-000}). We denote $f_{\lambda }(t):=\left\Vert
A^{\prime }((\nu y_{\varepsilon }^{\lambda }+(1-\nu )y_{\varepsilon }^{\ast
})(t))-A^{\prime }(y_{\varepsilon }^{\ast }(t))\right\Vert _{L(V,V^{\ast
})}^{2}$ and infer that $f_{\lambda }(t)\rightarrow 0$ a.e. $t\in (0,T),$
and that $\left\vert f_{\lambda }(t)\right\vert \leq C.$ This implies, by
the Lebesgue dominated convergence theorem, that $f_{\lambda }\rightarrow 0$
in $L^{2}(0,T).$ Thus, by (\ref{500}) 
\begin{equation*}
\zeta _{\lambda }\rightarrow 0\mbox{ strongly in }C([0,T];H)\cap
L^{2}(0,T;V),\mbox{ as }\lambda \rightarrow 0.
\end{equation*}%
This proves (\ref{g47-4}). \hfill $\square $

\medskip

\noindent Now, we introduce the adjoint system%
\begin{equation}
-p_{\varepsilon }^{\prime }(t)+(A^{\prime }(y_{\varepsilon }^{\ast
}(t)))^{\ast }p_{\varepsilon }(t)=0,\mbox{ a.e.\textit{\ }}t\in
(0,T_{\varepsilon }^{\ast }),  \label{g43}
\end{equation}%
\begin{equation}
p_{\varepsilon }(T_{\varepsilon }^{\ast })=\frac{1}{\varepsilon }%
(Py_{\varepsilon }^{\ast }(T_{\varepsilon }^{\ast })-Py^{tar}).  \label{g44}
\end{equation}

\medskip

\noindent \textbf{Proposition 5.3. }\textit{Let }$Py^{tar}\in P(D_{H})$ 
\textit{and assume} (\ref{p-regular}). \textit{Then, for each }$\varepsilon
>0,$ \textit{problem} (\ref{g43})-(\ref{g44}) \textit{has a unique solution} 
\begin{equation}
p_{\varepsilon }\in C([0,T_{\varepsilon }^{\ast }];H)\cap
W^{1,2}(0,T_{\varepsilon }^{\ast };H)\cap L^{2}(0,T_{\varepsilon }^{\ast
};D_{H}).  \label{g49}
\end{equation}

\medskip

\noindent \textbf{Proof. }By Theorem 3.2 we deduce that $p_{\varepsilon
}(T_{\varepsilon }^{\ast })\in D_{H},$ since $P(y_{\varepsilon }^{\ast
}(T_{\varepsilon }^{\ast })-y^{tar})\in P(D_{H}).$ If $P\neq I,$ the second
component of $p_{\varepsilon }$ is $0\in D_{H}.$ We use the transformation $%
t\rightarrow T_{\varepsilon }^{\ast }-t$ and so (\ref{g43})-(\ref{g44})
transforms into a forward equation. The operator $(A^{\prime
}(y_{\varepsilon }^{\ast }(t)))^{\ast }$ is continuous from $V$ to $V^{\ast
} $ for all $t\geq 0$ and satisfies the properties of Lions theorem, such
that we can deduce, as in Proposition 5.1, that (\ref{g43})-(\ref{g44}) has
a unique solution 
\begin{equation}
p_{\varepsilon }\in C([0,T_{\varepsilon }^{\ast }];H)\cap
W^{1,2}(0,T_{\varepsilon }^{\ast };V^{\ast })\cap L^{2}(0,T_{\varepsilon
}^{\ast };V).  \label{g49-1}
\end{equation}%
A first estimate is obtained by testing (\ref{g43}) by $p_{\varepsilon }(t)$
and integrating over $(t,T_{\varepsilon }^{\ast }).$ Using (\ref{g11}), this
yields%
\begin{equation}
\left\Vert p_{\varepsilon }(t)\right\Vert _{H}^{2}+\int_{0}^{t}\left\Vert
p_{\varepsilon }(\tau )\right\Vert _{V}^{2}d\tau \leq C\left\Vert
p_{\varepsilon }(T_{\varepsilon }^{\ast })\right\Vert _{H}^{2},%
\mbox{ for
all }t\in \lbrack 0,T_{\varepsilon }^{\ast }].  \label{g49-2}
\end{equation}%
To prove the additional regularity we multiply (\ref{g43}) by $\Gamma _{\nu
}p_{\varepsilon }(t),$ integrate over $(t,T_{\varepsilon }^{\ast })$ and use
(\ref{p-regular}). We have%
\begin{eqnarray*}
&&\frac{1}{2}\left\langle p_{\varepsilon }(t),\Gamma _{\nu }p_{\varepsilon
}(t)\right\rangle _{V^{\ast },V}+\gamma _{1}\int_{0}^{t}\left\Vert \Gamma
_{\nu }p_{\varepsilon }(\tau )\right\Vert _{H}^{2}d\tau \\
&\leq &\frac{1}{2}\left\langle p_{\varepsilon }(T_{\varepsilon }^{\ast
}),\Gamma _{\nu }p_{\varepsilon }(T_{\varepsilon }^{\ast })\right\rangle
_{V^{\ast },V}+\gamma _{2}\int_{0}^{t}\left\Vert p_{\varepsilon }(\tau
)\right\Vert _{V}^{2}(1+\gamma _{3}\left\Vert y_{\varepsilon }^{\ast }(\tau
)\right\Vert _{V}^{l})d\tau .
\end{eqnarray*}%
Taking into account (\ref{g49-2}) and (\ref{g15}), that is $\left\Vert
y_{\varepsilon }^{\ast }(\tau )\right\Vert _{V}\leq C_{T_{\varepsilon
}^{\ast }}<C_{T^{\ast }+1},$ we obtain 
\begin{eqnarray}
&&\left\Vert \Gamma _{\nu }^{1/2}p_{\varepsilon }(t)\right\Vert
_{H}^{2}+\int_{0}^{t}\left\Vert \Gamma _{\nu }p_{\varepsilon }(\tau
)\right\Vert _{H}^{2}d\tau  \label{g49-3} \\
&\leq &C\left( \left\Vert p_{\varepsilon }(T_{\varepsilon }^{\ast
})\right\Vert _{H}\left\Vert \Gamma _{H}p_{\varepsilon }(T_{\varepsilon
}^{\ast })\right\Vert _{H}+1\right) \leq C\left( \left\Vert p_{\varepsilon
}(T_{\varepsilon }^{\ast })\right\Vert _{H}\left\Vert p_{\varepsilon
}(T_{\varepsilon }^{\ast })\right\Vert _{D_{H}}+1\right) ,  \notag
\end{eqnarray}%
for all $t\in \lbrack 0,T_{\varepsilon }^{\ast }],$ since $\left\Vert \Gamma
_{H}p_{\varepsilon }(\tau )\right\Vert _{H}=\left\Vert p_{\varepsilon }(\tau
)\right\Vert _{D_{H}}.$

\noindent Here, we used the relation $\left\Vert \Gamma _{\nu
}p_{\varepsilon }(T_{\varepsilon }^{\ast })\right\Vert _{H}\leq \left\Vert
\Gamma _{H}p_{\varepsilon }(T_{\varepsilon }^{\ast })\right\Vert _{H}$ for $%
p_{\varepsilon }(T_{\varepsilon }^{\ast })\in D_{H}.$ Now, we can pass to
the limit as $\nu \rightarrow 0$ and obtain%
\begin{equation*}
\Gamma _{\nu }p_{\varepsilon }\rightarrow \Gamma _{H}p\mbox{ weakly in }%
L^{2}(0,T_{\varepsilon }^{\ast };H),\mbox{ weak-star in }L^{\infty
}(0,T_{\varepsilon }^{\ast };V),
\end{equation*}%
and so by (\ref{g49-3}) we get%
\begin{equation*}
\left\Vert p_{\varepsilon }(t)\right\Vert _{V}^{2}+\int_{0}^{t}\left\Vert
\Gamma _{H}p_{\varepsilon }(\tau )\right\Vert _{H}^{2}d\tau \leq C\left(
\left\Vert p_{\varepsilon }(T_{\varepsilon }^{\ast })\right\Vert
_{H}\left\Vert p_{\varepsilon }(T_{\varepsilon }^{\ast })\right\Vert
_{D_{H}}+1\right) ,\mbox{ for all }t\in \lbrack 0,T_{\varepsilon }^{\ast }].
\end{equation*}%
Thus, $p_{\varepsilon }\in L^{2}(0,T_{\varepsilon }^{\ast };D_{H})\cap
L^{\infty }(0,T_{\varepsilon }^{\ast };V).$ For a.a. $t\in (0,T_{\varepsilon
}^{\ast })$ we still have 
\begin{equation*}
\left\Vert (A_{H}^{\prime }(y_{\varepsilon }^{\ast }(t)))^{\ast
}p_{\varepsilon }(t)\right\Vert _{H}\leq \left\Vert p_{\varepsilon
}(t)\right\Vert _{D_{H}}(1+C\left\Vert y_{\varepsilon }^{\ast
}(t)\right\Vert _{V}^{\kappa }).
\end{equation*}%
By (\ref{g43}) it follows that $p_{\varepsilon }^{\prime }\in
L^{2}(0,T_{\varepsilon }^{\ast };H)$ and so (\ref{g49}) is proved. \hfill $%
\square $

\subsection{Approximating optimality conditions}

Let us introduce the sets 
\begin{equation*}
K=\{w\in U;\mbox{ }\left\Vert w\right\Vert _{U}\leq \rho \},\mbox{ }\mathcal{%
K}_{T}=\{z\in L^{2}(0,T;U);\mbox{ }z(t)\in K\mbox{ a.e. }t\in (0,T)\},
\end{equation*}%
and denote the normal cone to $K$ at $w$ by%
\begin{equation}
N_{K}(w)=\{\zeta \in U^{\ast };\mbox{ }\left\langle \zeta ,w-\overline{w}%
\right\rangle _{U^{\ast },U}\geq 0,\mbox{ for all }\overline{w}\in K\},%
\mbox{
}  \label{con}
\end{equation}%
and the normal cone to $\mathcal{K}_{T}$ at $\omega $ by%
\begin{equation}
N_{\mathcal{K}_{T}}(\omega )=\left\{ \chi \in L^{2}(0,T;U^{\ast });\mbox{ }%
\int_{0}^{T}\left\langle \chi (t),(\omega -\overline{\omega }%
)(t)\right\rangle _{U^{\ast },U}dt\geq 0,\mbox{ for all }\overline{\omega }%
\in \mathcal{K}_{T}\right\} \mbox{.}  \label{con-cal}
\end{equation}%
We recall (see e.g., \cite{VB-optim}) that $\chi \in N_{\mathcal{K}%
_{T}}(\omega )$ iff $\chi (t)\in N_{K}(\omega (t))$ a.e. $t\in (0,T).$

\noindent We denote by $F:U\rightarrow U^{\ast }$ the duality mapping of $U$
(see (\ref{A1}) in the Appendix) and recall that $h_{\varepsilon }^{\ast
}(t) $ was defined in (\ref{g22-4}), $h_{\varepsilon }^{\ast
}(t)=\int_{0}^{t}P(u_{\varepsilon }^{\ast }-u^{\ast })(\tau )d\tau .$

\medskip

\noindent \textbf{Proposition 5.4. }\textit{Assume }%
\begin{equation}
y_{0}\in V,\mbox{ }Py^{tar}\in P(D_{H}),\mbox{ }Py_{0}\neq Py^{tar}.
\label{g40}
\end{equation}%
\textit{Let }$(T_{\varepsilon }^{\ast },u_{\varepsilon }^{\ast })$\textit{\
be an optimal control in }$(P_{\varepsilon })$ \textit{with the optimal state%
} $y_{\varepsilon }^{\ast }$\textit{. Then, }%
\begin{equation}
Pu_{\varepsilon }^{\ast }(t)=-(\varepsilon F+N_{K})^{-1}\left( B^{\ast
}p_{\varepsilon }(t)+\int_{t}^{T_{\varepsilon }^{\ast }}F(h_{\varepsilon
}^{\ast }(\tau ))d\tau \right) ,\mbox{ \textit{for all} }t\in \lbrack
0,T_{\varepsilon }^{\ast }],  \label{g41}
\end{equation}%
\textit{and } 
\begin{eqnarray}
&&\rho \left\Vert B^{\ast }p_{\varepsilon }(t)+\int_{t}^{T_{\varepsilon
}^{\ast }}F(h_{\varepsilon }^{\ast }(\tau ))d\tau +\varepsilon
F(Pu_{\varepsilon }^{\ast }(t))\right\Vert _{U^{\ast }}+\left(
A_{H}y_{\varepsilon }^{\ast }(t),p_{\varepsilon }(t)\right) _{H}+  \notag \\
&&+\int_{t}^{T_{\varepsilon }^{\ast }}\left( Pu_{\varepsilon }^{\ast }(\tau
),F(h_{\varepsilon }^{\ast }(\tau ))\right) _{U,U^{\ast }}d\tau +\frac{%
\varepsilon }{2}\left\Vert Pu_{\varepsilon }^{\ast }(t)\right\Vert _{U}^{2} 
\notag \\
&=&1+\frac{1}{2}\left\Vert h_{\varepsilon }^{\ast }(T_{\varepsilon }^{\ast
})\right\Vert _{U}^{2},\mbox{ }t\in \lbrack 0,T_{\varepsilon }^{\ast }],
\label{g42}
\end{eqnarray}%
\textit{where }$p_{\varepsilon }$\ \textit{is the solution to the adjoint
equation }(\ref{g43})-(\ref{g44}). \textit{Moreover, }$t\rightarrow
u_{\varepsilon }^{\ast }(t)$\textit{\ turns out to be continuous on} $%
[0,T_{\varepsilon }^{\ast }].$

\medskip

\noindent \textbf{Proof. }Let $(T_{\varepsilon }^{\ast },u_{\varepsilon
}^{\ast })$ be an optimal controller in $(P_{\varepsilon }).$ We shall
compute separate variations with respect to $T_{\varepsilon }^{\ast }$ and $%
u_{\varepsilon }^{\ast }$. By the condition of optimality for $%
u_{\varepsilon }^{\ast }$ we have 
\begin{equation*}
J_{\varepsilon }(T_{\varepsilon }^{\ast },u_{\varepsilon }^{\ast })\leq
J_{\varepsilon }(T,u),\mbox{ for all }u(t)\in K,\mbox{ }\left\Vert
u(t)\right\Vert _{U}\leq \rho ,\mbox{ a.e. }t\geq 0.
\end{equation*}%
In particular, replacing $u$ by $u_{\varepsilon }^{\lambda }=Pu_{\varepsilon
}^{\ast }+\lambda v$ with $v=P(\overline{u}-u_{\varepsilon }^{\ast }),$ $%
\overline{u}(t)\in K,$ and performing some calculations, recalling (\ref%
{A4-1-0}) we get%
\begin{eqnarray}
&&\left( \frac{1}{\varepsilon }(Py_{\varepsilon }^{\ast }(T_{\varepsilon
}^{\ast })-Py^{tar}),PY(T_{\varepsilon }^{\ast })\right) _{H}+\varepsilon
\int_{0}^{T_{\varepsilon }^{\ast }}\left\langle F(Pu_{\varepsilon }^{\ast
}(t)),v(t)\right\rangle _{U^{\ast },U}dt  \label{g50} \\
&&+\int_{0}^{T_{\varepsilon }^{\ast }}\left\langle F(h_{\varepsilon }^{\ast
}(t)),\int_{0}^{t}v(s)ds\right\rangle _{U^{\ast },U}dt\geq 0.  \notag
\end{eqnarray}%
Observing that%
\begin{eqnarray}
&&\int_{0}^{T_{\varepsilon }^{\ast }}\left\langle F(h_{\varepsilon }^{\ast
}(t)),\int_{0}^{t}v(s)ds\right\rangle _{U^{\ast
},U}dt=\int_{0}^{T_{\varepsilon }^{\ast }}\int_{0}^{t}\left\langle
F(h_{\varepsilon }^{\ast }(t)),v(s)\right\rangle _{U^{\ast },U}dsdt
\label{g51} \\
&=&\int_{0}^{T_{\varepsilon }^{\ast }}\int_{s}^{T_{\varepsilon }^{\ast
}}\left\langle F(h_{\varepsilon }^{\ast }(t)),v(s)\right\rangle _{U^{\ast
},U}dtds=\int_{0}^{T_{\varepsilon }^{\ast }}\left\langle
\int_{s}^{T_{\varepsilon }^{\ast }}F(h_{\varepsilon }^{\ast
}(t))dt,v(s)\right\rangle _{U^{\ast },U}ds  \notag
\end{eqnarray}%
we obtain 
\begin{equation}
\left( \frac{1}{\varepsilon }(Py_{\varepsilon }^{\ast }(T_{\varepsilon
}^{\ast })-Py^{tar}),Y(T_{\varepsilon }^{\ast })\right)
_{H}+\int_{0}^{T_{\varepsilon }^{\ast }}\left\langle \varepsilon
F(Pu_{\varepsilon }^{\ast }(t))+\int_{t}^{T_{\varepsilon }^{\ast
}}F(h_{\varepsilon }^{\ast }(\tau ))d\tau ,v(t)\right\rangle _{U^{\ast
},U}dt\geq 0.  \label{g52}
\end{equation}%
Here we used that $P^{2}=P$ and the fact that $(Pw,P\overline{w})_{H}$ is
the same with $(Pw,\overline{w})$ when $Pw=(w_{1},0).$ We test (\ref{g47})
by $p_{\varepsilon }(t)$ and integrate over $(0,T_{\varepsilon }^{\ast }).$
By a straightforward calculation we obtain 
\begin{eqnarray*}
&&\int_{\Omega }\left( p_{\varepsilon }(T_{\varepsilon }^{\ast
})Y(T_{\varepsilon }^{\ast })-p_{\varepsilon }(0)Y(0)\right)
dx+\int_{0}^{T_{\varepsilon }^{\ast }}\left\langle (-p_{\varepsilon
}^{\prime }+(A^{\prime }(y_{\varepsilon }^{\ast }))^{\ast }p_{\varepsilon
})(t),Y(t)\right\rangle _{V^{\ast },V}dt\mbox{ \ \ } \\
&=&\int_{0}^{T_{\varepsilon }^{\ast }}(Bv(t),p_{\varepsilon }(t))_{H}dt.
\end{eqnarray*}%
Using again the adjoint system, this equation reduces to 
\begin{equation}
\left( \frac{1}{\varepsilon }P(y_{\varepsilon }^{\ast }(T_{\varepsilon
}^{\ast })-Py^{tar}),Y(T_{\varepsilon }^{\ast })\right)
_{H}=\int_{0}^{T_{\varepsilon }^{\ast }}\left\langle B^{\ast }p_{\varepsilon
}(t),v(t)\right\rangle _{U^{\ast },U}dt.  \label{g53}
\end{equation}%
We recall that $v=P(\overline{u}-u_{\varepsilon }^{\ast }).$ Replacing the
left-hand side of (\ref{g53}) into (\ref{g52}) we deduce that%
\begin{equation}
\int_{0}^{T_{\varepsilon }^{\ast }}\left( B^{\ast }p_{\varepsilon
}(t),v(t)\right) _{U^{\ast },U}dt+\int_{0}^{T_{\varepsilon }^{\ast }}\left(
\varepsilon F(Pu_{\varepsilon }^{\ast }(t))+\int_{t}^{T_{\varepsilon }^{\ast
}}F(h_{\varepsilon }^{\ast }(\tau ))d\tau ,v(t)\right) _{U^{\ast },U}dt\geq
0,  \label{g54}
\end{equation}%
for all $\overline{u}(t)\in K,$ that is $\overline{u}(t)\in U,$ $\left\Vert 
\overline{u}(t)\right\Vert _{U}\leq \rho $ a.e. $t\geq 0.$ This yields%
\begin{equation*}
\int_{0}^{T_{\varepsilon }^{\ast }}\left\langle -B^{\ast }p_{\varepsilon
}(t)-\varepsilon F(u_{\varepsilon }^{\ast }(t))-\int_{t}^{T_{\varepsilon
}^{\ast }}F(h_{\varepsilon }^{\ast }(\tau ))d\tau ,Pu_{\varepsilon }^{\ast
}(t)-P\overline{u}(t)\right\rangle _{U^{\ast },U}dt\geq 0,\mbox{ }
\end{equation*}%
for all $\overline{u}(t)\in K,$ a.e. $t\in (0,T_{\varepsilon }^{\ast }),$
and implies, by (\ref{con}), that%
\begin{equation}
z_{\varepsilon }(t):=-B^{\ast }p_{\varepsilon }(t)-\varepsilon
F(Pu_{\varepsilon }^{\ast }(t))-\int_{t}^{T_{\varepsilon }^{\ast
}}F(h_{\varepsilon }^{\ast }(\tau ))d\tau \in N_{K}(Pu_{\varepsilon }^{\ast
}(t)),\mbox{ a.e. }t\in (0,T_{\varepsilon }^{\ast }),  \label{g55}
\end{equation}%
or, equivalently, 
\begin{equation}
Pu_{\varepsilon }^{\ast }(t)=(\varepsilon F+N_{K})^{-1}\left( -B^{\ast
}p_{\varepsilon }(t)-\int_{t}^{T_{\varepsilon }^{\ast }}F(h_{\varepsilon
}^{\ast }(\tau ))d\tau \right) ,\mbox{ a.e. }t\in (0,T_{\varepsilon }^{\ast
}).  \label{g56}
\end{equation}%
Moreover, relation (\ref{g56}) implies that $t\rightarrow u_{\varepsilon
}^{\ast }(t)$ is continuous, because $(\varepsilon F+N_{K})^{-1}$ is
single-valued and Lipschitz continuous, the integral is continuous and $%
p_{\varepsilon }$ belongs to $C([0,T_{\varepsilon }^{\ast }];H),$ so that (%
\ref{g56}) is true for all $t\in \lbrack 0,T_{\varepsilon }^{\ast }].$ We
also note that 
\begin{equation*}
F(h_{\varepsilon }^{\ast }(\tau ))=F\left( \int_{0}^{\tau }P(u_{\varepsilon
}^{\ast }-u^{\ast })(s)ds\right) =PF\left( \int_{0}^{\tau }P(u_{\varepsilon
}^{\ast }-u^{\ast })(s)ds\right) .
\end{equation*}

We recall by (\ref{g31-1}) that $A_{H}y_{\varepsilon }^{\ast }(t)\in H,$ $%
(y_{\varepsilon }^{\ast })^{\prime }(t)\in H$ a.e. $t\in (0,T_{\varepsilon
}^{\ast }).$ Also, $y_{\varepsilon }^{\ast }\in C_{w}([0,T_{\varepsilon
}^{\ast }];V)$ and $Ay\in C_{w}([0,T_{\varepsilon }^{\ast }];V^{\ast })$
(see the observation after Theorem 3.2). By the state equation we have $%
(y_{\varepsilon }^{\ast })^{\prime }=-Ay_{\varepsilon }^{\ast
}+Bu_{\varepsilon }^{\ast }\in C_{w}([0,T_{\varepsilon }^{\ast }];V^{\ast
}). $ Indeed, $y_{\varepsilon }^{\ast }(T_{\varepsilon }^{\ast })\in V,$ so $%
Ay_{\varepsilon }^{\ast }(T_{\varepsilon }^{\ast })\in V^{\ast }$ and $%
Bu_{\varepsilon }^{\ast }(T_{\varepsilon }^{\ast })\in H.$ Thus, $%
(y_{\varepsilon }^{\ast })^{\prime }(T_{\varepsilon }^{\ast })\in V^{\ast },$
for all $\varepsilon >0.$ Recalling (\ref{g44}) and that $Py_{\varepsilon
}^{\ast }(T_{\varepsilon }^{\ast })-Py^{tar}\in P(D_{H}),$ we deduce that $%
p_{\varepsilon }(T_{\varepsilon }^{\ast })\in D_{H},$ for all $\varepsilon
>0.$

Next, we keep $u_{\varepsilon }^{\ast }$ fixed and give variations to $%
T_{\varepsilon }^{\ast }$. Since $T_{\varepsilon }^{\ast }$ realizes the
minimum in $(P_{\varepsilon })$ we can write 
\begin{equation*}
J_{\varepsilon }(T_{\varepsilon }^{\ast },u_{\varepsilon }^{\ast })\leq
J_{\varepsilon }(T_{\varepsilon }^{\ast }+\lambda ,u_{\varepsilon }^{\ast }),%
\mbox{ }\lambda >0,
\end{equation*}%
that is, 
\begin{eqnarray*}
&&J_{\varepsilon }(T_{\varepsilon }^{\ast },u_{\varepsilon }^{\ast
})=T_{\varepsilon }^{\ast }+\frac{1}{2\varepsilon }\left\Vert
P(y_{\varepsilon }^{T_{\varepsilon }^{\ast },u_{\varepsilon }^{\ast
}}(T_{\varepsilon }^{\ast })-Py^{tar})\right\Vert _{H}^{2} \\
&&+\frac{\varepsilon }{2}\int_{0}^{T_{\varepsilon }^{\ast }}\left\Vert
Pu_{\varepsilon }^{\ast }(t)\right\Vert _{U}^{2}dt+\frac{1}{2}%
\int_{0}^{T_{\varepsilon }^{\ast }}\left\Vert \int_{0}^{t}P(u_{\varepsilon
}^{\ast }(\tau )-u^{\ast }(\tau ))d\tau \right\Vert _{U}^{2}dt \\
&\leq &J_{\varepsilon }(T_{\varepsilon }^{\ast }+\lambda ,u_{\varepsilon
}^{\ast })=T_{\varepsilon }^{\ast }+\lambda +\frac{1}{2\varepsilon }%
\left\Vert P(y_{\varepsilon }^{T_{\varepsilon }^{\ast }+\lambda
,u_{\varepsilon }^{\ast }}(T_{\varepsilon }^{\ast }+\lambda
)-Py^{tar})\right\Vert _{H}^{2}+\frac{\varepsilon }{2}\int_{0}^{T_{%
\varepsilon }^{\ast }+\lambda }\left\Vert Pu_{\varepsilon }^{\ast
}(t)\right\Vert _{U}^{2}dt \\
&&+\frac{1}{2}\int_{0}^{T_{\varepsilon }^{\ast }+\lambda }\left\Vert
\int_{0}^{t}P(u_{\varepsilon }^{\ast }(\tau )-u^{\ast }(\tau ))d\tau
\right\Vert _{U}^{2}dt.
\end{eqnarray*}%
In these calculations we took into account that $u_{\varepsilon }^{\ast }$
and the solution to the approximating state are continuous with respect to $%
t\in \lbrack 0,\infty ).$ Then, the solution $y_{\varepsilon
}^{T_{\varepsilon }^{\ast }+\lambda ,u_{\varepsilon }^{\ast }}(t)$
calculated for $t\in (0,T_{\varepsilon }^{\ast }+\lambda )$ and $%
u_{\varepsilon }^{\ast },$ restricted to $(0,T_{\varepsilon }^{\ast })$
coincides with $y_{\varepsilon }^{T_{\varepsilon }^{\ast },u_{\varepsilon
}^{\ast }}(t)$ the solution calculated on $(0,T_{\varepsilon }^{\ast })$,
which was denoted by $y_{\varepsilon }^{\ast }(t).$ Performing some
calculations we get 
\begin{equation*}
1+\frac{1}{\varepsilon }\left\langle (y_{\varepsilon }^{\ast })^{\prime
}(T_{\varepsilon }^{\ast }),Py_{\varepsilon }^{\ast }(T_{\varepsilon }^{\ast
})-Py^{tar}\right\rangle _{V^{\ast },V}+\frac{1}{2}\left\Vert h_{\varepsilon
}^{\ast }(T_{\varepsilon }^{\ast })\right\Vert _{U}^{2}+\frac{\varepsilon }{2%
}\left\Vert Pu_{\varepsilon }^{\ast }(T_{\varepsilon }^{\ast })\right\Vert
_{U}^{2}\geq 0.
\end{equation*}%
Doing the same for $T_{\varepsilon }^{\ast }-\lambda $ and observing that
the solution $y_{\varepsilon }^{T_{\varepsilon }^{\ast },u_{\varepsilon
}^{\ast }}(t),$ calculated for $t\in (0,T_{\varepsilon }^{\ast })$ and $%
u_{\varepsilon }^{\ast },$ restricted to $(0,T_{\varepsilon }^{\ast
}-\lambda )$ is in fact $y_{\varepsilon }^{T_{\varepsilon }^{\ast }-\lambda
,u_{\varepsilon }^{\ast }}(t)$ the solution calculated on $(0,T_{\varepsilon
}^{\ast }-\lambda ),$ we get the reverse inequality. Finally, we obtain%
\begin{equation}
1+\frac{1}{\varepsilon }\left\langle (y_{\varepsilon }^{\ast })^{\prime
}(T_{\varepsilon }^{\ast }),Py_{\varepsilon }^{\ast }(T_{\varepsilon }^{\ast
})-Py^{tar}\right\rangle _{V^{\ast },V}+\frac{1}{2}\left\Vert h_{\varepsilon
}^{\ast }(T_{\varepsilon }^{\ast })\right\Vert _{U}^{2}+\frac{\varepsilon }{2%
}\left\Vert Pu_{\varepsilon }^{\ast }(T_{\varepsilon }^{\ast })\right\Vert
_{U}^{2}=0.  \label{g57}
\end{equation}%
Then, using the state system (\ref{g4}) for $y_{\varepsilon }^{\ast }$ and
the final conditions of the adjoint system, we can express the term 
\begin{equation*}
\frac{1}{\varepsilon }\left\langle (y_{\varepsilon }^{\ast })^{\prime
}(T_{\varepsilon }^{\ast }),Py_{\varepsilon }^{\ast }(T_{\varepsilon }^{\ast
})-Py^{tar}\right\rangle _{V^{\ast },V}=\left\langle Bu_{\varepsilon }^{\ast
}(T_{\varepsilon }^{\ast })-Ay_{\varepsilon }^{\ast }(T_{\varepsilon }^{\ast
}),p_{\varepsilon }(T_{\varepsilon }^{\ast })\right\rangle _{V^{\ast },V}.
\end{equation*}%
Plugging this in (\ref{g57}), we obtain 
\begin{equation}
1+\left\langle u_{\varepsilon }^{\ast }(T_{\varepsilon }^{\ast }),B^{\ast
}p_{\varepsilon }(T_{\varepsilon }^{\ast })\right\rangle _{U,U^{\ast
}}-\left\langle AX_{\varepsilon }^{\ast }(T_{\varepsilon }^{\ast
}),p_{\varepsilon }(T_{\varepsilon }^{\ast })\right\rangle _{V^{\ast },V}+%
\frac{1}{2}\left\Vert h_{\varepsilon }^{\ast }(T_{\varepsilon }^{\ast
})\right\Vert _{U}^{2}+\frac{\varepsilon }{2}\left\Vert Pu_{\varepsilon
}^{\ast }(T_{\varepsilon }^{\ast })\right\Vert _{U}^{2}=0.  \label{g58}
\end{equation}%
We replace $Pu_{\varepsilon }^{\ast }(T_{\varepsilon }^{\ast })$ from (\ref%
{g56}), 
\begin{equation*}
Pu_{\varepsilon }^{\ast }(T_{\varepsilon }^{\ast })=(\varepsilon
F+N_{K})^{-1}(-B^{\ast }p_{\varepsilon }(T_{\varepsilon }^{\ast })),
\end{equation*}%
which can be still written 
\begin{equation}
\varepsilon F(Pu_{\varepsilon }^{\ast }(T_{\varepsilon }^{\ast
}))+z_{\varepsilon }^{\ast }(T_{\varepsilon }^{\ast })=-B^{\ast
}p_{\varepsilon }(T_{\varepsilon }^{\ast }),\mbox{ where }z_{\varepsilon
}^{\ast }(T_{\varepsilon }^{\ast })\in N_{K}(Pu_{\varepsilon }^{\ast
}(T_{\varepsilon }^{\ast })).  \label{g58-0}
\end{equation}%
By using this and (\ref{A4-0}) we obtain for the second term in (\ref{g58})%
\begin{eqnarray*}
\left\langle u_{\varepsilon }^{\ast }(T_{\varepsilon }^{\ast }),B^{\ast
}p_{\varepsilon }(T_{\varepsilon }^{\ast })\right\rangle _{U,U^{\ast }}
&=&-\left\langle u_{\varepsilon }^{\ast }(T_{\varepsilon }^{\ast
}),\varepsilon F(Pu_{\varepsilon }^{\ast }(T_{\varepsilon }^{\ast
}))+z_{\varepsilon }^{\ast }(T_{\varepsilon }^{\ast })\right\rangle
_{U,U^{\ast }} \\
&=&-\varepsilon \left\Vert Pu_{\varepsilon }^{\ast }(T_{\varepsilon }^{\ast
})\right\Vert _{U}^{2}-\rho \left\Vert z_{\varepsilon }^{\ast
}(T_{\varepsilon }^{\ast })\right\Vert _{U^{\ast }} \\
&=&-\varepsilon \left\Vert Pu_{\varepsilon }^{\ast }(T_{\varepsilon }^{\ast
})\right\Vert _{U}^{2}-\rho \left\Vert B^{\ast }p_{\varepsilon
}(T_{\varepsilon }^{\ast })+\varepsilon F(Pu_{\varepsilon }^{\ast
}(T_{\varepsilon }^{\ast }))\right\Vert _{U^{\ast }}.
\end{eqnarray*}%
Therefore, (\ref{g58}) becomes 
\begin{eqnarray*}
&&1-\rho \left\Vert B^{\ast }p_{\varepsilon }(T_{\varepsilon }^{\ast
})+\varepsilon F(Pu_{\varepsilon }^{\ast }(T_{\varepsilon }^{\ast
}))\right\Vert _{U^{\ast }}-\varepsilon \left\Vert Pu_{\varepsilon }^{\ast
}(T_{\varepsilon }^{\ast })\right\Vert _{U}^{2}-\left\langle Ay_{\varepsilon
}^{\ast }(T_{\varepsilon }^{\ast }),p_{\varepsilon }(T_{\varepsilon }^{\ast
})\right\rangle _{V^{\ast },V} \\
&&+\frac{1}{2}\left\Vert h_{\varepsilon }^{\ast }(T_{\varepsilon }^{\ast
})\right\Vert _{U}^{2}+\frac{\varepsilon }{2}\left\Vert Pu_{\varepsilon
}^{\ast }(T_{\varepsilon }^{\ast })\right\Vert _{U}^{2}=0,
\end{eqnarray*}%
which finally can be written 
\begin{equation}
\rho \left\Vert B^{\ast }p_{\varepsilon }(T_{\varepsilon }^{\ast
})+\varepsilon F(Pu_{\varepsilon }^{\ast }(T_{\varepsilon }^{\ast
}))\right\Vert _{U^{\ast }}+\left\langle Ay_{\varepsilon }^{\ast
}(T_{\varepsilon }^{\ast }),p_{\varepsilon }(T_{\varepsilon }^{\ast
})\right\rangle _{V^{\ast },V}+\frac{\varepsilon }{2}\left\Vert
Pu_{\varepsilon }^{\ast }(T_{\varepsilon }^{\ast })\right\Vert _{U}^{2}=1+%
\frac{1}{2}\left\Vert h_{\varepsilon }^{\ast }(T_{\varepsilon }^{\ast
})\right\Vert _{U}^{2}.  \label{g59}
\end{equation}%
The next calculation can be performed due to the supplementary regularity of 
$p_{\varepsilon },$ that is $p_{\varepsilon }^{\prime }\in
L^{2}(0,T_{\varepsilon }^{\ast };H),$ given by (\ref{g49}). We multiply
scalarly the state equation by $p_{\varepsilon }^{\prime }(t)$, add with the
adjoint equation multiplied by $(y_{\varepsilon }^{\ast })^{\prime }(t)$,
getting%
\begin{equation*}
(A_{H}y_{\varepsilon }^{\ast }(t),p_{\varepsilon }^{\prime
}(t))_{H}+((A_{H}^{\prime }(y_{\varepsilon }^{\ast }(t)))^{\ast
}p_{\varepsilon }(t),(y_{\varepsilon }^{\ast })^{\prime
}(t))_{H}=(Bu_{\varepsilon }^{\ast }(t),p_{\varepsilon }^{\prime }(t))_{H},
\end{equation*}%
a.e. $t\in (0,T_{\varepsilon }^{\ast }),$ that reduces to%
\begin{equation*}
\left( A_{H}y_{\varepsilon }^{\ast }(t),p_{\varepsilon }(t)\right)
_{H}^{\prime }=\left\langle u_{\varepsilon }^{\ast }(t),B^{\ast
}p_{\varepsilon }^{\prime }(t)\right\rangle _{U,U^{\ast }},\mbox{ a.e. }t\in
(0,T_{\varepsilon }^{\ast }).
\end{equation*}%
We integrate on $(t,T_{\varepsilon }^{\ast })$ and obtain 
\begin{eqnarray}
&&(A_{H}y_{\varepsilon }^{\ast }(T_{\varepsilon }^{\ast }),p_{\varepsilon
}(T_{\varepsilon }^{\ast }))_{H}-(A_{H}y_{\varepsilon }^{\ast
}(t),p_{\varepsilon }(t))_{H}=\int_{t}^{T_{\varepsilon }^{\ast
}}\left\langle u_{\varepsilon }^{\ast }(s),B^{\ast }p_{\varepsilon }^{\prime
}(s)\right\rangle _{U,U^{\ast }}ds  \label{g60} \\
&=&\int_{t}^{T_{\varepsilon }^{\ast }}\left\langle Pu_{\varepsilon }^{\ast
}(s),\left( B^{\ast }p_{\varepsilon }(s)+\int_{s}^{T_{\varepsilon }^{\ast
}}F(h_{\varepsilon }^{\ast }(\tau ))d\tau )\right) ^{\prime }\right\rangle
_{U,U^{\ast }}ds  \notag \\
&&-\int_{t}^{T_{\varepsilon }^{\ast }}\left\langle Pu_{\varepsilon }^{\ast
}(s),\left( \int_{s}^{T_{\varepsilon }^{\ast }}F(h_{\varepsilon }^{\ast
}(\tau ))d\tau )\right) ^{\prime }\right\rangle _{U,U^{\ast }}ds,  \notag
\end{eqnarray}%
since $\left\langle u_{\varepsilon }^{\ast }(s),Pk(s)\right\rangle
_{U,U^{\ast }}=\left\langle Pu_{\varepsilon }^{\ast }(s),Pk(s)\right\rangle
_{U,U^{\ast }}$where $k(s)$ is either $p_{\varepsilon }(s)$ or $%
\int_{s}^{T_{\varepsilon }^{\ast }}F(h_{\varepsilon }^{\ast }(\tau ))d\tau $%
, both containing $P$ in their expressions. Denoting%
\begin{equation*}
\zeta _{\varepsilon }(t)=-B^{\ast }p_{\varepsilon
}(t)-\int_{t}^{T_{\varepsilon }^{\ast }}F(h_{\varepsilon }^{\ast }(\tau
))d\tau 
\end{equation*}%
we see by (\ref{g56}) that $\zeta _{\varepsilon }(t)=(\varepsilon
F+N_{K})(Pu_{\varepsilon }^{\ast }(t)).$ By (\ref{A4-3}), 
\begin{eqnarray}
j_{\varepsilon }^{\ast }(\zeta _{\varepsilon }(t)) &=&\frac{\varepsilon }{2}%
\left\Vert Pu_{\varepsilon }^{\ast }(t)\right\Vert _{U}^{2}+\rho \left\Vert
z_{\varepsilon }(t)\right\Vert _{U^{\ast }}  \label{g60-1} \\
&=&\frac{\varepsilon }{2}\left\Vert Pu_{\varepsilon }^{\ast }(t)\right\Vert
_{U}^{2}+\rho \left\Vert B^{\ast }p_{\varepsilon
}(t)+\int_{t}^{T_{\varepsilon }^{\ast }}F(h_{\varepsilon }^{\ast }(\tau
))d\tau +\varepsilon F(Pu_{\varepsilon }^{\ast }(t))\right\Vert _{U^{\ast }},%
\mbox{ }t\in \lbrack 0,T_{\varepsilon }^{\ast }].  \notag
\end{eqnarray}%
Thus, we can express (\ref{g56}) as 
\begin{equation}
Pu_{\varepsilon }^{\ast }(t)=(\varepsilon F+N_{K})^{-1}(\zeta _{\varepsilon
})=(\partial j_{\varepsilon })^{-1}(\zeta _{\varepsilon }(t))=\partial
j_{\varepsilon }^{\ast }(\zeta _{\varepsilon }(t)),  \label{g61}
\end{equation}%
(see (\ref{A4-1})-(\ref{A4-3})). Then, the integrand of the first term on
the right-hand side in (\ref{g60}) becomes 
\begin{eqnarray*}
&&\left\langle Pu_{\varepsilon }^{\ast }(t),\left( B^{\ast }p_{\varepsilon
}(t)+\int_{t}^{T_{\varepsilon }^{\ast }}F(h_{\varepsilon }^{\ast }(\tau
))d\tau \right) ^{\prime }\right\rangle _{U,U^{\ast }} \\
&=&-\left\langle Pu_{\varepsilon }^{\ast }(t),\left( -B^{\ast
}p_{\varepsilon }(t)-\int_{t}^{T_{\varepsilon }^{\ast }}F(h_{\varepsilon
}^{\ast }(\tau ))d\tau \right) ^{\prime }\right\rangle _{U,U^{\ast }} \\
&=&-\left( \partial j_{\varepsilon }^{\ast }(\zeta _{\varepsilon }(t)),\zeta
_{\varepsilon }^{\prime }(t)\right) _{U^{\ast }}=-\frac{dj_{\varepsilon
}^{\ast }}{dt}\left( \zeta _{\varepsilon }(t)\right) .
\end{eqnarray*}%
Plugging this in (\ref{g60}) we get%
\begin{eqnarray}
&&\left( A_{H}y_{\varepsilon }^{\ast }(T_{\varepsilon }^{\ast
}),p_{\varepsilon }(T_{\varepsilon }^{\ast })\right) _{H}+j_{\varepsilon
}^{\ast }\left( \zeta _{\varepsilon }(T_{\varepsilon }^{\ast })\right) 
\label{g62} \\
&=&\left( A_{H}y_{\varepsilon }^{\ast }(t),p_{\varepsilon }(t)\right)
_{H}+j_{\varepsilon }^{\ast }\left( \zeta _{\varepsilon }(t)\right)
+\int_{t}^{T_{\varepsilon }^{\ast }}\left( Pu_{\varepsilon }^{\ast
}(s),F(h_{\varepsilon }^{\ast }(s))\right) _{U,U^{\ast }}ds,  \notag
\end{eqnarray}%
for all $t\in \lbrack 0,T_{\varepsilon }^{\ast }].$ By comparison with (\ref%
{g59}), we obtain 
\begin{eqnarray*}
&&(A_{H}y_{\varepsilon }^{\ast }(t),p_{\varepsilon }(t))_{H}+j_{\varepsilon
}^{\ast }\left( \zeta _{\varepsilon }(t)\right) +\int_{t}^{T_{\varepsilon
}^{\ast }}\left( Pu_{\varepsilon }^{\ast }(s),F(h_{\varepsilon }^{\ast
}(s))\right) _{U,U^{\ast }}ds \\
&=&1+\frac{1}{2}\left\Vert h_{\varepsilon }^{\ast }(T_{\varepsilon }^{\ast
})\right\Vert _{U}^{2}-\rho \left\Vert B^{\ast }p_{\varepsilon
}(T_{\varepsilon }^{\ast })+\varepsilon F(Pu_{\varepsilon }^{\ast
}(T_{\varepsilon }^{\ast }))\right\Vert _{U^{\ast }}-\frac{\varepsilon }{2}%
\left\Vert Pu_{\varepsilon }^{\ast }(T_{\varepsilon }^{\ast })\right\Vert
_{H}^{2}+j_{\varepsilon }^{\ast }\left( \zeta _{\varepsilon }(T_{\varepsilon
}^{\ast })\right) .
\end{eqnarray*}%
Recalling (\ref{g60-1}), this yields 
\begin{eqnarray*}
&&(A_{H}y_{\varepsilon }^{\ast }(t),p_{\varepsilon }(t))_{H}+\frac{%
\varepsilon }{2}\left\Vert Pu_{\varepsilon }^{\ast }(t)\right\Vert
_{U}^{2}+\rho \left\Vert B^{\ast }p_{\varepsilon
}(t)+\int_{t}^{T_{\varepsilon }^{\ast }}F(h_{\varepsilon }^{\ast }(\tau
))d\tau +\varepsilon F(Pu_{\varepsilon }^{\ast }(t))\right\Vert _{U^{\ast }}
\\
&&+\int_{t}^{T_{\varepsilon }^{\ast }}\left( Pu_{\varepsilon }^{\ast
}(s),F(h_{\varepsilon }^{\ast }(s))\right) _{U,U^{\ast }}ds \\
&=&1+\frac{1}{2}\left\Vert h_{\varepsilon }^{\ast }(T_{\varepsilon }^{\ast
})\right\Vert _{U}^{2}-\rho \left\Vert B^{\ast }p_{\varepsilon
}(T_{\varepsilon }^{\ast })+\varepsilon F(Pu_{\varepsilon }^{\ast
}(T_{\varepsilon }^{\ast }))\right\Vert _{U^{\ast }}-\frac{\varepsilon }{2}%
\left\Vert Pu_{\varepsilon }^{\ast }(T_{\varepsilon }^{\ast })\right\Vert
_{H}^{2} \\
&&+\frac{\varepsilon }{2}\left\Vert Pu_{\varepsilon }^{\ast }(T_{\varepsilon
}^{\ast })\right\Vert _{U}^{2}+\rho \left\Vert B^{\ast }p_{\varepsilon
}(T_{\varepsilon }^{\ast })+\int_{T_{\varepsilon }^{\ast }}^{T_{\varepsilon
}^{\ast }}F(h_{\varepsilon }^{\ast }(\tau ))d\tau +\varepsilon
F(Pu_{\varepsilon }^{\ast }(T_{\varepsilon }^{\ast }))\right\Vert _{U^{\ast
}}
\end{eqnarray*}%
and so we obtain (\ref{g42}), as claimed.\hfill $\square $

\subsection{Optimality conditions for $(P)$}

In order to ensure the passing to the limit in the approximating optimality
conditions (\ref{g41})-(\ref{g42}) we complete the hypotheses $%
(a_{1})-(a_{6}),$ $(b_{1}),$ $(c_{1})$ with $(d_{1})-(d_{5})$.

\medskip

\noindent \textbf{Theorem 5.5. }\textit{Let }%
\begin{equation}
y_{0}\in V,\mbox{ }Py^{tar}\in P(D_{H}),\mbox{ }Py_{0}\neq Py^{tar}.\mbox{ }
\label{g70}
\end{equation}%
\textit{Let }$(T^{\ast },u^{\ast },y^{\ast })$\textit{\ be an optimal pair
in }$(P).$\textit{\ Then, the first order necessary conditions of optimality
are }%
\begin{equation}
Pu^{\ast }(t)\in (N_{K})^{-1}\left( -B^{\ast }p(t)\right) ,%
\mbox{
\textit{a.e.} }t\in (0,T^{\ast }),  \label{g72}
\end{equation}

\begin{equation}
\rho \left\Vert B^{\ast }p(t)\right\Vert _{U^{\ast }}+(A_{H}y^{\ast
}(t),p(t))_{H}=1,\mbox{ \textit{a.e.} }t\in (0,T^{\ast }),\mbox{ }
\label{g73}
\end{equation}%
\textit{where }$y^{\ast }$ \textit{is the solution to\ the state system} (%
\ref{g4})-(\ref{g4-1}) \textit{corresponding to }$(T^{\ast },u^{\ast }),$%
\textit{\ and }$p$ \textit{is a solution to } 
\begin{equation}
-p^{\prime }(t)+(A^{\prime }(y^{\ast }(t)))^{\ast }p(t)=0,%
\mbox{
\textit{a.e. }}t\in (0,T^{\ast }),  \label{g75}
\end{equation}%
\begin{equation}
p(T^{\ast })\in V^{\ast }.  \label{g75-0}
\end{equation}

\medskip

\noindent \textbf{Proof. }First, we prove that 
\begin{equation}
\left\Vert B^{\ast }p_{\varepsilon }(T_{\varepsilon }^{\ast })\right\Vert
_{U^{\ast }}\leq C,  \label{g76}
\end{equation}%
with\textbf{\ }$C$ independent of $\varepsilon .$

Let us begin with the case $P(y_{1},y_{2})=(y_{1},0),$ $%
B(u_{1},u_{2})=(u_{1},0).$ We recall that in this case $y^{tar}=(y_{1}^{%
\mbox{target}},z),$ $\forall z\in H,$ and $Py^{tar}=(y_{1}^{\mbox{target}%
},0).$ We start from (\ref{g59}) and express the third term on the left-hand
side as 
\begin{eqnarray*}
&&\left\langle Ay_{\varepsilon }^{\ast }(T_{\varepsilon }^{\ast
}),p_{\varepsilon }(T_{\varepsilon }^{\ast })\right\rangle _{V^{\ast },V}=%
\frac{1}{\varepsilon }\left\langle Ay_{\varepsilon }^{\ast }(T_{\varepsilon
}^{\ast }),Py_{\varepsilon }^{\ast }(T_{\varepsilon }^{\ast
})-Py^{tar}\right\rangle _{V^{\ast },V} \\
&=&\frac{1}{\varepsilon }\left\langle Ay_{\varepsilon }^{\ast
}(T_{\varepsilon }^{\ast })-A\widehat{y},Py_{\varepsilon }^{\ast
}(T_{\varepsilon }^{\ast })-Py^{tar}\right\rangle _{V^{\ast },V}+\frac{1}{%
\varepsilon }\left\langle A\widehat{y},Py_{\varepsilon }^{\ast
}(T_{\varepsilon }^{\ast })-Py^{tar}\right\rangle _{V^{\ast },V},
\end{eqnarray*}%
where $\widehat{y}=(y_{1}^{\mbox{target}},\widehat{z})$ set by (\ref{y^}),
satisfying (\ref{g5-000}), with the choice (\ref{Ay^}). We note that 
\begin{equation}
Py_{\varepsilon }^{\ast }(T_{\varepsilon }^{\ast
})-Py^{tar}=(Py_{\varepsilon }^{\ast }(T_{\varepsilon }^{\ast })-y_{1}^{%
\mbox{target}},0)=Py_{\varepsilon }(T_{\varepsilon }^{\ast })-P\widehat{y}.
\label{1010}
\end{equation}%
Then, by (\ref{g59}), we can write 
\begin{eqnarray}
&&\rho \left\Vert B^{\ast }p_{\varepsilon }(T_{\varepsilon }^{\ast
})\right\Vert _{U^{\ast }}+\left\langle Ay_{\varepsilon }^{\ast
}(T_{\varepsilon }^{\ast })-A\widehat{y},p_{\varepsilon }(T_{\varepsilon
}^{\ast })\right\rangle _{V^{\ast },V}  \label{1005} \\
&\leq &1+\frac{1}{2}\left\Vert h_{\varepsilon }^{\ast }(T_{\varepsilon
}^{\ast })\right\Vert _{U}^{2}+\varepsilon \rho \left\Vert F(Pu_{\varepsilon
}^{\ast }(T_{\varepsilon }^{\ast }))\right\Vert _{U^{\ast }}+\left\vert
(A_{H}\widehat{y},p_{\varepsilon }(T_{\varepsilon }^{\ast }))_{H}\right\vert 
\notag \\
&\leq &1+\left\Vert h_{\varepsilon }^{\ast }(T_{\varepsilon }^{\ast
})\right\Vert _{U}^{2}+\varepsilon \rho ^{2}+\left\vert \left\langle
Pp_{\varepsilon }(T_{\varepsilon }^{\ast }),A_{H}\widehat{y}\right\rangle
_{V^{\ast },V}\right\vert \leq 1+\left\Vert h_{\varepsilon }^{\ast
}(T_{\varepsilon }^{\ast })\right\Vert _{U}^{2}+\varepsilon \rho
^{2}+\left\Vert Pp_{\varepsilon }(T_{\varepsilon }^{\ast })\right\Vert
_{V^{\ast }}\left\Vert A_{H}\widehat{y}\right\Vert _{V}.  \notag
\end{eqnarray}%
Here, we took into account that 
\begin{equation}
p_{\varepsilon }(T_{\varepsilon }^{\ast })=\frac{1}{\varepsilon }%
(Py_{\varepsilon }(T_{\varepsilon }^{\ast })-Py^{tar})=\frac{1}{\varepsilon }%
P(Py_{\varepsilon }(T_{\varepsilon }^{\ast })-Py^{tar})=Pp_{\varepsilon
}(T_{\varepsilon }^{\ast }).  \label{1000}
\end{equation}%
Now, we use (\ref{g5-000}) which is assumed to take place for $t\in
(0,T_{\ast }+\delta )$, with $T_{\ast }$ the time specified in the
controllability hypothesis $(c_{1}),$ with $T_{\ast }\geq T^{\ast }.$ Recall
that $T_{\varepsilon }^{\ast }\rightarrow T^{\ast }.$ Hence, for $%
\varepsilon $ sufficiently small, $T_{\varepsilon }^{\ast }\in (0,T^{\ast
}+\delta )\subset (0,T_{\ast }+\delta ),$ with $\delta $ arbitrary small and
it follows that relation (\ref{g5-000}) can take place also for $%
t=T_{\varepsilon }^{\ast }$, that is 
\begin{eqnarray*}
&&\left\langle Ay_{\varepsilon }^{\ast }(T_{\varepsilon }^{\ast })-A\widehat{%
y},p_{\varepsilon }(T_{\varepsilon }^{\ast })\right\rangle _{V^{\ast },V}=%
\frac{1}{\varepsilon }\left\langle Ay_{\varepsilon }^{\ast }(T_{\varepsilon
}^{\ast })-A\widehat{y},P(y_{\varepsilon }(T_{\varepsilon }^{\ast })-%
\widehat{y})\right\rangle _{V^{\ast },V} \\
&\geq &-\frac{C_{3}}{\varepsilon }\left\Vert Py_{\varepsilon
}(T_{\varepsilon }^{\ast })-P\widehat{y}\right\Vert _{P(H)}^{2}=-\frac{C_{3}%
}{\varepsilon }\left\Vert Py_{\varepsilon }(T_{\varepsilon }^{\ast })-y_{1}^{%
\mbox{target}}\right\Vert _{P(H)}^{2}.
\end{eqnarray*}%
Here we used (\ref{1010}). Then, 
\begin{eqnarray*}
&&\rho \left\Vert B^{\ast }p_{\varepsilon }(T_{\varepsilon }^{\ast
})\right\Vert _{U^{\ast }}\leq 1+\left\Vert h_{\varepsilon }^{\ast
}(T_{\varepsilon }^{\ast })\right\Vert _{U}^{2}+\varepsilon C\rho
^{2}+\left\Vert Pp_{\varepsilon }(T_{\varepsilon }^{\ast })\right\Vert
_{V^{\ast }}\left\Vert A_{H}\widehat{y}\right\Vert _{V}+C_{3}\left\Vert
p_{\varepsilon }(T_{\varepsilon }^{\ast })\right\Vert _{H}^{2} \\
&\leq &1+\left\Vert h_{\varepsilon }^{\ast }(T_{\varepsilon }^{\ast
})\right\Vert _{U}^{2}+\varepsilon \rho ^{2}+C^{\ast }\left\Vert B^{\ast
}p_{\varepsilon }(T_{\varepsilon }^{\ast })\right\Vert _{U^{\ast
}}\left\Vert A_{H}\widehat{y}\right\Vert _{V}-\frac{C_{3}}{\varepsilon }%
\left\Vert Py_{\varepsilon }(T_{\varepsilon }^{\ast })-Py^{tar}\right\Vert
_{H}^{2},
\end{eqnarray*}%
where we took into account (\ref{g74-2}). We recall (\ref{g36-0}), 
\begin{equation*}
\limsup_{\varepsilon \rightarrow 0}\frac{1}{2\varepsilon }\left\Vert
Py_{\varepsilon }^{\ast }(T_{\varepsilon }^{\ast })-Py^{tar}\right\Vert
_{H}^{2}=0,
\end{equation*}%
and so, by (\ref{Ay^}), we can write%
\begin{equation*}
\rho \left\Vert B^{\ast }p_{\varepsilon }(T_{\varepsilon }^{\ast
})\right\Vert _{U^{\ast }}<1+C+\rho _{1}\left\Vert B^{\ast }p_{\varepsilon
}(T_{\varepsilon }^{\ast })\right\Vert _{U^{\ast }},
\end{equation*}%
because $\left\Vert h_{\varepsilon }^{\ast }(T_{\varepsilon }^{\ast
})\right\Vert _{U}^{2}+\varepsilon C\rho ^{2}\rightarrow 0.$ The convergence
of the first term is due to (\ref{g34-0}). This yields 
\begin{equation*}
(\rho -\rho _{1})\left\Vert B^{\ast }p_{\varepsilon }(T_{\varepsilon }^{\ast
})\right\Vert _{U^{\ast }}<1+C,
\end{equation*}%
and choosing $\rho >\rho _{1}$ we finally get (\ref{g76}).

As a matter of fact, in the proof of (\ref{g76}) the second component $%
\widehat{z}$ of $\widehat{y}$ can be generally set as the second component
of the approximating state solution.

If $P=I,$ we proceed in the same way, and use that $y^{tar}=(y_{1}^{%
\mbox{target}},y_{2}^{\mbox{target}})$ and $(d_{5})\ $and take $Ay^{tar}$
instead of $A\widehat{y}.$ We have%
\begin{eqnarray*}
&&\left\langle Ay_{\varepsilon }^{\ast }(T_{\varepsilon }^{\ast
})-Ay^{tar},p_{\varepsilon }(T_{\varepsilon }^{\ast })\right\rangle
_{V^{\ast },V}=\frac{1}{\varepsilon }\left\langle Ay_{\varepsilon }^{\ast
}(T_{\varepsilon }^{\ast })-Ay^{tar},y_{\varepsilon }(T_{\varepsilon }^{\ast
})-y^{tar}\right\rangle _{V^{\ast },V} \\
&\geq &-\frac{C_{3}}{\varepsilon }\left\Vert y_{\varepsilon }(T_{\varepsilon
}^{\ast })-y^{tar}\right\Vert _{H}^{2}=-\frac{C_{3}}{\varepsilon }\left\Vert
Py_{\varepsilon }(T_{\varepsilon }^{\ast })-Py^{tar}\right\Vert _{H}^{2}
\end{eqnarray*}%
which tends to zero by (\ref{g36-0}). Here, $Py_{\varepsilon
}(T_{\varepsilon }^{\ast })-Py^{tar}$ has both nonzero components.

We recall the adjoint system given by (\ref{g43})-(\ref{g44}). Since the
final data is bounded in $U^{\ast },$ according to (\ref{g76}), we expect to
obtain at limit a solution with a weaker regularity. We are going to obtain
some uniform estimates for the solution $p_{\varepsilon }.$

\noindent A first estimate is obtained by multiplying scalarly (\ref{g43})
by $\Gamma _{H}^{-1}p_{\varepsilon }(t)$ and integrating from $t$ to $%
T_{\varepsilon }^{\ast }$ 
\begin{equation*}
\frac{1}{2}\left\Vert p_{\varepsilon }(t)\right\Vert _{V^{\ast
}}^{2}+\int_{t}^{T_{\varepsilon }^{\ast }}\left( (A_{H}^{\prime
}(y_{\varepsilon }^{\ast }(\tau )))^{\ast }p_{\varepsilon }(\tau ),\Gamma
_{H}^{-1}p_{\varepsilon }(\tau )\right) _{H}d\tau =\frac{1}{2}\left\Vert
p_{\varepsilon }(T_{\varepsilon }^{\ast })\right\Vert _{V^{\ast }}^{2},
\end{equation*}%
where the first term on the right-hand side was obtained by using the
properties of the duality mapping (\ref{Gamma1}). According to (\ref{g74-10}%
) and (\ref{g74-2}) for $v=p_{\varepsilon }(T_{\varepsilon }^{\ast })\in
V^{\ast },$ we successively get%
\begin{eqnarray*}
&&\frac{1}{2}\left\Vert p_{\varepsilon }(t)\right\Vert _{V^{\ast
}}^{2}+C_{1}\int_{t}^{T_{\varepsilon }^{\ast }}\left\Vert p_{\varepsilon
}(\tau )\right\Vert _{H}^{2}d\tau \leq \frac{1}{2}\left\Vert p_{\varepsilon
}(T_{\varepsilon }^{\ast })\right\Vert _{V^{\ast
}}^{2}+C_{2}\int_{t}^{T_{\varepsilon }^{\ast }}\left\Vert p_{\varepsilon
}(\tau )\right\Vert _{V^{\ast }}^{2}(1+C_{3}\left\Vert y_{\varepsilon
}^{\ast }(\tau )\right\Vert _{V}^{l})d\tau \\
&\leq &\frac{1}{2}\left\Vert B^{\ast }p_{\varepsilon }(T_{\varepsilon
}^{\ast })\right\Vert _{U^{\ast }}^{2}+C_{4}\int_{t}^{T_{\varepsilon }^{\ast
}}\left\Vert p_{\varepsilon }(\tau )\right\Vert _{V^{\ast }}^{2}d\tau .
\end{eqnarray*}%
Here, we used (\ref{g15}). By Gronwall lemma and (\ref{g76}) we obtain 
\begin{equation}
\left\Vert p_{\varepsilon }(t)\right\Vert _{V^{\ast
}}^{2}+\int_{t}^{T_{\varepsilon }^{\ast }}\left\Vert p_{\varepsilon }(\tau
)\right\Vert _{H}^{2}d\tau \leq C,\mbox{ for all }t\in \lbrack
0,T_{\varepsilon }^{\ast }],  \label{g77}
\end{equation}%
independently on $\varepsilon .$

Next, we multiply scalarly (\ref{g43}) by $\Gamma _{H}^{-\alpha
}p_{\varepsilon }(t)$ (where $\alpha $ is chosen by (\ref{g74}) and (\ref%
{g74-1})) and integrate from $t$ to $T_{\varepsilon }^{\ast }.$ Applying (%
\ref{g74-1}) and (\ref{g74}) we obtain 
\begin{eqnarray}
&&\frac{1}{2}\left\Vert \Gamma ^{-\alpha /2}p_{\varepsilon }(t)\right\Vert
_{H}^{2}+C_{1}\int_{t}^{T_{\varepsilon }^{\ast }}\left\Vert \Gamma
_{H}^{(1-\alpha )/2}p_{\varepsilon }(\tau )\right\Vert _{H}^{2}d\tau
\label{g82} \\
&\leq &C_{2}\int_{t}^{T_{\varepsilon }^{\ast }}\left\Vert p_{\varepsilon
}(\tau )\right\Vert _{H}^{2}(1+C_{3}\left\Vert y_{\varepsilon }^{\ast }(\tau
)\right\Vert _{V}^{l})+\frac{1}{2}\left\Vert B^{\ast }p_{\varepsilon
}(T_{\varepsilon }^{\ast })\right\Vert _{U^{\ast }}^{2}.  \notag
\end{eqnarray}%
Recalling (\ref{g77}) we obtain $\left\Vert \Gamma _{H}^{(1-\alpha
)/2}p_{\varepsilon }\right\Vert _{L^{2}(0,T_{\varepsilon }^{\ast };H)}\leq
C, $ that is 
\begin{equation}
\left\Vert p_{\varepsilon }\right\Vert _{L^{2}\left( 0,T_{\varepsilon
}^{\ast };D\left( \Gamma _{H}^{(1-\alpha )/2}\right) \right) }\leq C.
\label{g83}
\end{equation}%
To use further these estimate we have to modify the functional framework in
the following sense. We extend the operator $(A_{H}^{\prime }(y_{\varepsilon
}^{\ast }(t)))^{\ast }$ to $H,$ for all $t\in \lbrack 0,T_{\varepsilon
}^{\ast }],$ namely we define $\tilde{A}_{H}^{\prime }(y_{\varepsilon
}^{\ast }(t)):H\subset D_{H}^{\ast }\rightarrow D_{H}^{\ast }$ by 
\begin{equation*}
\left\langle \tilde{A}_{H}^{\prime }(y_{\varepsilon }^{\ast }(t))v,\psi
\right\rangle _{D_{H}^{\ast },D_{H}}=(v,(A_{H}^{\prime }(y_{\varepsilon
}^{\ast }(t)))\psi )_{H},\mbox{ for }v\in H,\mbox{ }\psi \in D_{H},%
\mbox{
for all }t\in \lbrack 0,T_{\varepsilon }^{\ast }].
\end{equation*}%
The norm on $D_{H}^{\ast }$ is defined by $\left\Vert \theta \right\Vert
_{D_{H}^{\ast }}=\left\Vert (A_{H}^{\prime }(y_{\varepsilon }^{\ast
}(t)))^{-1}\theta \right\Vert _{H}.$ We have, by (\ref{g7-1}) 
\begin{eqnarray*}
&&\left\vert \left\langle \tilde{A}_{H}^{\prime }((y_{\varepsilon }^{\ast
}(t))v,\psi \right\rangle _{D_{H}^{\ast },D_{H}}\right\vert =\left\vert
\left( v,A_{H}^{\prime }((y_{\varepsilon }^{\ast }(t))\psi \right)
_{H}\right\vert \leq \left\Vert v\right\Vert _{H}\left\Vert A_{H}^{\prime
}((y_{\varepsilon }^{\ast }(t))\psi \right\Vert _{H} \\
&\leq &C\left\Vert v\right\Vert _{H}\left\Vert \psi \right\Vert
_{D_{H}}(1+\left\Vert y_{\varepsilon }^{\ast }(t)\right\Vert _{V}^{\kappa
})\leq C\left\Vert v\right\Vert _{H}\left\Vert \psi \right\Vert _{D_{H}},
\end{eqnarray*}%
which yields for $v=p_{\varepsilon }$, 
\begin{equation*}
\left\Vert \tilde{A}_{H}^{\prime }(y_{\varepsilon }^{\ast
}(t))p_{\varepsilon }\right\Vert _{L^{2}(0,T_{\varepsilon }^{\ast
};D_{H}^{\ast })}\leq C\left\Vert p_{\varepsilon }\right\Vert
_{L^{2}(0,T_{\varepsilon }^{\ast };H)}\leq C,
\end{equation*}%
since $\left\Vert y_{\varepsilon }^{\ast }(t)\right\Vert _{V}\leq C.$ By
comparison in the adjoint equation (\ref{g43}) we obtain 
\begin{equation}
\left\Vert p_{\varepsilon }^{\prime }\right\Vert _{L^{2}(0,T_{\varepsilon
}^{\ast };D_{H}^{\ast })}\leq C.  \label{g83-1}
\end{equation}%
We recall that $T_{\varepsilon }^{\ast }\rightarrow T^{\ast }$ and so $%
T^{\ast }-\delta \leq T_{\varepsilon }^{\ast },$ with $\delta $ arbitrary,
so that the estimates are true also on $(0,T^{\ast }-\delta ).$ By (\ref{g77}%
) and the latter, selecting a subsequence, denoted still by $\varepsilon ,$
we have 
\begin{equation*}
p_{\varepsilon }\rightarrow p\mbox{ weakly in }L^{2}(0,T^{\ast }-\delta ;H),%
\mbox{ weak-star in }L^{\infty }(0,T^{\ast }-\delta ;V^{\ast }),
\end{equation*}%
\begin{equation*}
p_{\varepsilon }^{\prime }\rightarrow p^{\prime }\mbox{ weakly in }%
L^{2}(0,T^{\ast }-\delta ;D_{H}^{\ast }).
\end{equation*}%
Because $\delta $ is arbitrary, the previous convergences take place on $%
\cap_{\delta >0}^{} (0,T^{\ast }-\delta )=(0,T^{\ast })$.

Since $D\left( \Gamma _{H}^{(1-\alpha )/2}\right) $ is compact in $H$ and $%
H\subset D_{H}^{\ast },$ we have by Aubin-Lions lemma that 
\begin{equation}
p_{\varepsilon }\rightarrow p\mbox{ strongly in }L^{2}(0,T^{\ast };H).
\label{g84}
\end{equation}%
Then, using the convergence $y_{\varepsilon }^{\ast }(t)\rightarrow y^{\ast
}(t)$ strongly in $V$ a.e. $t,$ and the continuity (\ref{g6}), we have%
\begin{eqnarray*}
&&\left\langle \tilde{A}_{H}^{\prime }(y_{\varepsilon }^{\ast
})p_{\varepsilon }(t)-\tilde{A}_{H}^{\prime }(y^{\ast })p(t),\psi
(t)\right\rangle _{D_{H}^{\ast },D_{H}} \\
&=&\left\langle (\tilde{A}_{H}^{\prime }(y_{\varepsilon }^{\ast })-\tilde{A}%
_{H}^{\prime }(y^{\ast }))p_{\varepsilon }(t),\psi (t)\right\rangle
_{D_{H}^{\ast },D_{H}}+\left\langle (\tilde{A}_{H}^{\prime }(y^{\ast
})(p_{\varepsilon }(t)-p(t)),\psi (t)\right\rangle _{D_{H}^{\ast },D_{H}} \\
&=&\left\langle p_{\varepsilon }(t),(A_{H}^{\prime }(y_{\varepsilon }^{\ast
})-A_{H}^{\prime }(y^{\ast }))\psi (t)\right\rangle _{H}+\left\langle
p_{\varepsilon }(t)-p(t),A_{H}^{\prime }(y^{\ast })\psi (t)\right\rangle
_{H},
\end{eqnarray*}%
for $\psi \in L^{2}(0,T_{\varepsilon }^{\ast };D_{H})$, which implies by the
previous convergences that 
\begin{equation*}
\tilde{A}_{H}^{\prime }(y_{\varepsilon }^{\ast })p_{\varepsilon }\rightarrow 
\tilde{A}_{H}^{\prime }(y^{\ast })p\mbox{ weakly in }L^{2}(0,T^{\ast
};D_{H}^{\ast }),\mbox{ as }\varepsilon \rightarrow 0.
\end{equation*}%
We also have 
\begin{equation}
B^{\ast }p_{\varepsilon }\rightarrow B^{\ast }p\mbox{ strongly in }%
L^{2}(0,T^{\ast };U^{\ast }).  \label{g85}
\end{equation}%
By these convergences we obtain (\ref{g75}) in the sense of distributions
and a.e.

We go back now to (\ref{g41}), and recall (\ref{g55}) which can be
equivalently written 
\begin{equation*}
N_{\mathcal{K}_{T_{\varepsilon }^{\ast }}}(Pu_{\varepsilon }^{\ast })\ni
z_{\varepsilon }=-B^{\ast }p_{\varepsilon }-\varepsilon F(Pu_{\varepsilon
}^{\ast })-\int_{\cdot }^{T_{\varepsilon }^{\ast }}F(h_{\varepsilon }^{\ast
}(\tau ))d\tau .
\end{equation*}%
We pass to the limit as $\varepsilon \rightarrow 0$ and have 
\begin{equation*}
u_{\varepsilon }^{\ast }\rightarrow u^{\ast }\mbox{ weakly in }%
L^{2}(0,T^{\ast };U),\mbox{ }B^{\ast }p_{\varepsilon }\rightarrow B^{\ast }p%
\mbox{ strongly in }L^{2}(0,T^{\ast };U^{\ast }),
\end{equation*}%
\begin{equation*}
F(h_{\varepsilon }^{\ast })\rightarrow 0\mbox{ strongly in }L^{2}(0,T^{\ast
};U^{\ast }),
\end{equation*}%
by (\ref{g34}), and 
\begin{equation*}
\int_{t}^{T_{\varepsilon }^{\ast }}F(h_{\varepsilon }^{\ast }(\tau ))d\tau
\rightarrow 0,\mbox{ strongly in }U^{\ast },\mbox{ for all }t\in
(0,T_{\varepsilon }^{\ast }).
\end{equation*}%
Therefore, 
\begin{equation*}
z_{\varepsilon }=-B^{\ast }p_{\varepsilon }-\varepsilon F(Pu_{\varepsilon
}^{\ast })-\int_{\cdot }^{T_{\varepsilon }^{\ast }}F(h_{\varepsilon }^{\ast
}(\tau ))d\tau \rightarrow -B^{\ast }p\mbox{ strongly in }L^{2}(0,T^{\ast
};U^{\ast }).
\end{equation*}%
But\ $N_{\mathcal{K}_{T^{\ast }}}$ is maximal monotone from $L^{2}(0,T^{\ast
};U)$ to $L^{2}(0,T^{\ast };U^{\ast })$, that is weakly-strongly closed and
since $Pu_{\varepsilon }^{\ast }\rightarrow Pu^{\ast }$ weakly in $%
L^{2}(0,T^{\ast };U),$ we get $-B^{\ast }p\in N_{\mathcal{K}_{T^{\ast
}}}(Pu^{\ast }),$ or equivalently (\ref{g72}).

Finally, we have to pass to the limit in (\ref{g42}). For this, we integrate
(\ref{g42}) from $s$ to $s^{\prime },$ $0<s<s^{\prime }<T^{\ast }$ and get%
\begin{eqnarray*}
&&\rho \int_{s}^{s^{\prime }}\left\Vert B^{\ast }p_{\varepsilon
}(t)+\int_{t}^{T_{\varepsilon }^{\ast }}F(h_{\varepsilon }^{\ast }(\tau
))d\tau +\varepsilon F(Pu_{\varepsilon }^{\ast }(t))\right\Vert _{U^{\ast
}}dt+\int_{s}^{s^{\prime }}\left( Ay_{\varepsilon }^{\ast
}(t),p_{\varepsilon }(t)\right) _{H}dt \\
&&+\int_{s}^{s^{\prime }}\int_{t}^{T_{\varepsilon }^{\ast }}\left\langle
Pu_{\varepsilon }^{\ast }(\tau ),F(h_{\varepsilon }^{\ast }(\tau
))\right\rangle _{U,U^{\ast }}d\tau dt+\frac{\varepsilon }{2}%
\int_{s}^{s^{\prime }}\left\Vert Pu_{\varepsilon }^{\ast }(t)\right\Vert
_{H}^{2}dt \\
&=&(s^{\prime }-s)\left( 1+\frac{1}{2}\left\Vert h_{\varepsilon }^{\ast
}(T_{\varepsilon }^{\ast })\right\Vert _{U}^{2}\right) .
\end{eqnarray*}%
We recall that $Ay_{\varepsilon }^{\ast }\rightarrow Ay^{\ast }$ weakly in $%
L^{2}(0,T^{\ast };H)$ and note that 
\begin{equation*}
\int_{t}^{T_{\varepsilon }^{\ast }}\left\langle Pu_{\varepsilon }^{\ast
}(\tau ),F(h_{\varepsilon }^{\ast }(\tau ))\right\rangle _{U,U^{\ast }}d\tau
\rightarrow 0,\mbox{ for all }t\in (0,T_{\varepsilon }^{\ast }).
\end{equation*}%
Finally, $\left\Vert h_{\varepsilon }^{\ast }(T_{\varepsilon }^{\ast
})\right\Vert _{U}^{2}\rightarrow 0,$ by (\ref{g34-0}). We pass to the limit
as $\varepsilon $ goes to 0 and get%
\begin{equation*}
\int_{s}^{s^{\prime }}\left\{ \rho \left\Vert B^{\ast }p(t)\right\Vert
_{U^{\ast }}+\left( Ay^{\ast }(t),p(t)\right) _{H}\right\} dt=(s^{\prime
}-s).
\end{equation*}%
Dividing by $(s^{\prime }-s)$ and passing to the limit as $s\rightarrow
s^{\prime }$ we obtain (\ref{g73}), for a.e. $t\in (0,T^{\ast })$. \hfill $%
\square $

\bigskip

In the case when $U=U^{\ast }=H$ we have a particular result for which we
assume the hypotheses $(a_{1})-(a_{6})$ and replace $(d_{1})-(d_{5})$ by
simpler ones.

\medskip

\noindent \textbf{Corollary 5.6. }\textit{Let }$U=U^{\ast }=H$ \textit{and
assume} (\ref{g70}), (\ref{g5-000}), \textit{and}%
\begin{equation}
((A_{H}^{\prime }(y))^{\ast }v,v)_{H}\geq C_{1}\left\Vert v\right\Vert
_{H}^{2}-C_{2}\left\Vert v\right\Vert _{V}^{2}(1+C_{3}\left\Vert
y\right\Vert _{V}^{l}),\mbox{ \textit{for all} }y,v\in V,\mbox{ }l\geq 0,
\label{1001}
\end{equation}%
\begin{equation}
\left\Vert Pv\right\Vert _{H}\leq C^{\ast }\left\Vert B^{\ast }v\right\Vert
_{H},\mbox{ \textit{for} }v\in H,  \label{1004}
\end{equation}%
\begin{equation}
\rho >\rho _{1},\mbox{ }\rho _{1}:=C^{\ast }\left\Vert A_{H}\widehat{y}%
\right\Vert _{H},  \label{1002}
\end{equation}%
(\textit{instead of} (\ref{g74-2})). \textit{Then}, (\ref{g72})-(\ref{g75}) 
\textit{take place and }$p(T^{\ast })\in H.$

\medskip

\noindent \textbf{Proof. }We resume the proof of the estimate for $%
\left\Vert B^{\ast }p_{\varepsilon }(T_{\varepsilon }^{\ast })\right\Vert
_{U^{\ast }}$ in Theorem 5.5 and have now in (\ref{1005}) 
\begin{eqnarray*}
&&\rho \left\Vert B^{\ast }p_{\varepsilon }(T_{\varepsilon }^{\ast
})\right\Vert _{U^{\ast }}+\left\langle Ay_{\varepsilon }^{\ast
}(T_{\varepsilon }^{\ast })-A\widehat{y},p_{\varepsilon }(T_{\varepsilon
}^{\ast })\right\rangle _{V^{\ast },V} \\
&\leq &1+\left\Vert h_{\varepsilon }^{\ast }(T_{\varepsilon }^{\ast
})\right\Vert _{U}^{2}+\varepsilon \rho ^{2}+\left\vert (A_{H}\widehat{y}%
,p_{\varepsilon }(T_{\varepsilon }^{\ast }))_{H}\right\vert \\
&\leq &1+C+\left\Vert Pp_{\varepsilon }(T_{\varepsilon }^{\ast })\right\Vert
_{H}\left\Vert A_{H}\widehat{y}\right\Vert _{H}<C_{1}+\rho _{1}\left\Vert
B^{\ast }p_{\varepsilon }(T_{\varepsilon }^{\ast })\right\Vert _{H}.
\end{eqnarray*}%
Since $U^{\ast }=H$ we get $\left\Vert B^{\ast }p_{\varepsilon
}(T_{\varepsilon }^{\ast })\right\Vert _{H}\leq C,$ which will ensure a more
regular solution for $p.$ We multiply (\ref{g43}) by $p_{\varepsilon }(t),$
integrate from $t$ to $T_{\varepsilon }^{\ast }$ and use (\ref{1001}) to
obtain 
\begin{equation}
\left\Vert p_{\varepsilon }(t)\right\Vert _{H}^{2}+\int_{t}^{T_{\varepsilon
}^{\ast }}\left\Vert p_{\varepsilon }(\tau )\right\Vert _{V}^{2}d\tau \leq C,%
\mbox{ for all }t\in \lbrack 0,T_{\varepsilon }^{\ast }].  \label{1006}
\end{equation}%
Then,%
\begin{equation*}
\int_{0}^{T_{\varepsilon }^{\ast }}\left\Vert (A_{H}^{\prime
}(y_{\varepsilon }^{\ast }))^{\ast }p_{\varepsilon }(t)\right\Vert _{V^{\ast
}}^{2}dt\leq \int_{0}^{T_{\varepsilon }^{\ast }}\left\Vert p_{\varepsilon
}(t)\right\Vert _{V}^{2}(1+C_{3}\left\Vert y_{\varepsilon }^{\ast
}(t)\right\Vert _{V}^{l})dt\leq C
\end{equation*}%
and by (\ref{g43}) we infer that 
\begin{equation*}
\int_{0}^{T_{\varepsilon }^{\ast }}\left\Vert p_{\varepsilon }^{\prime
}(t)\right\Vert _{V^{\ast }}^{2}dt\leq C.
\end{equation*}%
On a subsequence we obtain 
\begin{eqnarray*}
p_{\varepsilon } &\rightarrow &p\mbox{ weakly in }L^{2}(0,T^{\ast };V)\cap
W^{1,2}(0,T^{\ast };V^{\ast }),\mbox{ weak-star in }L^{\infty }(0,T^{\ast
};H), \\
&&\mbox{ strongly in }L^{2}(0,T^{\ast };H),
\end{eqnarray*}%
where $p$ turns out to be the solution to (\ref{g75}). The rest of the proof
can be led as in Theorem 5.5.\hfill $\square $

\medskip

\noindent \textbf{Remark 5.7.} Consider the case $P=I$ and $y^{tar}=0$ and
assume that, for each $\widetilde{y}\in C([0,T];H)\cap L^{2}(0,T;V)$, the
linearized problem 
\begin{equation*}
Y^{\prime }(t)+A^{\prime }(\widetilde{y})Y(t)=Bv(t),\mbox{ a.e. }t>0,\mbox{ }%
Y(0)=Y_{0}
\end{equation*}%
is exactly null controllable in the following sense: for each $Y_{0}\in H$
with $\left\Vert Y_{0}\right\Vert _{H}\leq 1,$ there is $v\in L^{2}(0,T;U),$
with $\left\Vert v\right\Vert _{L^{2}(0,T;U)}\leq \gamma ,$ such that $%
Y(T)=0.$ Then, Theorem 5.5 remains true, without assuming $(d_{1})-(d_{5}).$
Here there is the argument. By the above controllability hypothesis, we get
for the dual equation, $-p^{\prime }(t)+(A^{\prime }(\widetilde{y}))^{\ast
}p(t)=0,$ the following observability inequality:%
\begin{equation*}
\left\Vert p(0)\right\Vert _{H}\leq \gamma \left( \int_{0}^{T_{\varepsilon
}^{\ast }}\left\Vert B^{\ast }p(t)\right\Vert _{U^{\ast }}^{2}dt\right)
^{1/2}
\end{equation*}%
and therefore 
\begin{equation*}
\left\Vert p(t)\right\Vert _{H}\leq \gamma \left( \int_{t}^{T_{\varepsilon
}^{\ast }}\left\Vert B^{\ast }p(\tau )\right\Vert _{U^{\ast }}^{2}d\tau
\right) ^{1/2}.
\end{equation*}%
Then, substituting in (\ref{g42}) we get $\left\Vert p_{\varepsilon
}(t)\right\Vert _{H}\leq \gamma _{1}\int_{t}^{T_{\varepsilon }^{\ast
}}\left\Vert B^{\ast }p_{\varepsilon }(\tau )\right\Vert _{U}^{2}d\tau
+C(\varepsilon )+C\leq C_{t},$ for all $t\in \lbrack 0,T_{\varepsilon
}^{\ast }]$ and then $\int_{0}^{T_{\varepsilon }^{\ast }-\delta }\left\Vert
p_{\varepsilon }(t)\right\Vert _{H}^{2}dt\leq C_{\delta },$ $\delta >0.$
Thus, we may pass to the limit in (\ref{g43})-(\ref{g44}) to get (\ref{g72}%
)-(\ref{g75}).

\section{Examples}

\setcounter{equation}{0}

We particularize our results to some equations and systems modelling various
processes in physical applications. Let $\Omega $ be an open bounded subset
of $\mathbb{R}^{d},$ $d\leq 3$, with a sufficient regular boundary $\partial
\Omega $ and let $\nu $ be the outward normal to $\partial \Omega .$ Let $%
L^{r}(\Omega )$ be the space of $r$-summable functions, $y:\Omega
\rightarrow \mathbb{R},$ with the norm $\left\Vert y\right\Vert _{r}=\left(
\int_{\Omega }\left\vert y\right\vert ^{r}dx\right) ^{1/r}$, $1\leq r<\infty
,$ and $\left\Vert y\right\Vert _{\infty }=ess\sup_{x\in \Omega }\left\vert
y(x)\right\vert $ for $r=\infty .$ The spaces $H^{r}(\Omega
):=W^{1,r}(\Omega ),$ with $1\leq r<\infty ,$ and $H_{0}^{1}(\Omega )$ are
the standard Sobolev spaces, and $H^{-1}(\Omega )$ is the dual of $%
H_{0}^{1}(\Omega ).$

\paragraph{Example 1. Diffusion equation with a potential and drift term.}

Let us consider the problem%
\begin{eqnarray}
y_{t}-\Delta y+\beta (y)+a_{1}y-\nabla \cdot (by) &=&u\mbox{, \ in }%
(0,\infty )\times \Omega ,  \label{ex-1} \\
y(0) &=&y_{0},\mbox{ in }\Omega ,  \notag \\
\frac{\partial y}{\partial \nu }+\gamma y &=&0,\mbox{ \ on }(0,\infty
)\times \partial \Omega ,  \notag
\end{eqnarray}%
where 
\begin{eqnarray}
\beta &:&\mathbb{R}\rightarrow \mathbb{R}\mbox{, }\beta \in C^{1}(\mathbb{R})%
\mbox{, }\beta (0)=0,  \notag \\
0 &<&a_{0}\leq \beta ^{\prime }(r)\leq L(\left\vert r\right\vert ^{\kappa
}+1),\mbox{ for all }r\in \mathbb{R},\mbox{ }\kappa \in \left[ 0,2\right] ,
\label{100}
\end{eqnarray}%
\begin{equation}
a_{1}\in L^{\infty }(\Omega ),\mbox{ }b\in (W^{1,\infty }(\Omega ))^{3},%
\mbox{ }b\cdot \nu =0\mbox{ on }\partial \Omega ,\mbox{ }\gamma \in
L^{\infty }(\partial \Omega ),\mbox{ }\gamma \geq \gamma _{0}>0\mbox{ a.e.}
\label{101}
\end{equation}%
This problem characterizes the evolution of a diffusion process under the
influence of a potential $\beta $ and of a drift term $\nabla \cdot (by).$
For $\beta =0$ the model can describe the diffusion with transport of a
substance in a fluid. If $b=0,$ $\beta (y)=y^{3}$ and $a_{1}=-1$ we note
that this is the Allen-Cahn equation describing the phase transitions of a
material, which can exists in different phases, under the influence of a
double-well potential. Such a problem with different assumptions for $\beta $
was treated in \cite{VB-93}, Section 6.1.4.

We study problem $(\mathcal{P})$ for $u\in L^{\infty }(0,\infty
;L^{2}(\Omega ))$ with $\left\Vert u(t)\right\Vert _{H}\leq \rho $ a.e. $%
t\geq 0.$

\medskip

\noindent \textbf{Proposition 6.1. }\textit{Let}\textbf{\ }$y_{0}\in V,$ $y^{%
\mbox{target}}\in D_{H},$ $\mathit{d=}-\Delta y^{\mbox{target}}+\beta (y^{%
\mbox{target}})+a_{1}y^{\mbox{target}}-\nabla \cdot (by^{\mbox{target}})\in
L^{2}(\Omega )$ \textit{and} $\rho >\left\Vert d\right\Vert _{L^{2}(\Omega
)}.$\textit{\ Then, there exists }$(T^{\ast },u^{\ast })$\textit{\ solution
to }$(\mathcal{P})$ \textit{satisfying }(\ref{g72})-(\ref{g73}), \textit{%
where }$p$\textit{\ solves}%
\begin{eqnarray*}
p_{t}-\Delta p+\beta ^{\prime }(y)p+a_{1}p+b\cdot \nabla p &=&0,%
\mbox{ \ 
\textit{in} }(0,\infty )\times \Omega , \\
\frac{\partial p}{\partial \nu }+\gamma p &=&0,\mbox{ \ \textit{on} }%
(0,\infty )\times \partial \Omega , \\
p(T^{\ast }) &\in &H.
\end{eqnarray*}

\noindent \textbf{Proof. }Let us set:%
\begin{equation*}
H=L^{2}(\Omega ),\mbox{ }V=H^{1}(\Omega ),\mbox{ }V^{\ast }=(H^{1}(\Omega
))^{\ast},\mbox{ }D_{H}=\left\{ y\in H^{2}(\Omega );\mbox{ }\frac{\partial y%
}{\partial \nu }+\gamma y=0\mbox{\ on }\partial \Omega \right\} ,
\end{equation*}%
\begin{equation*}
\Gamma :V\rightarrow V^{\ast },\mbox{ }\left\langle \Gamma y,\psi
\right\rangle _{V^{\ast },V}=\int_{\Omega }\nabla y\cdot \nabla \psi
dx+\int_{\partial \Omega }\gamma (x)y\psi d\sigma ,\mbox{ for }\psi \in V,
\end{equation*}%
\begin{equation*}
\Gamma _{H}:D_{H}\subset H\rightarrow H\mbox{, }\Gamma _{H}=-\Delta ,
\end{equation*}%
\begin{equation*}
U=H,\mbox{ }B=I,\mbox{ }A:V\rightarrow V^{\ast },\mbox{ }
\end{equation*}%
\begin{equation*}
\left\langle Ay,\psi \right\rangle _{V^{\ast },V}=\int_{\Omega }(\nabla
y+by)\cdot \nabla \psi dx+\int_{\Omega }(\beta (y)+a_{1}y)\psi
dx+\int_{\partial \Omega }\gamma (x)y\psi d\sigma ,
\end{equation*}%
\begin{equation*}
A_{H}:D_{H}\subset H\rightarrow H,\mbox{ }A_{H}y=-\Delta y+\beta
(y)+a_{1}y-\nabla \cdot (by).
\end{equation*}%
We shall check first the hypotheses $(a_{1})-(c_{1})$ in Section 3. Since $%
\beta $ is maximal monotone we have 
\begin{eqnarray*}
&&\left\langle Ay-A\overline{y},y-\overline{y}\right\rangle _{V^{\ast
},V}\geq \int_{\Omega }(\left\vert \nabla (y-\overline{y})\right\vert
^{2}+(\beta (y)-\beta (\overline{y}))(y-\overline{y}))dx+\int_{\Omega
}a_{1}(x)(y-\overline{y})^{2}dx \\
&&+\int_{\Omega }b(y-\overline{y})\cdot \nabla (y-\overline{y}%
)dx+\int_{\partial \Omega }\gamma (x)(y-\overline{y})^{2}d\sigma \geq
C_{1}\left\Vert y-\overline{y}\right\Vert _{V}^{2}-C_{2}\left\Vert y-%
\overline{y}\right\Vert _{H}^{2},\mbox{ }y,\overline{y}\in V,
\end{eqnarray*}%
implying that $\lambda I+A$ is coercive for $\lambda $ large. Here, we used
the trace theorem, $\left\Vert y\right\Vert _{L^{2}(\partial \Omega )}\leq
C_{tr}\left\Vert y\right\Vert _{V},$ with $C_{tr}$ is a constant.

By (\ref{100}) it follows that $\left\vert \beta (r)\right\vert \leq
C\left\vert r\right\vert ^{\kappa +1}+\left\vert r\right\vert ,$ and so, for 
$\kappa \in \lbrack 0,2]$ we have that 
\begin{equation*}
\left\Vert \beta (y)\right\Vert _{H}\leq C\int_{\Omega }\left( \left\vert
y\right\vert ^{2(\kappa +1)}+\left\vert y\right\vert ^{2}\right) dx\leq
C_{1}\left( \left\Vert y\right\Vert _{V}^{2(\kappa +1)}+\left\Vert
y\right\Vert _{V}^{2}\right) \cdot
\end{equation*}%
Let $y_{n}\rightarrow y$ strongly in $V.$ Since $\left\Vert \beta
(y_{n})\right\Vert _{H}\leq C$ it follows that $\beta (y_{n})\rightarrow
\beta (y)$ weakly in $H$ because $\beta $ is strongly-weakly closed.
Moreover, we have $\beta (y_{n})\rightarrow \beta (y)$ a.e. on $\Omega ,$
and so $\beta (y_{n})\rightarrow \beta (y)$ strongly in $H,$ by Vitali's
theorem. Therefore, it follows that $A$ is continuous from $V$ to $V^{\ast
}. $ Then, 
\begin{eqnarray*}
&&\left\langle Ay,\psi \right\rangle _{V^{\ast },V}\leq C\left\Vert
y\right\Vert _{V}\left\Vert \psi \right\Vert _{V}+\left\langle \beta
(y)+a_{1}y,\psi \right\rangle _{V^{\ast },V}-\left\langle \nabla \cdot
(by),\psi \right\rangle _{V^{\ast },V} \\
&\leq &C\left\Vert y\right\Vert _{V}\left\Vert \psi \right\Vert
_{V}+\left\Vert \beta (y)\right\Vert _{V^{\ast }}\left\Vert \psi \right\Vert
_{V}+\left\Vert by\right\Vert _{2}\left\Vert \nabla \psi \right\Vert
_{2}\leq C\left( \left\Vert y\right\Vert _{V}+\left\Vert y\right\Vert
_{V}^{k+1}+\left\Vert b\right\Vert _{\infty }\left\Vert y\right\Vert
_{H}\right) \left\Vert \psi \right\Vert _{V},
\end{eqnarray*}%
hence $\left\Vert Ay\right\Vert _{V^{\ast }}$ is bounded on bounded subsets.

Relation (\ref{A0H}) is immediately verified, because 
\begin{eqnarray*}
&&(-\Delta y+\beta (y)+a_{1}y-\nabla \cdot (by),-\Delta y)_{H}\geq
\left\Vert \Delta y\right\Vert _{H}^{2}-\left\Vert a_{1}\right\Vert _{\infty
}\left\Vert \nabla y\right\Vert _{H}^{2}-\left\Vert -\Delta y\right\Vert
_{H}\left\Vert \nabla \cdot (by)\right\Vert _{H} \\
&\geq &\frac{1}{2}\left\Vert \Gamma _{H}y\right\Vert _{H}^{2}-\left\Vert
a_{1}\right\Vert _{\infty }\left\Vert \nabla y\right\Vert
_{H}^{2}-\sum_{i=1}^{3}\left\Vert b_{i}y\right\Vert _{V}^{2}\geq \frac{1}{2}%
\left\Vert \Gamma _{H}y\right\Vert _{H}^{2}-C\left\Vert y\right\Vert
_{V}^{2},
\end{eqnarray*}%
since $(\beta (y),-\Delta y)_{H}\geq 0$ by the monotonicity of $\beta .$

The controllability $(c_{1})$ follows by Proposition 7.1.

Next we verify $(a_{3})-(a_{5})$ in Section 5.1. We introduce $A^{\prime
}(y):V\rightarrow V^{\ast },$%
\begin{equation*}
\left\langle A^{\prime }(y)z,\psi \right\rangle _{V^{\ast },V}=\int_{\Omega
}(\nabla z+bz)\cdot \nabla \psi dx+\int_{\Omega }(\beta ^{\prime
}(y)z+a_{1}z)\psi dx+\int_{\partial \Omega }\gamma (x)z\psi d\sigma ,
\end{equation*}%
then $A_{H}^{\prime }(y):D_{H}\rightarrow H,$ $A_{H}^{\prime }(y)z=-\Delta
z+\beta ^{\prime }(y)z+a_{1}z-\nabla \cdot (bz),$ and%
\begin{equation*}
(A_{H}^{\prime }(y))^{\ast }z=-\Delta z+\beta ^{\prime }(y)z+a_{1}z-b\cdot
\nabla z\mbox{ with }\frac{\partial z}{\partial \nu }+\gamma z=0\mbox{ on }%
\partial \Omega ,
\end{equation*}%
and provide first some estimates. Using the H\"{o}lder inequality we have 
\begin{eqnarray*}
&&I_{1}^{2}=\left\Vert \beta ^{\prime }(y)z\right\Vert _{H}^{2}\leq
C\int_{\Omega }(\left\vert y\right\vert ^{2\kappa }+1)\left\vert
z\right\vert ^{2}dx\leq C\left( \int_{\Omega }\left\vert y\right\vert
^{2\kappa q}dx\right) ^{1/q}\left( \int_{\Omega }\left\vert z\right\vert
^{2q^{\prime }}dx\right) ^{1/q^{\prime }}+C_{1}\left\Vert z\right\Vert
_{H}^{2} \\
&\leq &C\left( \left\Vert y\right\Vert _{2\kappa q}^{2\kappa }\left\Vert
z\right\Vert _{2q^{\prime }}^{2}+\left\Vert z\right\Vert _{H}^{2}\right) ,%
\mbox{ for }y,z\in D_{H},
\end{eqnarray*}%
where $1/q+1/q^{\prime }=1.$ Now, we recall the embedding $W^{s,m}(\Omega
)\subset L^{r}(\Omega ),$ where $d>sm,$ $m\leq r\leq \frac{dm}{d-sm}$ (see 
\cite{Adams}, p. 217, Theorem 7.57) and apply it for $m=2,$ $r=2q^{\prime },$
$s=1-\alpha ,$ for $\alpha \in (0,1)$ to get 
\begin{equation*}
H^{1}(\Omega )\subset W^{1-\alpha ,2}(\Omega )=H^{1-\alpha }(\Omega )\subset
L^{2q^{\prime }}(\Omega ),
\end{equation*}%
with $q^{\prime }>1.$ Then, $\left\Vert y\right\Vert _{2\kappa q}\leq
C\left\Vert y\right\Vert _{V}$ if $2\kappa q\leq 6.$ Thus, we obtain 
\begin{equation}
I_{1}=\left\Vert \beta ^{\prime }(y)z\right\Vert _{H}\leq C\left( \left\Vert
y\right\Vert _{V}^{\kappa }\left\Vert z\right\Vert _{H^{1-\alpha }(\Omega
)}+\left\Vert z\right\Vert _{H}\right) \leq C\left\Vert z\right\Vert
_{V}\left( \left\Vert y\right\Vert _{V}^{\kappa }+1\right) .  \label{103}
\end{equation}%
To this end we must have $3>2(1-\alpha )$ which is satisfied for $\alpha \in
\lbrack 0,1]$ and 
\begin{equation*}
2\leq 2q^{\prime }\leq \frac{6}{3-2(1-\alpha )}\mbox{ implying }q^{\prime
}\leq \frac{3}{1+2\alpha }.
\end{equation*}%
In particular, these are true for $\kappa \leq 2$, $q^{\prime }\geq 3.$ Then%
\begin{equation*}
\left\Vert A^{\prime }(y)z\right\Vert _{V^{\ast }}\leq C_{1}\left\Vert
z\right\Vert _{V}+\left\Vert \beta ^{\prime }(y)z\right\Vert _{V^{\ast
}}\leq C_{2}\left\Vert z\right\Vert _{V}(1+C\left\Vert y\right\Vert
_{V}^{\kappa }),\mbox{ for }y,z\in V,
\end{equation*}%
\begin{equation*}
\left\Vert A_{H}^{\prime }(y)z\right\Vert _{H}\leq C_{1}\left\Vert
z\right\Vert _{D_{H}}(1+C\left\Vert y\right\Vert _{V}^{\kappa }),\mbox{ for }%
y,z\in D_{H}.
\end{equation*}

\noindent Moreover, $y\rightarrow A^{\prime }(y)z$ is continuous from $V$ to 
$L(V,V^{\ast })$. Indeed, let $y_{n}\in V,$ $y_{n}\rightarrow y$ strongly in 
$V.$ Then, as before, $\beta ^{\prime }(y_{n})\rightarrow \beta ^{\prime
}(y) $ strongly in $H.$ Therefore,%
\begin{equation*}
\left\Vert A^{\prime }(y_{n})z-A^{\prime }(y)z\right\Vert _{V^{\ast
}}=\int_{\Omega }(\beta ^{\prime }(y_{n})-\beta ^{\prime }(y))z\psi
dx\rightarrow 0.
\end{equation*}%
Similarly, let $y_{n},$ $y\in D_{H},$ $y_{n}\rightarrow y$ strongly in $V$
and $z\in D_{H}\subset C(\overline{\Omega }).$ Then, 
\begin{equation*}
\left\Vert A^{\prime }(y_{n})z-A^{\prime }(y)z\right\Vert _{H}=\left\Vert
(\beta ^{\prime }(y_{n})-\beta ^{\prime }(y))z\right\Vert _{H}\rightarrow 0.
\end{equation*}

\noindent To prove hypothesis $(a_{6}),$ equivalently (\ref{p-regular}), we
calculate%
\begin{eqnarray*}
&&\left\langle -\Delta z+\beta ^{\prime }(y)z+a_{1}z-b\cdot \nabla z,\Gamma
_{\nu }z\right\rangle _{V^{\ast },V} \\
&\geq &\left\Vert \Gamma _{\nu }z\right\Vert _{H}^{2}-\left\Vert \beta
^{\prime }(y)z\right\Vert _{H}\left\Vert \Gamma _{\nu }z\right\Vert
_{H}-\left\Vert a_{1}\right\Vert _{\infty }\left\Vert z\right\Vert
_{H}\left\Vert \Gamma _{\nu }z\right\Vert _{H}-\left\Vert b\cdot \nabla
z\right\Vert _{H}\left\Vert \Gamma _{\nu }z\right\Vert _{H} \\
&\geq &C\left\Vert \Gamma _{\nu }z\right\Vert _{H}^{2}-C_{2}(\left\Vert
y\right\Vert _{V}^{\kappa }+1)\left\Vert z\right\Vert _{V}^{2},\mbox{ for }%
z\in V.
\end{eqnarray*}%
Here, we used the last inequality in (\ref{103}) and the following relations%
\begin{equation}
\left\langle \Gamma z,\Gamma _{\nu }z\right\rangle _{V^{\ast },V}\geq
\left\Vert \Gamma _{\nu }z\right\Vert _{H}^{2},\mbox{ }z\in V,  \label{700}
\end{equation}%
\begin{equation}
\left\Vert \Gamma _{H}z\right\Vert _{H}\geq \left\Vert \Gamma _{\nu
}z\right\Vert _{H},\mbox{ }z\in D_{H},  \label{700-0}
\end{equation}%
\begin{equation}
\left\Vert z\right\Vert _{D_{H}}=\left\Vert \Gamma _{H}z\right\Vert _{H}.
\label{700-1}
\end{equation}%
Since $U=H$ it remains to check the hypotheses (\ref{1001}), 
\begin{equation*}
((A_{H}^{\prime }(y))^{\ast }z,z)_{H}\geq \left\Vert \nabla z\right\Vert
_{H}^{2}-\left\Vert a_{1}\right\Vert _{\infty }\left\Vert z\right\Vert
_{H}^{2}-\left\Vert b\cdot \nabla z\right\Vert _{H}\left\Vert z\right\Vert
_{H}\geq C_{1}\left\Vert z\right\Vert _{V}^{2}-C_{2}\left\Vert z\right\Vert
_{H}^{2},
\end{equation*}%
and (\ref{1004}) which is automatically verified with $C^{\ast },$ for $\rho 
$ large enough. Thus, Corollary 5.6 can be applied.\hfill $\square $

\medskip

We remark, that in virtue of Remark 5.7, Proposition 6.1 applies to equation
(\ref{ex-1}) with an internal controller%
\begin{equation*}
y_{t}-\Delta y+\beta (y)+a_{1}y-\nabla \cdot (by)=1_{\Omega _{0}}u,%
\mbox{ in 
}(0,\infty )\times \Omega ,
\end{equation*}%
where $\Omega _{0}$ is an open subset of $\Omega $ and $1_{\Omega _{0}}$ is
the characteristic function of $\Omega _{0}.$ Indeed, by \cite{Fursikov-Imm}%
, the corresponding linearized system is exactly null controllable.

\paragraph{Example 2. Porous media equation.}

Let us consider the porous media equation 
\begin{eqnarray}
y_{t}-\Delta \beta (y) &=&u,\mbox{ in }(0,\infty )\times \Omega ,  \notag \\
y &=&0\mbox{ \ on }(0,\infty )\times \partial \Omega ,  \label{ex-6} \\
y(0) &=&y_{0},  \notag
\end{eqnarray}%
where%
\begin{eqnarray}
\beta &:&\mathbb{R}\rightarrow \mathbb{R}\mbox{, \ }\beta \in C^{2}(\mathbb{R%
})\mbox{, \ }\beta (0)=0,  \label{hyp-beta} \\
0 &<&a_{0}\leq \beta ^{\prime }(r)\leq c_{1}\left\vert r\right\vert ^{\kappa
}+c_{2},\mbox{ for }r\in \mathbb{R},\mbox{ }c_{1},c_{2}>0,\mbox{ }0\leq
\kappa <1.  \notag
\end{eqnarray}%
The hypothesis for $\kappa $ places the equation in the slow diffusion case.
We study problem $(\mathcal{P})$ for $u\in L^{\infty }(0,\infty
;H^{-1}(\Omega )),$ $\left\Vert u(t)\right\Vert _{H^{-1}(\Omega )}\leq \rho
, $ a.e. $t\geq 0.$

\medskip

\noindent \textbf{Proposition 6.2. }\textit{Let}\textbf{\ }$y^{\mbox{target}%
}\in H_{0}^{1}(\Omega ),$ $\Delta \beta (y^{\mbox{target}})\in H^{-1}(\Omega
),$\textit{\ }$y_{0}\in H_{0}^{1}(\Omega ),$ $\int_{0}^{y_{0}}\beta (s)ds\in
L^{1}(\Omega ).$\textit{\ Then, there exists }$T^{\ast },$\textit{\ }$%
u^{\ast }$ \textit{and} $y^{\ast }$\textit{\ solution to }$(\mathcal{P})$ 
\textit{satisfying }(\ref{g72})-(\ref{g73}), \textit{where }$U,H,V$ \textit{%
are chosen below and }$p\in C_{w}([0,T^{\ast }];H^{-1}(\Omega ))\cap
L^{2}(0,T^{\ast };L^{2}(\Omega ))$\textit{\ is the solution to} 
\begin{eqnarray}
-p_{t}-\Delta (\beta (y^{\ast })p) &=&0,\mbox{ \textit{in} }(0,T^{\ast
})\times \Omega ,  \notag \\
p &=&0\mbox{ \ \textit{on} }(0,T^{\ast })\times \partial \Omega ,
\label{dual-6} \\
p(T^{\ast }) &\in &H^{-1}(\Omega ).  \notag
\end{eqnarray}

\medskip

\noindent \textbf{Proof. }The proof is led in three steps. First, we prove
an intermediate result for $\beta $ having the properties%
\begin{equation}
0<a_{0}\leq \beta ^{\prime }(r)\leq M_{1},\mbox{ }\left\vert \beta ^{\prime
\prime }(r)\right\vert \leq M_{2},\mbox{ for all }r\in \mathbb{R}.
\label{hyp-bound}
\end{equation}%
Then, we consider (\ref{ex-6}) by replacing $\beta $ by the Yosida
approximation $\beta _{\nu }$ which has the properties (\ref{hyp-bound}) and
obtain the minimum time controllability for the approximating solution $%
y_{\nu }.$ Third, we pass to the limit as $\nu \rightarrow 0.$ To this end,
we choose 
\begin{equation*}
D_{H}=H_{0}^{1}(\Omega ),\mbox{ }V=L^{2}(\Omega ),\mbox{ }H=H^{-1}(\Omega
)\equiv (H^{-1}(\Omega ))^{\ast },\mbox{ }V^{\ast }=(L^{2}(\Omega ))^{\ast },
\end{equation*}%
where $(L^{2}(\Omega ))^{\ast }$ is the dual of $L^{2}(\Omega )$ in the
pairing with $H^{-1}(\Omega )$ as pivot space. Moreover, 
\begin{equation*}
P=I,\mbox{ }B=I,\mbox{ }U=H^{-1}(\Omega )=U^{\ast }\mbox{ and }\Gamma
_{H}:D_{H}\subset H\rightarrow H,\mbox{ }\Gamma _{H}=-\Delta .
\end{equation*}%
We define the operator $A:V\rightarrow V^{\ast }$ by 
\begin{equation*}
\left\langle Ay,\psi \right\rangle _{V^{\ast },V}=(\beta (y),\psi )_{V},%
\mbox{ for }y,\psi \in V=L^{2}(\Omega ),
\end{equation*}%
and $A_{H}:D_{H}\subset H\rightarrow H$ by $A_{H}y=-\Delta \beta (y).$

The norm on $V^{\ast }=(L^{2}(\Omega ))^{\ast }$ is given by $(\theta
,\theta )_{V^{\ast }}=\left\Vert \psi \right\Vert _{L^{2}(\Omega )},$ where $%
\theta =A\psi .$

The controllability $(c_{1})$ follows by Proposition 7.1. We begin to check
the hypotheses of Corollary 5.6. First, 
\begin{equation*}
\left\langle Ay-A\overline{y},y-\overline{y}\right\rangle _{V^{\ast
},V}=(\beta (y)-\beta (\overline{y}),y-\overline{y})_{V}\geq a_{0}\left\Vert
y-\overline{y}\right\Vert _{V}^{2},
\end{equation*}%
which implies the coercivity, too. Then,%
\begin{equation*}
(A_{H}y,\Gamma _{H}y)_{H}=\left\langle -\Delta \beta (y),y\right\rangle
_{H,D_{H}}=\int_{\Omega }\beta ^{\prime }(y)\left\vert \nabla y\right\vert
^{2}dx\geq a_{0}\left\Vert y\right\Vert _{D_{H}}^{2}.
\end{equation*}

\noindent We have $A_{H}^{\prime }(y)z=-\Delta (\beta ^{\prime }(y)z)$ and $%
(A^{\prime }(y))^{\ast }z=-\Delta (\beta ^{\prime }(y)z)$, where $\beta
^{\prime }(y)p\in V=L^{2}(\Omega ).$ Next,%
\begin{equation*}
\left\Vert A^{\prime }(y)z\right\Vert _{V^{\ast }}=\left\Vert \beta ^{\prime
}(y)z\right\Vert _{V}\leq M_{1}\left\Vert z\right\Vert _{V},\mbox{ for }%
y,z\in V=L^{2}(\Omega ),
\end{equation*}%
and if $y_{n}\rightarrow y$ strongly in $V=L^{2}(\Omega ),$ we have%
\begin{equation*}
\left\Vert A^{\prime }(y_{n})z-A^{\prime }(y)z\right\Vert _{V^{\ast
}}=\left\Vert (\beta ^{\prime }(y_{n})-\beta ^{\prime }(y))z\right\Vert
_{L^{2}(\Omega )}\rightarrow 0\mbox{ as }n\rightarrow \infty .
\end{equation*}%
This follows by Lebesgue dominated convergence theorem since $\beta ^{\prime
}(y_{n})z\rightarrow \beta ^{\prime }(y)z$ a.e. on $\Omega $ and $\left\vert
(\beta ^{\prime }(y_{n})-\beta ^{\prime }(y))z\right\vert _{L^{2}(\Omega
)}\leq 2M_{1}\left\vert z\right\vert .$ Then, since we can write $\beta
(r)=a_{0}r+\beta _{1}^{\prime }(r)$ we have 
\begin{equation*}
\left\langle A^{\prime }(y)z,\Gamma _{\nu }(z)\right\rangle _{V^{\ast
},V}\geq \left\Vert \Gamma _{\nu }z\right\Vert _{H}^{2}-(\beta _{1}^{\prime
}(y)z,\Gamma _{\nu }z)_{L^{2}(\Omega )}\geq C_{1}\left\Vert \Gamma _{\nu
}z\right\Vert _{H}^{2}-C_{2}\left\Vert z\right\Vert _{L^{2}(\Omega )}^{2}.
\end{equation*}%
Finally, we have to check (\ref{1001}), that is 
\begin{equation*}
(-\Delta (\beta ^{\prime }(y)z),z)_{H}=\left\langle -\Delta (\beta ^{\prime
}(y)z),z\right\rangle _{V^{\ast },V}=\int_{\Omega }\beta ^{\prime
}(y)z^{2}dx\geq a_{0}\left\Vert z\right\Vert _{V}^{2},
\end{equation*}%
while (\ref{1004}) which is automatically verified. Thus, we get a minimum
time and a controller satisfying the thesis of Corollary 5.6.

In the second step we replace $\beta $ by $\beta _{\nu }$ in (\ref{ex-6}).
Both $\beta _{\nu }^{\prime }$ and $\beta _{\nu }^{\prime \prime }$ are
bounded by constants $C_{\nu },$ for each $\nu >0.$ On the basis of the
previous result we obtain that there exists $T_{\nu }^{\ast },$ $u_{\nu
}^{\ast }$ and $y_{\nu }^{\ast }$ satisfying 
\begin{equation}
Pu_{\nu }^{\ast }(t)\in (N_{K})^{-1}\left( -B^{\ast }p_{\nu }(t)\right) ,%
\mbox{ a.e. }t\in (0,T_{\nu }^{\ast }),  \label{opt1}
\end{equation}%
\begin{equation}
\rho \left\Vert p_{\nu }(t)\right\Vert _{U^{\ast }}+\int_{\Omega }\beta
_{\nu }(y_{\nu }^{\ast }(t)p_{\nu }(t)dx=1,\mbox{ a.e. }t\in (0,T_{\nu
}^{\ast }),\mbox{ }  \label{opt2}
\end{equation}%
where $y_{\nu }^{\ast }$ is the solution to\ the approximating state system (%
\ref{ex-6}) (with $\beta _{\nu }),$ corresponding to $(T_{\nu }^{\ast
},u_{\nu }^{\ast }),$\ and $p_{\nu }$ is a solution to\textit{\ } 
\begin{equation}
-p_{\nu }^{\prime }(t)-\Delta (\beta _{\nu }^{\prime }(y_{\nu }^{\ast
})p_{\nu })=0,\mbox{ a.e.\textit{\ }}t\in (0,T_{\nu }^{\ast }),  \label{opt3}
\end{equation}%
\begin{equation}
\left\Vert p_{\nu }(T_{\nu }^{\ast })\right\Vert _{H^{-1}(\Omega )}\leq C.
\label{opt4}
\end{equation}%
A first estimate for $y_{\nu }^{\ast }$ reads%
\begin{equation}
\left\Vert y_{\nu }^{\ast }\right\Vert _{L^{\infty }(0,T_{\nu }^{\ast
};L^{2}(\Omega ))\cap L^{2}(0,T_{\nu }^{\ast };H_{0}^{1}(\Omega ))}\leq C,
\label{est1}
\end{equation}%
where $C$ denote several constants. By multiplying the approximating
equation (\ref{ex-6}) by $\beta _{\nu }(y_{\nu }^{\ast }(t))$ and
integrating on $(0,t)$ we obtain%
\begin{equation*}
\int_{\Omega }j_{\nu }(y_{\nu }^{\ast }(t))dx+\int_{0}^{t}\left\Vert \nabla
\beta _{\nu }(y_{\nu }^{\ast }(\tau ))\right\Vert _{L^{2}(\Omega )}^{2}d\tau
\leq \int_{\Omega }j_{\nu }(y_{0})dx+\int_{0}^{t}\left\Vert u_{\nu }(\tau
)\right\Vert _{H^{-1}(\Omega )}\left\Vert \beta _{\nu }(y_{\nu }^{\ast
}(\tau ))\right\Vert _{H_{0}^{1}(\Omega )}d\tau ,
\end{equation*}%
where $\partial j_{\nu }(r)=\beta _{\nu }(r)$ and $\partial j(r)=\beta (r)$
for all $r\in \mathbb{R}.$ This implies 
\begin{equation}
\int_{\Omega }j_{\nu }(y_{\nu }^{\ast }(t))dx+\int_{0}^{t}\left\Vert \nabla
\beta _{\nu }(y_{\nu }^{\ast }(\tau ))\right\Vert _{L^{2}(\Omega )}^{2}d\tau
\leq C\left( \int_{\Omega }j(y_{0})dx+T_{\nu }^{\ast }\rho ^{2}\right) .
\label{est2}
\end{equation}%
Since $j(r)=\int_{0}^{r}\beta (s)ds$ and $j(y_{0})\in L^{1}(\Omega )$ it
follows that the right-hand side in (\ref{est2}) is bounded independently of 
$\nu .$ This yields%
\begin{equation}
\left\Vert \beta _{\nu }(y_{\nu }^{\ast })\right\Vert _{L^{2}(0,T_{\nu
}^{\ast };H_{0}^{1}((\Omega ))}\leq C.  \label{est2-0}
\end{equation}%
Then, we multiply (\ref{opt3}) by $p_{\nu }(t)$ and integrate over $%
(t,T_{\nu }^{\ast }),$ getting 
\begin{equation}
\left\Vert p_{\nu }\right\Vert _{L^{\infty }(0,T_{\nu }^{\ast
};H^{-1}(\Omega ))\cap L^{2}(0,T_{\nu }^{\ast };L^{2}(\Omega ))}\leq C.
\label{est3}
\end{equation}%
Next, we determine an estimate for $A^{\prime }(y_{\nu }^{\ast })p_{\nu }$
and begin by computing 
\begin{eqnarray*}
\int_{\Omega }\left\vert \beta _{\nu }^{\prime }(y_{\nu }^{\ast }(t)p_{\nu
}(t))\right\vert ^{q}dx &\leq &C_{1}\left( \int_{\Omega }\left\vert p_{\nu
}(t)\right\vert ^{2}dx\right) ^{q/2}\left( \left( \int_{\Omega }\left\vert
y_{\nu }^{\ast }(t)\right\vert ^{\kappa q^{\prime }}dx\right) ^{1/q^{\prime
}}+1\right)  \\
&\leq &C_{1}\left\Vert p_{\nu }(t)\right\Vert _{L^{2}(\Omega
)}^{q}(\left\Vert y_{\nu }(t)\right\Vert _{\kappa q^{\prime }}^{\kappa }+1),%
\mbox{ for a.e. }t,
\end{eqnarray*}%
where $\frac{1}{q^{\prime }}=1-\frac{q}{2},$ that is $q^{\prime }=\frac{2}{%
2-q},$ for $1<q<2$ and $\kappa q^{\prime }=\frac{2\kappa q}{2-q}\leq 2,$
meaning that $q\leq \frac{2}{\kappa +1}$, which is true if $\kappa <1.$
Therefore, by (\ref{est1}) and (\ref{est3}) we obtain%
\begin{equation}
\left\Vert \beta _{\nu }^{\prime }(y_{\nu }^{\ast })p_{\nu }\right\Vert
_{L^{q}(0,T_{\nu }^{\ast };L^{q}(\Omega ))}\leq C_{2},\mbox{ }1<q<2.
\label{est4}
\end{equation}%
This implies that 
\begin{equation}
\left\Vert \Delta \beta _{\nu }^{\prime }(y_{\nu }^{\ast })p_{\nu
}\right\Vert _{L^{q}(0,T_{\nu }^{\ast };X)}+\left\Vert p_{\nu }^{\prime
}\right\Vert _{L^{q}(0,T_{\nu }^{\ast };X)}\leq C_{3}  \label{est5}
\end{equation}%
where $X$ is the image of $L^{q}(\Omega )$ by the operator $-\Delta .$ More
precisely, $X$ is the completion of $L^{q}(\Omega )$ in the norm $\left\vert
\left\Vert w\right\Vert \right\vert _{X}=\left\Vert A^{-1}w\right\Vert
_{L^{q}(\Omega )}.$ Moreover, applying the same argument as in Theorem 5.5
we can deduce that $T_{\nu }^{\ast }\rightarrow T^{\ast }$, and on a
subsequence, it follows that%
\begin{equation}
y_{\nu }^{\ast }\rightarrow y^{\ast }\mbox{ weakly in }W^{1,2}(0,T^{\ast
};H^{-1}(\Omega ))\cap L^{2}(0,T^{\ast };H_{0}^{1}(\Omega )),%
\mbox{ strongly
in }L^{2}(0,T^{\ast };L^{2}(\Omega )),  \notag
\end{equation}%
\begin{equation*}
\beta _{\nu }(y_{\nu }^{\ast })\rightarrow \beta (y^{\ast })%
\mbox{ strongly
in }L^{2}(0,T^{\ast };L^{2}(\Omega )),
\end{equation*}%
since $y_{\nu }\rightarrow y$ strongly, $\beta _{\nu }(y_{\nu }^{\ast
})\rightarrow \eta $ weakly in $L^{2}(0,T_{\nu }^{\ast };L^{2}(\Omega ))$
and $\beta $ is strongly-weakly closed. Then, 
\begin{eqnarray}
p_{\nu } &\rightarrow &p\mbox{ weakly in }W^{1,2}(0,T^{\ast };X)\cap
L^{2}(0,T^{\ast };L^{2}(\Omega )),  \notag \\
&&\mbox{weak-star in }L^{\infty }(0,T^{\ast };H^{-1}(\Omega )),%
\mbox{
strongly in }L^{2}(0,T^{\ast };H^{-1}(\Omega )),  \label{est7}
\end{eqnarray}%
and%
\begin{equation}
\beta _{\nu }^{\prime }(y_{\nu }^{\ast })p_{\nu }\rightarrow \zeta =\beta
^{\prime }(y^{\ast })p\mbox{ weakly in }L^{q}(0,T^{\ast };L^{q}(\Omega )).
\label{est8}
\end{equation}%
Indeed, $\beta _{\nu }^{\prime }(y_{\nu }^{\ast })\rightarrow \beta ^{\prime
}(y^{\ast })$ a.e., $\int_{0}^{T_{\nu }^{\ast }}\left\Vert \beta _{\nu
}^{\prime }(y_{\nu }^{\ast }(t))\right\Vert _{L^{2}(\Omega )}dt\leq
\int_{0}^{T_{\nu }^{\ast }}\left\Vert y_{\nu }^{\ast }(t)\right\Vert
_{H_{0}^{1}(\Omega )}^{2}dt\leq C$ and so 
\begin{equation*}
\beta _{\nu }^{\prime }(y_{\nu }^{\ast })\rightarrow \beta ^{\prime
}(y^{\ast })\mbox{ weakly in }L^{2}(0,T^{\ast };L^{2}(\Omega )).
\end{equation*}%
Then, by (\ref{est7}) 
\begin{equation*}
\beta _{\nu }^{\prime }(y_{\nu }^{\ast })p_{\nu }\rightarrow \beta ^{\prime
}(y^{\ast })p\mbox{ weakly in }L^{1}(0,T^{\ast };L^{1}(\Omega ))
\end{equation*}%
and choosing $\varphi \in C_{0}^{\infty }((0,T)\times \Omega )$ with $%
T>T^{\ast }$ we have%
\begin{equation*}
\int_{0}^{T_{\nu }^{\ast }}\int_{\Omega }(\beta _{\nu }^{\prime }(y_{\nu
}^{\ast })p_{\nu }-\beta ^{\prime }(y^{\ast })p)\varphi dxdt\rightarrow 0,
\end{equation*}%
which yields that $\zeta =\beta ^{\prime }(y^{\ast })p$ a.e. Thus, (\ref%
{est8}) holds true. Now, we can pass to the limit in (\ref{opt3}) and (\ref%
{opt4}) to get (\ref{dual-6}) and in (\ref{opt2}) to deduce%
\begin{equation*}
\rho \left\Vert B^{\ast }p(t)\right\Vert _{H^{-1}(\Omega )}+(Ay^{\ast
}(t),p(t))_{H^{-1}(\Omega )}=1,\mbox{ a.e. }t\in (0,T^{\ast }).
\end{equation*}%
Finally, we pass to the limit in (\ref{opt1}), written as $-B^{\ast }p_{\nu
}\in N_{\mathcal{K}_{T^{\ast }}}(Pu_{\nu }^{\ast }),$ taking into account
that $Pu_{\nu }^{\ast }\rightarrow Pu^{\ast }$ weakly in $L^{2}(0,T^{\ast
};H^{-1}(\Omega )),$ $-B^{\ast }p_{\nu }\rightarrow -B^{\ast }p$ strongly in 
$L^{2}(0,T^{\ast };H^{-1}(\Omega )),$ and $N_{\mathcal{K}_{T^{\ast }}}$ is
weakly-strongly closed. Here, $\mathcal{K}_{T^{\ast }}=\{w\in
L^{2}(0,T^{\ast };H^{-1}(\Omega ));$ $\left\Vert w(t)\right\Vert
_{H^{-1}(\Omega )}\leq \rho $ a.e. $t\}.\hfill \square $

\paragraph{Example 3. Sliding mode control for reaction-diffusion systems
with nonlinear perturbations.}

Let us consider the system%
\begin{eqnarray}
y_{t}-D_{1}\Delta y+f(y,z) &=&u\mbox{, \ \ \ \ \ \ \ \ in }(0,\infty )\times
\Omega ,  \label{ex-2} \\
z_{t}-D_{2}\Delta z+g(y,z) &=&0,\mbox{ \ \ \ \ \ \ \ \ in }(0,\infty )\times
\Omega ,  \notag \\
\frac{\partial y}{\partial \nu } &=&\frac{\partial z}{\partial \nu }=0%
\mbox{, \ \ \ \ \ \ \ \ \ \ on }(0,\infty )\times \partial \Omega ,  \notag
\\
y(0) &=&y_{0},\mbox{ }z(0)=z_{0},\mbox{ \ \ in }\Omega .  \notag
\end{eqnarray}%
For certain expressions of $f$ and $g$, equation (\ref{ex-2}) can model
different reaction-diffusion processes, as for instance the diffusion, in a
habitat $\Omega ,$ of two populations with the densities $y$ and $z,$
interacting between them according to the laws expressed by $f$ and $g.$

In some situations, (\ref{g5-000}) can be satisfied and so one can control
the first component of the state $y,$ with one controller, letting $z$
uncontrolled. In this example we shall focus on the situation when $V\subset
U$ and prove the minimum time sliding mode control for this system.

\medskip

\noindent \textit{Case I.} Let us consider that $f,$ $g$ are generally
nonlinear, $f,$ $g\in C^{2}(\mathbb{R\times R)},$ such that 
\begin{equation}
\sup_{(r_{1},r_{2})\in \mathbb{R\times R}}(\left\vert \nabla
f(r_{1},r_{2})\right\vert +\left\vert \nabla g(r_{1},r_{2})\right\vert )\leq
M  \label{rd1}
\end{equation}%
and $D_{i}>0,$ $i=1,2.$ We study problem $(\mathcal{P})$ with $%
U=L^{4}(\Omega ).$ We set 
\begin{eqnarray*}
H &=&L^{2}(\Omega )\times L^{2}(\Omega ),\mbox{ }V=H^{1}(\Omega )\times
H^{1}(\Omega ),\mbox{ }V^{\ast }=(H^{1}(\Omega ))^{\ast }\times
(H^{1}(\Omega ))^{\ast },\mbox{ } \\
D &=&\left\{ w\in H^{2}(\Omega );\mbox{ }\frac{\partial w}{\partial \nu }=0%
\mbox{ on }\partial \Omega \right\} ,\mbox{ }D_{H}=D\times D, \\
U &=&(L^{4}(\Omega ),L^{4}(\Omega )),\mbox{ }U^{\ast }=(L^{4/3}(\Omega
),L^{4/3}(\Omega ))),\mbox{ }(B(u,u_{2})=(u,0),\mbox{ } \\
P(y,z) &=&(y,0),\mbox{ }y^{tar}=(y_{1}^{\mbox{target}},z_{1}),\mbox{ }%
\forall z_{1}\in L^{2}(\Omega ),
\end{eqnarray*}%
\begin{eqnarray*}
A &:&V\rightarrow V^{\ast },\mbox{ }\left\langle A(y,z),(\psi _{1},\psi
_{2})\right\rangle _{V^{\ast },V}=D_{1}\int_{\Omega }\nabla y\cdot \nabla
\psi _{1}dx+D_{2}\int_{\Omega }\nabla z\cdot \nabla \psi _{2}dx \\
&&+\int_{\Omega }(f(y,z)\psi _{1}+g(y,z)\psi _{2})dx,\mbox{ for }(y,z)\in V,
\end{eqnarray*}%
and $A_{H}:D_{H}\subset H\rightarrow H,$ 
\begin{equation*}
A_{H}w=\left[ 
\begin{tabular}{l}
$-D_{1}\Delta y+f(y,z)$ \\ 
$-D_{2}\Delta z+g(y,z)$%
\end{tabular}%
\right] ,\mbox{ }w=(y,z).
\end{equation*}%
In this case,%
\begin{equation*}
\mathbf{\Gamma }_{H}:D_{H}\subset H\rightarrow H\mbox{, }\mathbf{\Gamma }%
_{H}(y,z)=((I-\Delta )y,(I-\Delta )z)=(\Gamma _{H}y,\Gamma _{H}z)\mbox{,}
\end{equation*}%
with the homogeneous Neumann boundary condition. The controllability $%
(c_{1}) $ can follow as in \cite{BCGMR-sli}. Namely, first it is proved that
there exists a controller $u\in L^{\infty }((0,\infty )\times \Omega )$ and $%
T_{\ast },$ such that for $\rho $ large enough $y(T_{\ast })=y_{1}^{%
\mbox{target}},$ where $T_{\ast }$ depends on $\rho .$ This controller
belongs also to $L^{\infty }(0,\infty ;L^{4}(\Omega ))$ but the
controllability follows with a different $\rho $ calculated from a relation
between the norms in $L^{4}(\Omega )$ and $L^{\infty }(\Omega ).$ Moreover,
the time $T_{\ast }$ is smaller as $\rho $ is greater.

Let $K=\{u\in U=L^{4}(\Omega )\times L^{4}(\Omega );$ $\left\Vert
u\right\Vert _{U}\leq \rho \}.$

\medskip

\noindent \textbf{Proposition 6.3. }\textit{Assume }(\ref{rd1}) \textit{and} 
$y_{0}\in H^{1}(\Omega ),$ $y_{1}^{\mbox{target}}\in D,$ $\Delta y_{1}^{%
\mbox{target}}\in H^{1}(\Omega ).$\textit{\ Then, there exists }$T^{\ast },$%
\textit{\ }$u^{\ast }$\textit{\ solution to }$(\mathcal{P})$ \textit{%
satisfying }(\ref{g72})-(\ref{g73}), 
\begin{equation*}
Pu^{\ast }(t)\in (N_{K})^{-1}\left( -p(t)\right) ,\mbox{ \textit{a.e.} }t\in
(0,T^{\ast }),
\end{equation*}%
\begin{equation*}
\rho \left\Vert p(t)\right\Vert _{U^{\ast }}+(A_{H}y^{\ast }(t),v(t))_{H}=1,%
\mbox{ \textit{a.e.} }t\in (0,T^{\ast }),\mbox{ }
\end{equation*}%
\textit{where }$v=(p,q)$\textit{\ is the solution to} 
\begin{eqnarray*}
-p_{t}-D_{1}\Delta p+f_{y}(y^{\ast },z^{\ast })p+g_{y}(y^{\ast },z^{\ast })q
&=&0\mbox{, \ \ \ \ \ \ \ \ \ \ \ \ \textit{in} }(0,\infty )\times \Omega ,
\\
-q_{t}-D_{2}\Delta q+f_{z}(y^{\ast },z^{\ast })p+g_{z}(y^{\ast },z^{\ast })q
&=&0,\mbox{, \ \ \ \ \ \ \ \ \ \ \ \ \textit{in} }(0,\infty )\times \Omega ,
\\
\frac{\partial p}{\partial \nu } &=&\frac{\partial q}{\partial \nu }=0%
\mbox{, \ \ \ \ \ \ \ \ \textit{on} }(0,\infty )\times \partial \Omega , \\
p(T^{\ast }) &\in &V^{\ast },\mbox{ }q(T^{\ast })=0,\mbox{ \textit{in} }%
\Omega .
\end{eqnarray*}%
\textit{Moreover, }$y(t)=y_{1}^{\mbox{target}}$\textit{\ for }$t>T^{\ast }.$

\medskip

\noindent \textbf{Proof. }It is obvious that $A$ is continuous, monotone and
coercive\ and (\ref{g5})-(\ref{D0}), $(a_{2})$ are satisfied. We have for $%
w:=(y,z)$ and $v:=(p,q)$ 
\begin{equation*}
A_{H}^{\prime }(y,z)=\left[ 
\begin{tabular}{ll}
$-D_{1}\Delta +f_{y}(y,z)$ & $f_{z}(y,z)$ \\ 
$g_{y}(y,z)$ & $-D_{2}\Delta +g_{z}(y,z)$%
\end{tabular}%
\right] ,
\end{equation*}%
\begin{equation*}
(A_{H}^{\prime }(y,z))^{\ast }=\left[ 
\begin{tabular}{ll}
$-D_{1}\Delta +f_{y}(y,z)$ & $g_{y}(y,z)$ \\ 
$f_{z}(y,z)$ & $-D_{2}\Delta +g_{z}(y,z)$%
\end{tabular}%
\right] ,
\end{equation*}%
where $f_{y},$ $f_{z},$ $g_{y},$ $g_{z}$ denote the partial derivatives of $%
f $ and $g$ with respect to their arguments and they belong to $L^{\infty }(%
\mathbb{R}\times \mathbb{R}).$ Then, 
\begin{equation*}
\left( A_{H}(y,z),\mathbf{\Gamma }_{H}(y,z)\right) _{H}\geq C(\left\Vert
\Gamma _{H}y\right\Vert _{2}^{2}+\left\Vert \Gamma _{H}z\right\Vert
_{2}^{2}-C_{1}\left\Vert (y,z)\right\Vert _{V}^{2}
\end{equation*}%
and it is easy to see that $(a_{3})-(a_{5})$ are satisfied. Regarding $%
(a_{6})$ or (\ref{p-regular}) we have for $v=(p,q),$%
\begin{equation*}
\left\langle (A_{H}^{\prime }(y,z))^{\ast }v,\mathbf{\Gamma }_{\nu
}v\right\rangle _{V^{\ast },V}\geq \left\Vert \mathbf{\Gamma }_{\nu
}v\right\Vert _{H}^{2}-\left\Vert v\right\Vert _{H}\left\Vert \mathbf{\Gamma 
}_{\nu }v\right\Vert _{H}\geq C\left\Vert \mathbf{\Gamma }_{\nu
}v\right\Vert _{H}^{2}-C_{1}\left\Vert v\right\Vert _{V}^{2},
\end{equation*}%
where $\mathbf{\Gamma }_{\nu }(p,q)=(\Gamma _{\nu }p,\Gamma _{\nu }q).$

Next, we verify the hypotheses in Section 5.3. Relation (\ref{g74-2}) is
satisfied because 
\begin{equation*}
\left\Vert Pv\right\Vert _{V^{\ast }}=\left\Vert p\right\Vert
_{(H^{1}(\Omega ))^{\ast }}\leq C^{\ast }\left\Vert p\right\Vert
_{4/3}=C^{\ast }\left\Vert B^{\ast }p\right\Vert _{U^{\ast }},
\end{equation*}%
since $H^{1}(\Omega )\subset L^{4}(\Omega )$ for $d\leq 3.$ Also, 
\begin{equation*}
\left\Vert P\mathbf{\Gamma }_{H}^{-\alpha /2}v\right\Vert _{H}=\left\Vert
\Gamma _{H}^{-\alpha /2}p\right\Vert _{H}\leq C\left\Vert p\right\Vert
_{H^{-\alpha }(\Omega )}\leq C_{1}\left\Vert p\right\Vert
_{4/3}=C_{1}\left\Vert B^{\ast }p\right\Vert _{U^{\ast }},\mbox{ for }v\in H,
\end{equation*}%
\begin{eqnarray*}
\left( A_{H}^{\prime }(y)v,\mathbf{\Gamma }_{H}^{-1}v\right) _{H} &\geq
&C\left\Vert v\right\Vert _{H}^{2}-\left\Vert f_{y}(y,z)p\right\Vert
_{H}\left\Vert \Gamma _{H}^{-1}p\right\Vert _{H}-\left\Vert
g_{z}(y,z)q\right\Vert _{H}\left\Vert \Gamma _{H}^{-1}q\right\Vert _{H} \\
&\geq &C\left\Vert v\right\Vert _{H}^{2}-M\left\Vert v\right\Vert
_{H}\left\Vert v\right\Vert _{V^{\ast }}\geq C_{1}\left\Vert v\right\Vert
_{H}^{2}-C_{2}\left\Vert v\right\Vert _{V^{\ast }}^{2}
\end{eqnarray*}%
and 
\begin{equation*}
\left( A_{H}^{\prime }(y)v,\mathbf{\Gamma }_{H}^{-\alpha }v\right) _{H}\geq
C\left\Vert \mathbf{\Gamma }_{H}^{(1-\alpha )/2}v\right\Vert
_{H}^{2}-C_{1}\left\Vert v\right\Vert _{H}^{2},
\end{equation*}%
where $\mathbf{\Gamma }_{H}^{-\alpha }(y,z)=(\Gamma _{H}^{-\alpha }y,\Gamma
_{H}^{-\alpha }z).$ Here we used the estimate 
\begin{equation*}
\left\Vert \Gamma _{H}^{-\alpha /2}p\right\Vert _{H}=\left\Vert p\right\Vert
_{H^{-\alpha }(\Omega )}\leq \left\Vert v\right\Vert _{H}
\end{equation*}%
since $L^{2}(\Omega )\subset H^{-\alpha }(\Omega ).$ Moreover, $D(\Gamma
_{H}^{-\alpha /2})=H^{-\alpha }(\Omega ),$ according to the characterization
of the domains of the fractionary powers of $-\Delta $ given in \cite%
{Fujiwara}, Theorem 2, for $\alpha <1.$

Now, we have to check (\ref{Ay^}) and (\ref{g5-000}). To this end, we assume
that $\Delta y_{1}^{\mbox{target}}\in H^{1}(\Omega )$ and choose $\widehat{z}
$ to be exactly the second component $z$ of the solution to (\ref{ex-2}).
First, we prove that it has the necessary regularity. We recall (\ref{g15})
and (\ref{g15-10}), that is 
\begin{equation*}
\left\Vert (y,z)(t)\right\Vert _{V}^{2}+\left\Vert (y,z)\right\Vert
_{L^{2}(0,T;H)\cap L^{2}(0,T;D_{H})}\leq C\left( \left\Vert
(y_{0},z_{0})\right\Vert _{V}^{2}+\int_{0}^{t}\left\Vert u(\tau )\right\Vert
_{U}^{2}d\tau \right) e^{Ct},
\end{equation*}%
for all $t\in \lbrack 0,T],$ $T>0.$ Since $y_{t}\in L^{2}(0,T;H)$ we expect
to have a more regular $z.$ We consider the equation 
\begin{equation}
w_{t}-D_{2}\Delta w+g_{z}(y,z)w=-g_{y}(y,z)y_{t},  \label{200}
\end{equation}%
with $\frac{\partial w}{\partial \nu }=0$ and $w(0)=z_{t}(0)=D_{2}\Delta
z_{0}-g(y_{0},z_{0})\in P(V)=H^{1}(\Omega ).$ This is computed directly from
the equation for $z$, observing that by the hypotheses for $g$ and the
initial data we have $g(y_{0}z_{0})\in (I-P)(V)=H^{1}(\Omega ).$ We note
that (\ref{200}) represents also the equation for $z_{t},$ obtained by
formally differentiating the equation in $z.$ Since in (\ref{200}) all
coefficients on the left-hand side are in $L^{\infty }(\Omega )$ and $%
g_{y}(y,z)y_{t}\in L^{2}(0,T;H),$ it follows that it has an unique solution 
\begin{equation*}
w\in C_{w}([0,T];V)\cap L^{\infty }(0,T;V)\cap L^{2}(0,T;W)\cap
W^{1,2}(0,T;H)
\end{equation*}%
and so we can deduce that $w=z_{t}$ and that $z_{t}$ belongs to the same
spaces. Moreover, by multiplying (\ref{200}) by $w_{t}$ we obtain, by some
calculations similar to those in Theorem 3.2, that 
\begin{equation*}
\left\Vert z_{t}(t)\right\Vert _{C_{w}([0,T];V)\cap L^{2}(0,T;W)\cap
W^{1,2}(0,T;H)}\leq C_{T}^{1},\mbox{ for all }t\in \lbrack 0,T],
\end{equation*}%
where this constant depends also on $z_{t}(0)\in P(V),$ namely on $%
\left\Vert \Delta z_{0}\right\Vert _{V},$ 
\begin{equation*}
C_{T}^{1}:=C\left( \left\Vert (y_{0}\right\Vert _{H^{1}(\Omega
)}^{2}+\left\Vert \Delta z_{0}\right\Vert _{H^{1}(\Omega )}+T\rho
^{2}\right) e^{CT}.
\end{equation*}%
Going back to the equation in $z$ we have 
\begin{equation*}
-D_{2}\Delta z(t)+g(y(t),z(t))=-z_{t}(t)\in V\mbox{,}
\end{equation*}%
because $\left\Vert g(y(t),z(t))\right\Vert _{H^{1}(\Omega )}\leq C_{T}.$
Indeed, e.g., $\left\Vert g_{y}(y(t),z(t))\nabla y(t)\right\Vert
_{H^{1}(\Omega )}\leq MC_{T},$ by (\ref{g15}). Next, in the same way we see
that if $\Delta y_{1}^{\mbox{target}}\in H^{1}(\Omega ),$ then $f(y_{1}^{%
\mbox{target}},z(t))\in H^{1}(\Omega )$ and so $-D_{1}\Delta y_{1}^{%
\mbox{target}}+f(y_{1}^{\mbox{target}},z(t))\in H^{1}(\Omega )$. Moreover, 
\begin{eqnarray*}
&&\left\Vert -D_{1}\Delta y_{1}^{\mbox{target}}+f(y(t),z(t))\right\Vert
_{H^{1}(\Omega )}+\left\Vert -D_{2}\Delta z(t)+g(y(t),z(t))z(t)\right\Vert
_{H^{1}(\Omega )} \\
&\leq &C_{T}+C_{T}^{1}\leq C(C_{0}+\sqrt{T}\rho )e^{CT},\mbox{ }t\geq 0.
\end{eqnarray*}%
Here, $C_{0}=\left\Vert y_{0}\right\Vert _{H^{1}(\Omega )}+\left\Vert \Delta
z_{0}\right\Vert _{V}.$ This implies that 
\begin{equation*}
\left\Vert A(y_{1}^{\mbox{target}},\widehat{z})\right\Vert _{V}\leq C(C_{0}+%
\sqrt{T}\rho )e^{CT}.
\end{equation*}%
In order to satisfy (\ref{Ay^}) we have to impose that%
\begin{equation}
\rho >C^{\ast }C(C_{0}+\sqrt{T}\rho )e^{CT}.  \label{203}
\end{equation}%
We can check that if 
\begin{equation}
\sqrt{T}e^{CT}<\frac{1}{C}  \label{201}
\end{equation}%
then%
\begin{equation}
\rho >\frac{CC_{0}C^{\ast }e^{CT}}{1-C\sqrt{T}e^{CT}}  \label{202}
\end{equation}%
and consequently (\ref{203}) are satisfied. We note that, for any positive
constant $C,$ the equation $\sqrt{T}e^{CT}=\frac{1}{C}$ has a unique
solution $T_{\ast \ast }$ and so any $T\in \lbrack 0,T_{\ast \ast })$
verifies (\ref{201}). We can choose $\rho $ sufficiently large, such that
the time $T_{\ast }$ in hypothesis $(c_{1})$ becomes smaller enough, such
that to remain in $(0,T_{\ast \ast }).$ We have to check (\ref{g5-000}),
that is%
\begin{eqnarray*}
&&\left\langle A(y(t),z(t))-A(y_{1}^{\mbox{target}},\widehat{z}),y-y_{1}^{%
\mbox{target}}\right\rangle _{V^{\ast },V}, \\
&=&D_{1}\left\Vert \nabla (y(t)-y_{1}^{\mbox{target}})\right\Vert
_{L^{2}(\Omega )}+\int_{\Omega }f(y(t),z(t))-f(y_{1}^{\mbox{target}},%
\widehat{z}))(y(t)-y_{1}^{\mbox{target}})dx \\
&\geq &D_{1}\left\Vert \nabla (y(t)-y_{1}^{\mbox{target}})\right\Vert
_{L^{2}(\Omega )}-L_{f}\left\Vert y(t)-y_{1}^{\mbox{target}}\right\Vert
_{L^{2}(\Omega )}^{2},
\end{eqnarray*}%
which is true for any $t\geq 0,$ in particular for $t\in (0,T_{\ast }+\delta
).$

Finally, we prove that $Py(t)=y_{1}^{\mbox{target}}$ for $t\geq T^{\ast }.$
Let us denote the solution to (\ref{ex-2}) for $t\geq T^{\ast }$ by $(%
\widetilde{y}(t),\widetilde{z}(t)).$ Then, it satisfies%
\begin{eqnarray}
\widetilde{y}_{t}-D_{1}\Delta \widetilde{y}+f(\widetilde{y},\widetilde{z})
&=&u\mbox{, \ \ \ \ \ \ \ \ in }(T^{\ast },\infty )\times \Omega ,
\label{1007} \\
\widetilde{z}_{t}-D_{2}\Delta \widetilde{z}+g(\widetilde{y},\widetilde{z})
&=&0,\mbox{ \ \ \ \ \ \ \ \ in }(T^{\ast },\infty )\times \Omega ,  \notag
\end{eqnarray}%
with homogeneous Neumann boundary conditions and the initial data at $%
t=T^{\ast },$ $\widetilde{y}(T^{\ast })=y_{1}^{\mbox{target}},$ $\widetilde{z%
}(T^{\ast })=z(T^{\ast })$. The second equation with $\widetilde{y}=y_{1}^{%
\mbox{target }}$has a unique solution well defined. If $u$ is replaced in (%
\ref{1007}) by%
\begin{equation*}
u(t)=\left\{ 
\begin{array}{l}
u^{\ast }(t),\mbox{ }t\leq T^{\ast } \\ 
\widetilde{u}(t),\mbox{ \ }t>T^{\ast }%
\end{array}%
\right.
\end{equation*}%
where $\widetilde{u}(t)=-D_{1}\Delta y_{1}^{\mbox{target}}+f(y_{1}^{%
\mbox{target}},\widetilde{z}(t)),$ then $(\widetilde{y}(t),\widetilde{z}%
(t))=(y_{1}^{\mbox{target}},\widetilde{z}(t))$ verifies the first equation
and this proves that the solution slides on the manifold $y_{1}^{%
\mbox{target}}$ for all $t\geq 0.$ Thus Theorem 5.5 can be applied to obtain
the conclusion of the Proposition 6.3. $\hfill \square $

\medskip

\noindent \textit{Case II.} We can also put into evidence a particular case
in which the choice of $\widehat{z}$ is independent of $u,$ $T$ and the
system solution. Let us assume that $y_{1}^{\mbox{target}}=0,$ and 
\begin{equation}
f(0,z)=0\mbox{ for all }z,\mbox{ }f(y,z)y\geq 0\mbox{ for all }(y,z)\in H.
\label{rd2}
\end{equation}%
We have to check (\ref{g5-000}). We set $\widehat{y}=(0,\widehat{z}),$ where 
$\widehat{z}$ in this case can be taken any value such that $A(0,\widehat{z}%
)\in V,$ in particular $\widehat{z}=0.$ We have 
\begin{equation*}
\left\langle A(y,z)-A(0,\widehat{z}),y\right\rangle _{V^{\ast
},V}=(-D_{1}\Delta y+f(y,z),y))_{H}\geq D_{1}\left\Vert \nabla y\right\Vert
_{H}^{2}\geq 0.
\end{equation*}%
A particular situation is $f(y,z)=yf_{2}(z),$ with $f_{2}(z)$ Lipschitz and
positive, for example $f_{2}(z)=\frac{z^{2}}{1+z^{2}}.$

\medskip

\noindent \textit{Case III. Reaction-diffusion systems with linear
perturbations. }Let us consider (\ref{ex-2}) with $f(y,z)=a_{1}y+b_{1}z$ and 
$g(y,z)=a_{2}y+b_{2}z.$ The functional framework is the same as in the
precedent example and all hypotheses are satisfied. We shall check only $%
(d_{4}),$ by setting $\widehat{y}=(y_{1}^{\mbox{target}},z),$ where $z$ is
the solution to (\ref{ex-2}) corresponding to $T$ and $u.$ We have%
\begin{eqnarray*}
&&\left\langle A(y,z)-A(\overline{y},\widehat{z}),y-\overline{y}%
\right\rangle _{V^{\ast },V}, \\
&=&(-D_{1}\Delta (y-\overline{y})+a_{1}(y-\overline{y})+b_{1}(z-z),(y-%
\overline{y})))_{H}\geq C_{1}\left\Vert (y-\overline{y})\right\Vert
_{V}^{2}\geq 0.
\end{eqnarray*}

\medskip

\noindent \textit{Case IV. FitzHugh-Nagumo reaction-diffusion model. }For $%
f(r_{1},r_{2})=\alpha _{0}r_{1}+r_{2},$ $g(r_{1},r_{2})=-\sigma r_{1}+\gamma
r_{2}$ and $D_{2}=0,$ the system (\ref{ex-2}) becomes the well-known
FitzHugh-Nagumo model (studied e.g. in \cite{Kunish-Wang-2012}). In this
case, the hypotheses are verified with the choice 
\begin{eqnarray*}
H &=&L^{2}(\Omega )\times L^{2}(\Omega ),\mbox{ }V=H^{1}(\Omega )\times
L^{2}(\Omega ),\mbox{ }V^{\ast }=(H^{1}(\Omega ))^{\ast }\times L^{2}(\Omega
),\mbox{ } \\
D_{H} &=&\left\{ y\in H^{2}(\Omega );\mbox{ }\frac{\partial y}{\partial \nu }%
=0\mbox{ on }\partial \Omega \right\} \times L^{2}(\Omega ).
\end{eqnarray*}

\paragraph{Example 4. Phase field systems.}

Let us consider the phase-field system of Caginalp type, for the phase
function $\varphi $ and the energy $\sigma $ written in the following form
(see e.g., \cite{BCGMR-sli}) 
\begin{equation*}
\sigma _{t}-k\Delta \sigma +kl\Delta \varphi =f+u,\mbox{ in }(0,\infty
)\times \Omega ,
\end{equation*}%
\begin{equation*}
\varphi _{t}-\nu \Delta \varphi +\beta (\varphi )+\pi (\varphi )=\gamma
\sigma -\gamma l\varphi ,\mbox{ in }(0,\infty )\times \Omega ,
\end{equation*}%
\begin{equation*}
\frac{\partial \varphi }{\partial \nu }=\frac{\partial \sigma }{\partial \nu 
}=0,\mbox{ in }(0,\infty )\times \partial \Omega ,
\end{equation*}%
\begin{equation*}
\varphi (0)=\varphi _{0},\mbox{ }\sigma (0)=\sigma _{0},\mbox{ in }\Omega
\end{equation*}%
and study $(\mathcal{P})$ with $U=(L^{2}(\Omega ),L^{2}(\Omega )),$ $%
Pu=(u,0),$ $B(u,u_{2})=(u,0)$ and the rest of the spaces as in Example 3.
Here $\beta (r)=r^{3}$ and $\pi (r)=-r$ for $r\in \mathbb{R}$, the function $%
\beta +\pi $ representing the double-well potential. The controllability
hypothesis $(c_{1})$ can be proved in a similar way with the proof of \cite%
{BCGMR-sli} for the case $u\in L^{\infty }(0,T;L^{2}(\Omega )),$ $\left\Vert
u(t)\right\Vert _{2}\leq \rho $ a.e. $t.$ The regularity of the second state
component $\varphi $ is proved as in the precedent example and so $\widehat{z%
}:=\varphi (t)$ which is the appropriate choice for checking (\ref{Ay^}) and
(\ref{g5-000}). We mention that the proof of the minimum time for the
Caginalp system with the singular logarithmic potential is considered in 
\cite{CGM-min-time}\textit{.}

\paragraph{Example 5. Diffusion with nonlocal controllers.}

We note that the theory works too if\textbf{\ }$B$ is a nonlocal operator.
This is the case when $B_{i}u_{i}\neq u_{i}.$ Let us consider, for instance
Example 1, where $B:L^{2}(\Omega _{1})\rightarrow L^{2}(\Omega )$ is defined
by 
\begin{equation*}
(Bu)(z)=\int_{\Omega _{1}}K(x,z)u(z)dz,\mbox{ }x\in \Omega 
\end{equation*}%
and $K\in L^{2}(\Omega \times \Omega _{1}).$ If the kernel $K$ is such that%
\begin{equation*}
\left\Vert B^{\ast }v\right\Vert _{L^{2}(\Omega _{1})}\geq \gamma \left\Vert
v\right\Vert _{L^{2}(\Omega )},\mbox{ for }v\in L^{2}(\Omega ),
\end{equation*}%
it follows that the controllability assumption holds and all conditions are
satisfied. Here, $B^{\ast }:L^{2}(\Omega )\rightarrow L^{2}(\Omega _{1})$ is
defined by $B^{\ast }v(x)=\int_{\Omega }K(x,z)v(x)dx,$ because%
\begin{eqnarray*}
(Bu,v)_{L^{2}(\Omega )} &=&\int_{\Omega }v(x)\int_{\Omega _{1}}K(x,z)u(z)dzdx
\\
&=&\int_{\Omega _{1}}u(z)\int_{\Omega }K(x,z)v(x)dxdz=(u,B^{\ast
}v)_{L^{2}(\Omega _{1})}.
\end{eqnarray*}

\section{Appendix}

\setcounter{equation}{0}

\paragraph{Some definitions and results related to operators in Hilbert
spaces.}

Let $H,$ $V$ be Hilbert spaces, $V^{\ast }$ the dual of $V$, $V\subset
H\subset V^{\ast }$ with compact injections. Let $A:V\rightarrow V^{\ast }$.

The operator $A$ is \textit{demicontinuous} if $y_{n}\rightarrow y$ strongly
in $V$ implies $Ay_{n}\rightarrow Ay$ weakly in $V^{\ast },$ as $%
n\rightarrow \infty .$

Let $\left\langle \cdot ,\cdot \right\rangle _{V^{\ast },V}$ denote the
pairing between $V^{\ast }$ and $V.$ The operator $A$ is \textit{coercive}
if 
\begin{equation*}
\lim_{\left\Vert y\right\Vert _{V}\rightarrow \infty }\frac{\left\langle
Ay,y-y^{0}\right\rangle _{V^{\ast },V}}{\left\Vert y\right\Vert _{V}}%
=+\infty ,\mbox{ for some }y^{0}\in V.
\end{equation*}

Let $A$ be an operator on the Hilbert space $H.$ It is called $m$-\textit{%
accretive} if it is \textit{accretive}, 
\begin{equation*}
\left( Ay-A\overline{y},y-\overline{y}\right) _{H}\geq 0,\mbox{ for all }y,%
\overline{y}\in D(A)
\end{equation*}%
and $m$-accretive if $R(I+A)=H,$ where $I$ is the identity operator and $R$
is the range. The operator $A$ is \textit{quasi }$m$\textit{-accretive} if $%
(\lambda I+A)$ is $m$-accretive for $\lambda $ sufficiently large.

The operator $A_{H}:D(A_{H})\subset H\rightarrow H$ is the \textit{%
restriction} of $A$ on $H$ defined as $A_{H}y=Ay$ for $y\in D(A_{H})=\{y\in
V;$ $Ay\in H\}.$

Let $A:V\rightarrow V^{\ast }$ be single-valued, monotone, demicontinuous
and coercive. Then, it follows that it is\textit{\ surjective} (see e.g. 
\cite{vb-springer-2010}, p. 36, Corollary 2.2) and $A_{H}$ is $m$-accretive
on $H\times H$.

Let $A\in C^{1}(V,V^{\ast }).$ The G\^{a}teaux derivative of $A$ is the
linear operator $A^{\prime }(y):V\rightarrow V^{\ast }$ defined by 
\begin{equation}
A^{\prime }(y)z=\lim_{\lambda \rightarrow 0}\frac{A(y+\lambda z)-Ay}{\lambda 
}\mbox{ strongly in }V^{\ast },\mbox{ for all }y,z\in V.  \label{g0}
\end{equation}

\noindent If $A_{H}\in C^{1}(D(A_{H}),H)$ we similarly define $%
(A_{H})^{\prime }(y):D(A_{H})^{\prime }=D(A_{H})\subset H\rightarrow H$ by 
\begin{equation}
(A_{H})^{\prime }(y)z=\lim_{\lambda \rightarrow 0}\frac{A_{H}(y+\lambda
z)-A_{H}y}{\lambda }\mbox{ strongly in }H,\mbox{ for all }y,z\in D(A_{H})
\label{g-00}
\end{equation}%
and observe that $(A_{H})^{\prime }(y)=(A^{\prime })_{H}(y),$ for $y\in
D(A_{H}),$ so that in the paper we use the notation $A_{H}^{\prime }(y).$

Let $X$ and $Y$ be Banach spaces. The operator $G:X\rightarrow Y$ is said to
be \textit{strongly continuous} from $X$ to $L_{s}(X,Y)$ if for $%
y_{n}\rightarrow y$ strongly in $X,$ as $n\rightarrow \infty $ it follows 
\begin{equation}
\left\Vert G(y_{n})\psi -G(y)\psi \right\Vert _{Y}\rightarrow 0%
\mbox{ for
all }\psi \in X.  \label{g-000}
\end{equation}

\paragraph{Duality mapping.}

Let $U$ be Banach spaces with the dual $U^{\ast }$ uniformly convex,
implying that $U^{\ast }$ and $U$ is reflexive (see e.g., \cite%
{vb-springer-2010}, p. 2). Also, it follows that the norm in $U$ is G\^{a}%
teaux differentiable.

Let $F:U\rightarrow U^{\ast }$ be the duality mapping of $U$, which is
single valued and continuous (see e.g., \cite{vb-springer-2010}, p. 2,
Theorem 1.2). We recall that 
\begin{equation}
\left\langle Fu,u\right\rangle _{U^{\ast },U}=\left\Vert u\right\Vert
_{U}^{2},\mbox{ }\left\Vert Fu\right\Vert _{U^{\ast }}=\left\Vert
u\right\Vert _{U}.  \label{A1}
\end{equation}%
Let $K=\{u\in U;$ $\left\Vert u\right\Vert _{U}\leq \rho \},$ let $I_{K}$ be
the indicator function of $K$ and define 
\begin{equation}
j:U\rightarrow \mathbb{R},\mbox{ }j(u)=I_{K}(u).  \label{A6-0}
\end{equation}%
Then, 
\begin{equation}
\partial j(u)=\partial I_{K}(u)=N_{K}(u)=\left\{ 
\begin{array}{l}
\{\lambda Fu;\lambda >0\}\mbox{ \ \ \ }\left\Vert u\right\Vert _{U}=\rho \\ 
0,\mbox{ \ \ \ \ \ \ \ \ \ \ \ \ \ \ \ \ \ }\left\Vert u\right\Vert _{U}<\rho
\\ 
\varnothing ,\mbox{ \ \ \ \ \ \ \ \ \ \ \ \ \ \ \ \ }\left\Vert u\right\Vert
_{U}>\rho%
\end{array}%
\right.  \label{A5}
\end{equation}%
where $\partial j:U\rightarrow U^{\ast }$ is the subdifferential of $j,$ $%
N_{K}$ is the normal cone to $K$ in $U^{\ast }$ and $\lambda >0.$ The first
line in (\ref{A5}) should be understood in the multivalued sense.

The conjugate of $j$ is $j^{\ast }:U^{\ast }\rightarrow \mathbb{R}$,%
\begin{equation}
j^{\ast }(z)=\sup_{u\in K}\left\{ \left\langle z,v\right\rangle _{U^{\ast
},U}-j(v)\right\} =\sup \left\{ \left\langle z,v\right\rangle _{U^{\ast },U};%
\mbox{ }\left\Vert v\right\Vert _{U}\leq \rho \right\} =\rho \left\Vert
z\right\Vert _{U^{\ast }}.  \label{A2}
\end{equation}%
Then, 
\begin{equation}
\partial j^{\ast }(z)=(\partial j)^{-1}(z)=N_{K}^{-1}(z)=\rho \frac{F^{-1}z}{%
\left\Vert z\right\Vert _{U^{\ast }}},\mbox{ }z\in U^{\ast }.  \label{A3}
\end{equation}%
Since $U$ is reflexive, $F^{-1}$ is just the duality mapping of $U^{\ast }$
and so $D(F^{-1})=U^{\ast }$.

\noindent If $z\in N_{K}(u),$ then $u\in N_{K}^{-1}(z)$ and we have 
\begin{equation}
\left\langle z,u\right\rangle _{U^{\ast },U}=\rho \left\Vert z\right\Vert
_{U^{\ast }}.  \label{A4-0}
\end{equation}%
Let $\varepsilon $ be positive and define%
\begin{equation}
j_{\varepsilon }(u)=\frac{\varepsilon }{2}\left\Vert u\right\Vert
_{U}^{2}+I_{K}(u).  \label{A4-1}
\end{equation}%
We recall that the subdifferential 
\begin{equation}
\partial \left( \frac{1}{2}\left\Vert u\right\Vert _{U}^{2}\right) =Fu,
\label{A4-1-0}
\end{equation}%
whence%
\begin{equation}
\partial j_{\varepsilon }(u)=\varepsilon Fu+N_{K}(u),\mbox{ for all }u\in K.
\label{A4-2}
\end{equation}%
Then,%
\begin{eqnarray*}
j_{\varepsilon }^{\ast }(\zeta ) &=&\sup_{v\in K}\left\{ \left\langle \zeta
,v\right\rangle _{U^{\ast },U}-j_{\varepsilon }(v)\right\} =-\inf_{v\in
K}\left\{ I_{K}(v)+\frac{\varepsilon }{2}\left\Vert v\right\Vert
_{U}^{2}-\left\langle \zeta ,v\right\rangle _{U^{\ast },U}\right\} \\
&=&\left\langle \varepsilon Fv_{\varepsilon }+z_{\varepsilon
},v_{\varepsilon }\right\rangle _{U^{\ast },U}-\frac{\varepsilon }{2}%
\left\Vert v_{\varepsilon }\right\Vert _{U}^{2}=\frac{\varepsilon }{2}%
\left\Vert v_{\varepsilon }\right\Vert _{U}^{2}+\left\langle z_{\varepsilon
},v_{\varepsilon }\right\rangle _{U^{\ast },U}.
\end{eqnarray*}%
We specify that in the lines before the infimum is realized at $%
v_{\varepsilon }$ which is the solution to the equation $\varepsilon
Fv_{\varepsilon }+N_{K}(v_{\varepsilon })\ni \zeta ,$ that is $\varepsilon
Fv_{\varepsilon }+z_{\varepsilon }=\zeta ,$ where $z_{\varepsilon }\in
N_{K}(v_{\varepsilon }).$ Therefore, 
\begin{equation*}
j_{\varepsilon }^{\ast }(\zeta )=\left\langle \varepsilon Fv_{\varepsilon
}+z_{\varepsilon },v_{\varepsilon }\right\rangle _{U^{\ast },U}-\frac{%
\varepsilon }{2}\left\Vert v_{\varepsilon }\right\Vert _{U}^{2}=\frac{%
\varepsilon }{2}\left\Vert v_{\varepsilon }\right\Vert _{U}^{2}+\left\langle
z_{\varepsilon },v_{\varepsilon }\right\rangle _{U^{\ast },U},
\end{equation*}%
implying by (\ref{A4-0}) that 
\begin{equation*}
j_{\varepsilon }^{\ast }(\zeta )=\frac{\varepsilon }{2}\left\Vert
v_{\varepsilon }\right\Vert _{U}^{2}+\rho \left\Vert z_{\varepsilon
}\right\Vert _{U^{\ast }},\mbox{ }z_{\varepsilon }\in N_{K}(v_{\varepsilon
}).
\end{equation*}%
Finally, if $u\in (\varepsilon F+N_{K})^{-1}(\zeta ),$ it follows that 
\begin{equation}
(\partial j_{\varepsilon })^{-1}(\zeta )=(\varepsilon F+N_{K})^{-1}(\zeta
)=\partial j_{\varepsilon }^{\ast }(\zeta )  \label{A4-4}
\end{equation}%
and 
\begin{equation}
j_{\varepsilon }^{\ast }(\zeta )=\frac{\varepsilon }{2}\left\Vert
u_{\varepsilon }\right\Vert _{U}^{2}+\rho \left\Vert z_{\varepsilon
}\right\Vert _{U^{\ast }},\mbox{ }z_{\varepsilon }\in N_{K}(u_{\varepsilon
}).  \label{A4-3}
\end{equation}%
\medskip

\paragraph{The canonical isomorphism}

Assume now that $H$ and $V$ are Hilbert spaces, and $V$ has the dual $%
V^{\ast }.$ The \textit{duality mapping}, which we denote by $\Gamma
:V\rightarrow V^{\ast }$ is the \textit{canonical isomorphism} of $V$ onto $%
V^{\ast }$ (see e.g., \cite{vb-springer-2010}, p. 1). We have 
\begin{equation}
\Gamma \in L(V,V^{\ast }),\mbox{ }\left\langle \Gamma v,v\right\rangle
_{V^{\ast },V}=\left\Vert v\right\Vert _{V}^{2},\mbox{ }\left\Vert \Gamma
v\right\Vert _{V^{\ast }}=\left\Vert v\right\Vert _{V}.  \label{Gamma1}
\end{equation}%
In addition, $\Gamma _{H},$ the \textit{restriction} of $\Gamma $ to $H,$ is 
$m$-accretive on $H\times H,$ with the \textit{linear domain denoted} $D_{H}$
which is densely, continuously and compactly embedded in $V,$ 
\begin{equation}
D(\Gamma _{H}):=D_{H}\subset V.  \label{DH}
\end{equation}%
For $\nu >0,$ we denote by $\Gamma _{\nu }$ the Yosida approximation of $%
\Gamma _{H},$ that is 
\begin{equation}
\Gamma _{\nu }y=\frac{1}{\nu }(I-(I+\nu \Gamma _{H})^{-1})y=\Gamma
_{H}(I+\nu \Gamma _{H})^{-1}y,\mbox{ }y\in H.  \label{Yosida}
\end{equation}

\paragraph{Comments on the hypothesis of controllability.}

Hypothesis $(c_{1})$ ensures that the admissible set for problem $(\mathcal{P%
})$ is not empty$.$ For example, in the case of Caginalp phase field models
the proof of the controllability was provided in \cite{BCGMR-sli}. Further,
we shall argue for the reliability of such an hypothesis, giving a brief
proof of the controllability of (\ref{g4})-(\ref{g4-1}) in some cases.
First, let us set 
\begin{equation*}
u(t)=-\rho \,\mbox{Sign\thinspace }(B^{\ast }P(y(t)-y^{tar})),
\end{equation*}%
where Sign$:U^{\ast }\rightarrow 2^{U^{\ast }}$ is defined by 
\begin{equation}
\mbox{Sign\thinspace }v=\left\{ 
\begin{array}{c}
\frac{v}{\left\Vert v\right\Vert _{U^{\ast }}},\mbox{ \ }y\neq 0 \\ 
B(0,\rho ),\mbox{ }y=0.%
\end{array}%
\right.  \label{g17}
\end{equation}%
Here, $B(0,\rho )$ is the ball of center 0 and radius $\rho $ in $U^{\ast }.$
It is well known that $v\rightarrow $Sign\thinspace $v$ is $m$-accretive on $%
U^{\ast }$.

Let us consider the problem 
\begin{eqnarray}
y^{\prime }(t)+Ay(t)) &\ni &-\rho B\mbox{Sign\thinspace }(B^{\ast
}P(y(t)-y^{tar})),\mbox{ a.e. }t\in (0,T)  \label{g17-1} \\
y(0) &=&y_{0}.  \notag
\end{eqnarray}%
We refer to the case when one state component is controlled by one
controller, that is 
\begin{equation*}
Py=(y_{1},0),\mbox{ }U=U_{1}\times U_{2},\mbox{ }U_{1}=U_{1}^{\ast },\mbox{ }%
B=(B_{1},0),\mbox{ }B_{1}:U_{1}\rightarrow H_{1},
\end{equation*}%
and assume 
\begin{equation}
R(B_{1})=H_{1}.  \label{g8-0}
\end{equation}%
(When $P=I,$ $Bu=(B_{1}u_{1},B_{2}u_{2}),$ we impose the condition $R(B)=H.$
The proof is the same, by replacing $H_{1}$ by $H.)$

Hypothesis (\ref{g8-0}) implies, by the Banach closed range theorem (see 
\cite{Yosida}, p. 208, Corollary 1) that $(B_{1}^{\ast })^{-1}$ is
continuous from $U_{1}$ to $H_{1}$ and $B_{1}^{-1}\in L(H_{1},U_{1}).$ This
means $\left\Vert (B_{1}^{\ast })^{-1}z\right\Vert _{H_{1}}\leq C\left\Vert
z\right\Vert _{U_{1}}$ for $z\in U_{1},$ or, equivalently 
\begin{equation}
\left\Vert w\right\Vert _{H_{1}}\leq C\left\Vert B_{1}^{\ast }w\right\Vert
_{U_{1},}\mbox{ for }w\in H_{1}.  \label{g12-2}
\end{equation}%
We also assume that 
\begin{equation}
(A_{H}y-A_{H}y^{tar},P(y-y^{tar}))_{H}\geq -C_{1}\left\Vert
P(y-y^{tar})\right\Vert _{H}^{2},\mbox{ for all }y\in D_{H},\mbox{ }%
Py^{tar}\in P(D_{H}).  \label{g12-10}
\end{equation}

\noindent It is clear that when $B_{1}=I,$ then $U_{1}=H_{1}.$ Otherwise, we
have the situation in Example 5.

\medskip

\noindent \textbf{Proposition 7.1.} \textit{Let }$y_{0},$ $Py^{tar}\in
P(D_{H})$ \textit{and let }$(a_{1}),$\textit{\ }$(a_{2}),$ $(b_{1}),$ (\ref%
{g8-0}),\textit{\ }(\ref{g12-10})\textit{\ and }%
\begin{equation*}
\rho >\left\Vert A_{H}y^{tar}\right\Vert _{H}+C_{1}\left\Vert
P(y_{0}-y^{tar})\right\Vert _{H}
\end{equation*}%
\textit{hold.} \textit{Then, there exists }$T_{\ast }\in (0,T)$\textit{\
such that, for }$\rho $\textit{\ large enough, }$Py(T_{\ast })=Py^{tar},$ 
\textit{where }$y$\textit{\ is the solution to} (\ref{g17-1})$.$

\medskip

\noindent \textbf{Proof.} The operator $B$Sign\thinspace $(B^{\ast }Pv)$ is $%
m$-accretive on $H_{1}.$ Indeed, for $v,\bar{v}\in H_{1}$ we have 
\begin{equation*}
(B\mbox{Sign}\,(B^{\ast }Pv)-B\mbox{Sign}\,(B^{\ast }P\overline{v}),v-%
\overline{v})_{H_{1}}=(\mbox{Sign}\,(B^{\ast }Pv)-\mbox{Sign}\,(B^{\ast }P%
\overline{v}),B^{\ast }(v-\overline{v}))_{H_{1}}\geq 0,
\end{equation*}%
because $B^{\ast }Pv=B^{\ast }v,$ $w\rightarrow $Sign\thinspace $w$ is $m$%
-accretive and $P^{2}=P$. For the $m$-accretivity let us consider the
equation 
\begin{equation}
y+\rho B\mbox{Sign\thinspace }(B^{\ast }P(y-y^{tar}))=f\in H_{1}
\label{g12-3}
\end{equation}%
which, by denoting $z=y-y^{tar},$ is equivalent with $z+\rho B$%
Sign\thinspace $(B^{\ast }z)=f-y^{tar}.$ We set $B^{\ast }z=v\in U_{1}$ and
get 
\begin{equation}
B^{-1}(B^{\ast })^{-1}v+\rho \mbox{Sign\thinspace }v=B^{-1}f_{1}\in U_{1}.
\label{g12-4}
\end{equation}%
Denoting $G=B^{-1}(B^{\ast })^{-1}$ we see that $G\in L(U_{1},U_{1})$ and $%
(Gv,v)_{U_{1}}=\left\Vert (B^{\ast })^{-1}v\right\Vert _{H_{1}}^{2}\geq
\left\Vert v\right\Vert _{U_{1}}^{2}$ (because $B^{\ast }$ is continuous,
see (\ref{g12-2})$).$ Now, Sign\thinspace $w$ is $m$-accretive on $%
U_{1}\times U_{1}$, hence $R(G+$Sign$)=U_{1}$ (see e.g., \cite%
{vb-springer-2010}, p. 44, Corollary 2.6).

Then, we prove that (\ref{g17-1}) has a unique solution.

Since $A_{H}$ is quasi $m$-accretive and $S=\rho B$Sign\thinspace $(B^{\ast
}P(y(t)-y^{tar}))$ is $m$-accretive with $D(S)=H_{1}$ and $D_{H}\cap \overset%
{\circ }{D(S)}=D_{H}\neq \varnothing ,$ it follows by that $A_{H}+S$ is
quasi $m$-accretive on $H\times H$ (see \cite{vb-springer-2010}, p. 43,
Theorem 2.6)$.$ Therefore, (\ref{g17-1}) has a unique solution $y\in
L^{\infty }(0,T;D_{H})\cap W^{1,\infty }(0,T;H)$ satisfying estimate (\ref%
{g15}) (see the proof of Theorem 3.2).

\noindent Now, we can justify the controllability assertion. Let us write (%
\ref{g17-1}) in the equivalent form 
\begin{equation*}
(y-y^{tar})_{t}+A_{H}y-A_{H}y^{tar}+\rho B\mbox{Sign\thinspace }(B^{\ast
}P(y-y^{tar}))\ni -A_{H}y^{tar}
\end{equation*}%
and multiply it by $P(y(t)-y^{tar}).$ By (\ref{g12-10}) we get 
\begin{eqnarray*}
&&\frac{1}{2}\frac{d}{dt}\left\Vert P(y(t)-y^{tar})\right\Vert
_{H}^{2}-C_{1}\left\Vert P(y(t)-y^{tar})\right\Vert _{H}^{2}+\rho \left( B%
\mbox{Sign\thinspace }(B^{\ast }P(y(t)-y^{tar})),P(y(t)-y^{tar}\right) )_{H}
\\
&\leq &\left\Vert A_{H}y^{tar}\right\Vert _{H}\left\Vert
P(y(t)-y^{tar})\right\Vert _{H}.
\end{eqnarray*}%
Further, we obtain 
\begin{eqnarray*}
&&\left\Vert P(y(t)-y^{tar})\right\Vert _{H}\frac{d}{dt}\left\Vert
P(y(t)-y^{tar})\right\Vert _{H}+\rho \left\Vert B^{\ast
}P(y(t)-y^{tar})\right\Vert _{U_{1}^{\ast }} \\
&\leq &\left\Vert A_{H}y^{tar}\right\Vert _{H}\left\Vert
P(y(t)-y^{tar})\right\Vert _{H}+C_{1}\left\Vert P(y(t)-y^{tar})\right\Vert
_{H}^{2}.
\end{eqnarray*}%
Next, we use (\ref{g12-2}) for $w=P(y(t)-y^{tar})$ which implies $\left\Vert
B^{\ast }P(y(t)-y^{tar})\right\Vert _{U_{1}}\geq C\left\Vert
P(y(t)-y^{tar})\right\Vert _{H_{1}}$ and so 
\begin{equation*}
\frac{d}{dt}\left\Vert P(y(t)-y^{tar})\right\Vert _{H}-C_{1}\left\Vert
P(y(t)-y^{tar})\right\Vert _{H}+\rho \leq \left\Vert A_{H}y^{tar}\right\Vert
_{H}.
\end{equation*}%
For $\rho >\left\Vert A_{H}y^{tar}\right\Vert _{H},$ this yields%
\begin{equation*}
\left\Vert P(y(t)-y^{tar})\right\Vert _{H}<e^{C_{1}t}\left\Vert
P(y_{0}-y^{tar})\right\Vert _{H}-\frac{(\rho -\left\Vert
A_{H}y^{tar}\right\Vert _{H})}{C_{1}}(e^{C_{1}t}-1).
\end{equation*}%
Finally, we obtain that $t\rightarrow \left\Vert P(y(t)-y^{tar})\right\Vert
_{H}$ is strictly decreasing, vanishes at $t=T_{\ast }$ below and the
previous relation takes place for 
\begin{equation*}
t\leq T_{\ast }=\ln \frac{\rho -\left\Vert A_{H}y^{tar}\right\Vert _{H}}{%
\rho -(\left\Vert A_{H}y^{tar}\right\Vert _{H}+C_{1}\left\Vert
P(y_{0}-y^{tar})\right\Vert _{H})},
\end{equation*}%
for $\rho >\left\Vert A_{H}y^{tar}\right\Vert _{H}+C_{1}\left\Vert
P(y_{0}-y^{tar})\right\Vert _{H}.$ We also observe that $T^{\ast }$
decreases as $\rho $ increases and in fact $T_{\ast }\rightarrow 0$ as $\rho
\rightarrow \infty .$ This end the proof.\hfill $\square $

\medskip

\noindent \textbf{Acknowledgement. }The present paper benefits from the
support of the of Ministry of Research and Innovation, CNCS --UEFISCDI,
project number PN-III-P4-ID-PCE-2016-0011.

\end{document}